%% file: Lagrangian_PSS_and_comparison.tex
\documentclass[11pt,intlimits,reqno]{amsart}

\usepackage[hypertex]{hyperref}
\usepackage{graphicx}
\usepackage{color}
\usepackage{amsfonts,amssymb}
\usepackage{amsthm}
\usepackage{txfonts}
\usepackage[all]{xy}

\usepackage{mathrsfs}

\usepackage{a4wide}

\newtheorem{neu}{}[section]
\newtheorem{Cor}[neu]{Corollary}
\newtheorem*{Cor*}{Corollary}
\newtheorem{Thm}[neu]{Theorem}
\newtheorem*{Thm*}{Theorem}
\theoremstyle{definition}
\newtheorem{Prop}[neu]{Proposition}
\newtheorem*{Prop*}{Proposition}
\newtheorem{Lemma}[neu]{Lemma}
\newtheorem*{Rmk*}{Remark}
\newtheorem{Rmk}[neu]{Remark}

\newtheorem*{Ex*}{Example}

\newtheorem{Claim}[neu]{Claim}
\newtheorem{Def}[neu]{Definition}
\newtheorem{Conv}[neu]{Convention}

\newcommand{\Z}{\mathbb{Z}}
\newcommand{\R}{\mathbb{R}}
\newcommand{\C}{\mathbb{C}}

\newcommand{\zw}{\hspace*{0.75em}}
\newcommand{\pf}{\longrightarrow}

\newcommand{\hpf}{\hookrightarrow}
\newcommand{\wrt}{with respect to }

\newcommand{\CZ}{\mu_{\mathrm{CZ}}}
\newcommand{\Mas}{\mu_{\mathrm{Maslov}}}
\newcommand{\Morse}{\mu_{\mathrm{Morse}}}

\newcommand{\id}{\mathrm{id}}
\newcommand{\ind}{\mathrm{ind\,}}

\newcommand{\om}{\omega}
\newcommand{\ev}{\mathrm{ev}}

\newcommand{\Poincare}{Poincar\'{e} }

\renewcommand{\O}{\mathcal{O}}
\newcommand{\A}{\mathcal{A}}

\renewcommand{\P}{\mathcal{P}}

\newcommand{\F}{\mathcal{F}}

\newcommand{\D}{\mathbb{D}}

\newcommand{\M}{\mathcal{M}}
\newcommand{\Mh}{\widehat{\mathcal{M}}}
\newcommand{\J}{\mathcal{J}}
\newcommand{\Jh}{\widehat{\mathcal{J}}}
\newcommand{\B}{\mathcal{B}}

\renewcommand{\H}{\mathrm{H}}

\newcommand{\PSS}{\mathrm{PSS}}
\newcommand{\Ham}{\mathrm{Ham}}

\newcommand{\PD}{\mathrm{PD}}
\newcommand{\CF}{\mathrm{CF}}

\newcommand{\HF}{\mathrm{HF}}

\newcommand{\CM}{\mathrm{CM}}

\newcommand{\Crit}{\mathrm{Crit}}

\newcommand{\beq}{\begin{equation}}
\newcommand{\beqn}{\begin{equation}\nonumber}
\newcommand{\eeq}{\end{equation}}

\newcommand{\bea}{\begin{equation}\begin{aligned}}
\newcommand{\bean}{\begin{equation}\begin{aligned}\nonumber}
\newcommand{\eea}{\end{aligned}\end{equation}}

\numberwithin{equation}{section}

\definecolor{Peter}{rgb}{1,0,0}

\definecolor{PeterB}{rgb}{0,0,1}

\begin{document}
\title[A Lagrangian PSS morphism and two comparison homomorphisms in Floer homology]
{A Lagrangian Piunikhin-Salamon-Schwarz morphism and\\[1ex] two comparison homomorphisms in Floer homology}
\author{Peter Albers}
\address{Peter Albers\\\;Courant Institute\\\;New York University\vspace*{.3ex}}
\email{\;\htmladdnormallink{albers@cims.nyu.edu}{mailto:albers@cims.nyu.edu}\vspace*{.3ex}}
\urladdr{\;\htmladdnormallink{http://www.cims.nyu.edu/$\sim$albers}{http://www.cims.nyu.edu/~albers}}
\date{September 2007}
\begin{abstract}
This article address two issues.
First, we explore to what extent the techniques of Piunikhin, Salamon and Schwarz in
\cite{Piunikhin_Salamon_Schwarz_Symplectic_Floer_Donaldson_theory_and_quantum_cohomology} can be carried over to Lagrangian Floer homology. In
\cite{Piunikhin_Salamon_Schwarz_Symplectic_Floer_Donaldson_theory_and_quantum_cohomology} an isomorphism between Hamiltonian
Floer homology and singular homology is established. In contrast, Lagrangian Floer homology is not isomorphic to the singular
homology of the Lagrangian submanifold, in general. Depending on the minimal Maslov number, we construct for certain degrees two homomorphisms between
Lagrangian Floer homology and singular homology. In degrees, where both maps are defined, we prove them to be isomorphisms.
Examples show that this statement is sharp.

Second, we construct two comparison homomorphisms between Lagrangian and Hamiltonian Floer homology. They are defined without degree restrictions and
are proven to be the natural analogs to the homomorphisms in singular homology induced by the inclusion map of the Lagrangian submanifold into
the ambient symplectic manifold.
\end{abstract}
\maketitle
\vspace*{8ex}
\tableofcontents

\newpage
\section{Main results}\label{section:main_results}
\begin{Thm}\label{thm:existence_of_Lagrangian_PSS}
We consider a $2n$-dimensional, closed, symplectic manifold $(M,\om)$ and a closed, monotone Lagrangian
submanifold $L\subset M$ of minimal Maslov number $N_L\geq2$. Then there exist homomorphisms
\begin{gather}
\varphi_k:\HF_k(L,\phi_H(L))\pf\H^{n-k}(L;\Z/2)\quad\text{for }k\leq N_L-2\;,\\[1ex]
\rho_k:\H^{n-k}(L;\Z/2)\pf\HF_k(L,\phi_H(L))\quad\text{for }k\geq n-N_L+2\;,
\end{gather}
where $H:S^1\times M\pf\R$ is a Hamiltonian function and $\phi_H$ the corresponding Hamiltonian diffeomorphism. 
For $n-N_L+2\leq k\leq N_L-2$ these maps are inverse to each other
\beq
\varphi_k\circ\rho_k=\id_{\H^{n-k}(L;\Z/2)}\quad\text{and}\quad
        \rho_k\circ\varphi_k=\id_{\HF_k(L,\phi_H(L))}\;.
\eeq
\end{Thm}\noindent
These homomorphisms are constructed using the ideas introduced by \emph{Piunikhin, Salamon} and \emph{Schwarz} in
\cite{Piunikhin_Salamon_Schwarz_Symplectic_Floer_Donaldson_theory_and_quantum_cohomology}. We call them
\textsc{Lagrangian PSS morphisms}.\\[1ex]
Let us put this theorem into context.
In the situation of theorem \ref{thm:existence_of_Lagrangian_PSS}
both \emph{Lagrangian Floer homology} $\HF_*(L,\phi_H(L))$ and \emph{Hamiltonian Floer homology} $\HF_*(H)$ are defined. Originally, Floer
\cite{Floer_Morse_theory_for_Lagrangian_intersections} defined his theories under more restrictive conditions, which were generalized later by himself
\cite{Floer_symplectic_fixed_points_and_holomorphic_spheres} and
Oh \cite{Oh_Floer_cohomology_of_Lagrangian_intersections_and_pseudo-holomorphic_disks_I} to the above assumption of monotonicity (see definition
\ref{def:minimal_Maslov_and_Chern_nb}).

Floer \cite{Floer_symplectic_fixed_points_and_holomorphic_spheres} proved that Hamiltonian Floer homology $\HF_*(H)$ is isomorphic to the singular
homology of the symplectic manifold $M$. The construction  of Hamiltonian Floer homology together with this isomorphism property was
generalized to semi-positive symplectic manifolds by Hofer and Salamon \cite{Hofer_Salamon_Floer_homology_and_Novikov_rings}. Another construction of
an isomorphism was given later by Piunikhin, Salamon and Schwarz in
\cite{Piunikhin_Salamon_Schwarz_Symplectic_Floer_Donaldson_theory_and_quantum_cohomology}.

In \cite{Floer_Morse_theory_for_Lagrangian_intersections}
Floer defined Lagrangian Floer homology under the assumption $\om|_{\pi_2(M,L)}=0$, in which case Lagrangian Floer homology in fact
is isomorphic to the singular homology of the Lagrangian submanifold $L$. But under the more general assumption of monotonicity of $L$ this in no
longer true. Lagrangian Floer homology
might even vanish completely, for example in case $L\cap\phi_H(L)=\emptyset$, that is, if $L$ is displaceable (cf.~definition \ref{def:displaceable}).

In \cite{Oh_Floer_cohomology_spectral_sequences_and_the_Maslov_class}, Oh constructs a spectral sequence relating Lagrangian Floer homology and
singular homology, in particular proving that Lagrangian Floer homology still is isomorphic to the singular homology of the Lagrangian submanifold $L$
given that the minimal Maslov number $N_L$ (see definition \ref{def:minimal_Maslov_and_Chern_nb}) satisfies $N_L\geq\dim L+2$.

Serious technical difficulties in the construction of Floer homology are caused by holomorphic spheres in $(M,\om)$ in the case of Hamiltonian
Floer homology and both, holomorphic spheres and holomorphic disks with boundary on $L$, in the case of Lagrangian Floer homology. In the monotone
case holomorphic spheres and disks can be handled by a precise understanding of the behavior of the Fredholm theory whenever holomorphic objects appear.
When extending the construction of Floer homology beyond the monotone case, Hofer and Salamon dealt with holomorphic spheres by proving very subtle
transversality results for moduli spaces of holomorphic spheres. This uses the dichotomy between multiply covered and somewhere injective holomorphic
spheres (see \cite{Hofer_Salamon_Floer_homology_and_Novikov_rings,McDuff_Salamon_J_holomorphic_curves_and_symplectic_topology}). Likewise, the results
in \cite{Piunikhin_Salamon_Schwarz_Symplectic_Floer_Donaldson_theory_and_quantum_cohomology} rely heavily on transversality results for holomorphic
spheres. The mentioned dichotomy fails in the Lagrangian case since there exist holomorphic disks which are neither multiply covered nor somewhere
injective.

The obvious modification of the idea of Piunikhin, Salamon and Schwarz to Lagrangian Floer homology is, in general, not well-defined due
to bubbling-off of holomorphic disks. In this article we handle this problem using arguments involving only the Fredholm
index. This leads to the degree restrictions of the maps $\varphi$ and $\rho$ (see theorem \ref{thm:existence_of_Lagrangian_PSS}).

Surprisingly, we still can prove that these maps are inverse to each other (in case they are defined simultaneously). For this we cannot employ the
methods of Piunikhin, Salamon and Schwarz due to the above mentioned problems with holomorphic disks. We use a tricky mix out of Fredholm index
arguments and
transversality methods to conclude. All proofs break down beyond certain obvious degrees restrictions. This failure is not just a technical shortcoming
but a necessity, as is illustrated by (counter-)examples (see remark below).

Theorem \ref{thm:existence_of_Lagrangian_PSS} leads to a number of applications, e.g.~new bounds on the Maslov index of displaceable Lagrangian
submanifolds and existence of holomorphic disks. It
should be mentioned that not only the isomorphisms statement but the mere existence of the Lagrangian PSS morphisms has applications
(see section \ref{section:applications}).

Both Floer theories are graded. The degree in Hamiltonian Floer homology is given by the Conley-Zehnder index $\CZ$ (see definition
\ref{def:minimal_Maslov_and_Chern_nb}). This gives rise to an integer grading modulo twice
the minimal Chern number $N_M$. The grading of Lagrangian Floer homology is only well-defined modulo the minimal Maslov
number $N_L$ and up to an overall shift (see definition \ref{def:minimal_Maslov_and_Chern_nb} and section \ref{section:recollection_of_Floer_homology}
for details).

\begin{Rmk}
The statement of theorem \ref{thm:existence_of_Lagrangian_PSS} is sharp. In section \ref{section:applications} we recall examples provided
by Polterovich \cite{Polterovich_Monotone_Lagrange_submanifolds_of_linear_spaces}, where theorem \ref{thm:existence_of_Lagrangian_PSS} fails
precisely beyond the range of degrees specified in theorem \ref{thm:existence_of_Lagrangian_PSS}. Furthermore, this provides
existence results for holomorphic disks, see section \ref{section:applications}.
\end{Rmk}

\begin{Rmk}
Singular homology is modeled via Morse homology, see \cite{Schwarz_Matthias_Morse_homology}.
Both maps $\varphi_k$ and $\rho_k$ are defined on chain level, that is on $\CF_k(L,\phi_H(L))$ and on $\CM^{n-k}(f)$ for some Morse
function $f:L\rightarrow\R$. Furthermore, they are natural with respect to the continuation homomorphisms induced by changing the Hamiltonian
and the Morse function, respectively. The same holds for the comparison homomorphisms $\chi_k$ and $\tau_k$ in theorem \ref{thm:diagram_commutes} below.
\end{Rmk}

\begin{Rmk}
We note that the maps $\varphi_k$ and $\rho_k$ in theorem \ref{thm:existence_of_Lagrangian_PSS} vanish for $k<0$ respectively $k>n$ by
construction. As mentioned above the grading
of Lagrangian Floer homology is modulo the minimal Maslov number $N_L$. Since there is exactly one $k$ modulo $N_L$ between $0$ and $N_L-2$ respectively
between $n-N_L+2$ and $n$, the statement of theorem \ref{thm:existence_of_Lagrangian_PSS} is unambiguous.
\end{Rmk}
\noindent
Even though Hamiltonian Floer homology if graded modulo the twice the minimal Chern number $N_M$ of $M$ we need to reduce this grading to modulo the
minimal Maslov number $N_L$ in the next statement. We note that $2N_M$ is a multiple of $N_L$.

\begin{Thm}\label{thm:diagram_commutes}
Under the same assumptions as in theorem \ref{thm:existence_of_Lagrangian_PSS} there exist (for all degrees $k$) two \textsc{comparison homomorphisms}
\begin{gather}
\chi_k:\HF_k(L,\phi_H(L))\pf\HF_{k-n}(H)\;,\\[1ex]
\tau_k:\HF_{k}(H)\pf\HF_{k}(L,\phi_H(L))\;,
\end{gather}
such that the following two diagrams commute whenever $\varphi_k$ or $\rho_k$ are defined
\beqn
\xymatrix{
{\HF_k(L,\phi_H(L))}\ar[d]^{\varphi_k}\ar[r]^-{\chi_k}&{\HF_{k-n}(H)}&
             &{\HF_{k}(H)}\ar[r]^-{\tau_k}&{\HF_{k}(L,\phi_H(L))}\ar[d]^{\varphi_k}\\
{\H^{n-k}(L;\Z/2)}\ar@<1ex>[u]^{\rho_k}\ar[r]^{\iota^!}&{\H^{2n-k}(M;\Z/2)}\ar[u]^{\PSS}_\cong&
            &{\H^{n-k}(M;\Z/2)}\ar[u]^{\PSS}_\cong\ar[r]^{\iota^*}&{\H^{n-k}(L;\Z/2)}\ar@<1ex>[u]^{\rho_k} \\
}\eeq
where the inclusion map $\iota:L\hpf M$ induces on singular homology the pull-back map $\iota^*$ and the shriek map
$\iota^!=\PD\circ\iota_*\circ\PD^{-1}$. Here $\PD$ denotes \Poincare duality.
\end{Thm}

The following remark about the homomorphism $\chi_k:\HF_k(L,\phi_H(L))\pf\HF_{k-n}(H)$ is in order.
The author constructed the map $\chi$ in summer 2005. Later, it became apparent that the
same construction was simultaneously used by A.~Abbondandolo and M.~Schwarz in a different context. In
\cite{Abbo_Schwarz_Notes_on_Floer_homology_and_loop_space_homology} the authors describe the
construction of an isomorphism between the Floer homology of a quadratically growing Hamiltonian function
in a cotangent bundle and the singular homology of the corresponding loop space. This isomorphism is proved to be compatible with
the pair-of-pants product on Floer homology and the loop product defined by Chas and Sullivan on loop space homology. Abbondandolo and
Schwarz provide a map $h$ between the Lagrangian Floer homology of a fixed cotangent fibre and Hamiltonian Floer homology. The map $h$
corresponds under their isomorphism to the map induced by the inclusion of the based loop space into the free loop space and is equal to $\chi$
in this context.

Floer homology for cotangent bundles for quadratically growing Hamiltonian functions is analytically much more demanding than Floer homology for
compactly supported Hamiltonian functions. On the other hand no bubbling issues are present since the symplectic
structure is exact. The key point in the work of Abbondandolo-Schwarz
\cite{Abbo_Schwarz_Notes_on_Floer_homology_and_loop_space_homology} is to relate the
product structures in Floer homology and those in loop space homology.\\[1ex]
\noindent Finally, we point out that the idea of applying the ideas of Piunikhin, Salamon and Schwarz from
\cite{Piunikhin_Salamon_Schwarz_Symplectic_Floer_Donaldson_theory_and_quantum_cohomology} to Lagrangian Floer homology is certainly not original.
When writing this article was nearly finished the author learned about the recent work of Kati{\'c} and Milinkovi{\'c}
\cite{Katic_Milinkovic_PSS_ismomorphism_for_Lagrangian_intersections} where the special case of zero-sections in cotangent bundles is treated.
Since the zero-section $L$ is exact, i.e.~the symplectic form on cotangent bundles $T^*L$ is exact, $\om=d\lambda$, and furthermore $\lambda|_{L}=0$,
neither holomorphic disks nor holomorphic spheres are present. In particular, the core part of this article, namely
understanding the limitations of the PSS-construction caused by bubbling-off of
holomorphic disks, is void for the case considered in \cite{Katic_Milinkovic_PSS_ismomorphism_for_Lagrangian_intersections}.

Another appearance of the techniques of Piunikhin, Salamon and Schwarz can be found in the work of
Barraud and Cornea \cite{Barraud_Cornea_Homotopical_dynamics_in_symplectic_topology} and Cornea and Lalonde
\cite{Cornea_Lalonde_Cluster_Homology:_an_overview_of_the_construction_and_results}. In the latter the framework for Lagrangian Floer homology
is extended to deal with bubbling by using cluster homology, cf. remark \ref{rmk:cornea_lalonde_use_PSS}.
Moreover, we want to remark that the purpose of this article is to demonstrate what can be achieved with ``classical'' tools
in (Lagrangian) Floer homology.

It should be mentioned that Lagrangian Floer homology can be defined using Novikov rings, see for instance
\cite[section 3.2]{Barraud_Cornea_Homotopical_dynamics_in_symplectic_topology} for details.

We close with a brief remark about the coefficient ring $\Z/2$. In certain cases it is possible to choose $\Z$ as coefficient ring, e.g.~if
the Lagrangian submanifold is relative spin, cf.~\cite{FOOO}. We will not pursue this direction is the present version of this article. The
same applies to non-compact symplectic manifolds which are convex at infinity or geometrically bounded.\\[1ex]
\textbf{Acknowledgements.}\quad
The author thanks Octav Cornea for interesting and helpful discussions. He also would like to thank the referees for their thorough and dedicated
work.

The author was financially supported by the German Research Foundation (DFG) through Priority Programm 1154
"Global Differential Geometry", grants SCHW 892/2-1 and AL 904/1-1, and by NSF Grants DMS-0102298 and DMS-0603957.
Part of this work was done when the author was affiliated with the University of Leipzig.
\section{Applications}\label{section:applications}
Theorem \ref{thm:existence_of_Lagrangian_PSS} has the following immediate corollaries.
\begin{Cor}\label{cor:partial_isomorphism}
In the situation of theorem \ref{thm:existence_of_Lagrangian_PSS}
\beq
\HF_k(L,\phi_H(L))\cong \H^{n-k}(L;\Z/2)
\eeq
for all values of $k$ satisfying $n-N_L+2\leq k\leq N_L-2$. In particular, we obtain a full isomorphism whenever $N_L\geq n+2$,
e.g.~in case $\om|_{\pi_2(M,L)}=0$, where $N_L=\infty$. This recovers
Floer's original result \cite{Floer_Morse_theory_for_Lagrangian_intersections} as well as Oh's
\cite[theorem II]{Oh_Floer_cohomology_spectral_sequences_and_the_Maslov_class}.
\end{Cor}

\begin{Rmk}
Since $L$ is closed corollary \ref{cor:partial_isomorphism} can be stated as
\beq
\HF_k(L,\phi_H(L))\cong \H_k(L;\Z/2)
\eeq
for all $n-N_L+2\leq k\leq N_L-2$, since the bounds on the degree are symmetric \wrt the change $k\mapsto n-k$.
\end{Rmk}

\begin{Cor}
In the situation of theorem \ref{thm:existence_of_Lagrangian_PSS}\beq
\#\big(L\pitchfork\phi_H(L)\big)\geq\sum_{k\,=\,n-N_L+2}^{N_L-2} b_k(L)\;,
\eeq
where $b_k(L)$ are the $\Z/2$-Betti numbers of $L$.
\end{Cor}

\begin{Def}\label{def:displaceable}
A closed Lagrangian submanifold is called $L$ \textit{displaceable} (or Hamiltonianly displaceable) in the closed symplectic manifold $(M,\om)$
if there exists a Hamiltonian diffeomorphism $\phi_H\in\Ham(M,\om)$ such that $L\cap\phi_H(L)=\emptyset$.
The \textit{displacement energy} of $L$ is $e(L):=\inf\{||H||\mid L\cap\phi_H(L)=\emptyset\}$.
Here $||H||$ denotes the \emph{Hofer norm} of $H$.
\end{Def}

\begin{Cor}
If $L$ is (Hamiltonianly) displaceable, then
\beq
\H_i(L;\Z/2)=0\quad\forall n-N_L+2\leq i\leq N_L-2
\eeq
\end{Cor}

\begin{Rmk}\label{rmk:pol_examples}
Following an construction by Audin, Polterovich \cite[theorem 4]{Polterovich_Monotone_Lagrange_submanifolds_of_linear_spaces}
provides Lagrangian embeddings of products of spheres $S^k$ into $\R^{2n}$.
More precisely, if $\sigma_k$ denotes the antipodal map, then for $r=2,\ldots,n$ the manifolds
\beq
\Big(\big(S^{r-1}\times S^1\big)/\big(\sigma_{r-1}\times\sigma_1\big)\Big)\times S^{n-r}
\eeq
embed as monotone Lagrangian submanifolds $L_r$ into $\R^{2n}$. Moreover, $L_r$ has minimal Maslov number $N_{L_r}=r$.
For $r\geq\frac{n}{2}+2$  they satisfy
\beq
\H_k(L_r;\Z/2)\;\;
    \begin{cases}
    \;=0&\text{if } n-r+2\leq k\leq r-2\\
    \;\not=0&\text{for }k=n-r+1\text{ and }k=r-1\;.
    \end{cases}
\eeq
Since any Lagrangian submanifold of $\R^{2n}$ is displaceable, this shows that the isomorphism statement in
theorem \ref{thm:existence_of_Lagrangian_PSS} is sharp.
\end{Rmk}

\begin{Def}
We set $K_L:=\max\big\{k\mid\H_k(L;\Z/2)\not=0,\;k\leq \lfloor\frac{n}{2}\rfloor\big\}$.
In particular, $0\leq K_L\leq\lfloor\frac{n}{2}\rfloor$ and $\H_{K_L+1}(L)=0,\ldots,\H_{n-K_L-1}(L)=0$ holds.
\end{Def}

Theorem \ref{thm:existence_of_Lagrangian_PSS} implies new restrictions for monotone Lagrangian submanifolds to be displaceable. We recall Oh's
result \cite[theorem II]{Oh_Floer_cohomology_spectral_sequences_and_the_Maslov_class} (cp.~corollary \ref{cor:partial_isomorphism} above)
asserting that the minimal Maslov number of a monotone, displaceable Lagrangian submanifolds has to satisfy $N_L\leq n+1$, where $n=\dim L$ .
This is sharpened as follows.

\begin{Cor}\label{cor:K_L+N_L_leq_n+1}
If $\HF_*(L,\phi_H(L))=0$ (e.g.~$L$ is displaceable) then $N_L+K_L\leq n+1$ holds.
\end{Cor}

\begin{proof}
We note that by definition $K_L\leq\frac{n}{2}$. If $N_L<\frac{n}{2}+2$ the assertion is trivially true.
If $N_L\geq\frac{n}{2}+2$ we know that $\H_k(L;\Z/2)=0$ for all $n-N_L+2\leq k\leq N_L-2$ by corollary \ref{cor:partial_isomorphism}.
In particular, $K_L<n-N_L+2$.
\end{proof}

\begin{Cor}
Let $L\subset M$ be a displaceable monotone Lagrangian submanifold which satisfies $\H_{\lfloor\frac{n}{2}\rfloor}(L;\Z/2)\not=0$,
then $N_L<\frac{n}{2}+2$.
\end{Cor}

According to Audin's conjecture the minimal Maslov number of all (not just monotone) Lagrangian embeddings
of tori into $\R^{2n}$ equals 2. Oh verified the conjecture for monotone Lagrangian tori $T^n$ up to $n\leq24$,
see \cite[theorem III]{Oh_Floer_cohomology_spectral_sequences_and_the_Maslov_class}. Although not overly restrictive, the above corollary is
the first global result, i.e.~valid for all dimensions, in this direction (besides $N_L\leq n+1$ which holds for any monotone Lagrangian submanifold
in $\R^{2n}$).

A closer inspection of the proof of theorem \ref{thm:existence_of_Lagrangian_PSS} yields that the above mentioned results for the minimal Maslov
number $N_L$ can be strengthened as follows. The assertion that $N_L=r$ implies (by definition) that there exists a disk $d:\D^2\pf M$ with boundary
on $L$ such that $\Mas(d)=r$. In all above cases this disk can be chosen to be a \emph{holomorphic disk}. Moreover, the energy $\om(d)$
of these holomorphic disks is smaller than the displacement energy $e(L)$ of the Lagrangian submanifold $L$
(cf.~definition \ref{def:displaceable}).

Corollary \ref{cor:partial_isomorphism} is complemented by

\begin{Prop}
We assume that we are in the situation of theorem \ref{thm:existence_of_Lagrangian_PSS} and that $N_L\geq\frac{n}{2}+2$.
If $\HF_{N_L-1}(L)\not\cong\H_{N_L-1}(L)$, then there exists a holomorphic disk $d$ realizing
the minimal Maslov number $N_L$, i.e.~$\Mas(d)=N_L$.
If there exists no holomorphic disk of minimal Maslov number $N_L$ the isomorphism in theorem \ref{thm:existence_of_Lagrangian_PSS} holds for all degrees
$n-2N_L+2\leq k\leq 2N_L-2$. Moreover, the homomorphism $\varphi_k$ and $\rho_k$ are defined for $k\leq2N_L-2$ and $k\geq n-2N_L+2$, respectively.
\end{Prop}

For the examples from remark \ref{rmk:pol_examples} we can therefore conclude that there exist holomorphic disks realizing the minimal Maslov number.
This of course can be iterated. Namely, if there exists no holomorphic disk of minimal Maslov number $jN_L$ then we obtain isomorphisms
for all degrees $n-jN_L+2\leq k\leq jN_L-2$.

If $L$ is displaceable there cannot exist an isomorphism in the top and bottom degree. This leads to the existence of a holomorphic
disk of Maslov index less than $n+1$. This already follows from Oh's paper \cite{Oh_Floer_cohomology_spectral_sequences_and_the_Maslov_class}.
Much more general results in this direction are known. Indeed, if
$L$ is displaceable and monotone then through \emph{each point} of $L$ passes a holomorphic disk of Maslov index less than $n+1$, cf.
\cite{Albers_On_the_extrinsic_topology_of_Lagrangian_submanifolds,Cornea_Lalonde_Cluster_Homology:_an_overview_of_the_construction_and_results}.
Moreover, in \cite{Cornea_Lalonde_Cluster_Homology:_an_overview_of_the_construction_and_results} a very similar result is given for general
(i.e.~possibly non-monotone) Lagrangian submanifolds.

We end this section with an application of theorem \ref{thm:diagram_commutes}.

\begin{Cor}
If the monotone, closed Lagrangian submanifold $L^n\subset(M^{2n},\om)$ has vanishing Floer homology, $\HF_*(L,\phi_H(L))=0$,
then the following holds for the homomorphisms induced by the inclusion map $\iota:L\hpf M$
\bea
\iota^k=0&:\H^k(M;\Z/2)\pf\H^k(L;\Z/2)\qquad\text{for }\;k\leq N_L-2\,,\\
\iota_k=0&:\H_k(L;\Z/2)\pf\H_k(M;\Z/2)\qquad\text{for }\;k\geq n-N_L+2\;.
\eea
\end{Cor}

This improves theorem \cite[theorem 1.1]{Albers_On_the_extrinsic_topology_of_Lagrangian_submanifolds} which asserts the statement (only for
$\iota_k$) in case that $L$ is displaceable. We note that the displaceablility is essential in the proof given in
\cite{Albers_On_the_extrinsic_topology_of_Lagrangian_submanifolds}.
\section{Elements of Floer theory}\label{section:recollection_of_Floer_homology}
We briefly recall the construction of the Floer complex $(\CF_i(L,\phi_H(L)),\partial_F)$ for a closed, monotone Lagrangian submanifold $L$ in
a closed, symplectic manifold $(M,\om)$.
\begin{Def}\label{def:minimal_Maslov_and_Chern_nb}
A Lagrangian submanifold $L$ of the symplectic manifolds $(M,\om)$ is called \emph{monotone}, if there exists a constant $\lambda>0$, such that
$\om|_{\pi_2(M,L)}=\lambda\cdot\Mas$, where $\Mas:\pi_2(M,L)\pf\Z$ is the \emph{Maslov index}.
This implies that the symplectic manifolds $(M,\om)$ is monotone as well, that is, $\om|_{\pi_2(M)}=\tilde{\lambda}\cdot c_1|_{\pi_2(M)}$, where
$\tilde{\lambda}>0$ and $c_1$ is the first Chern class of $(M,\om)$.

We define the \textit{minimal Maslov number} $N_L$ of $L$ as the positive generator of the image of the Maslov index
$\Mas(\pi_2(M,L))\subset\Z$. We set $N_L=+\infty$ in case $\Mas$ vanishes (this implies that $\om|_{\pi_2(M,L)}=0$). The
\emph{minimal Chern number} $N_M$ of $M$ is defined analogously.
\end{Def}
\noindent
For a Hamiltonian function $H:S^1\times M\pf\R$ the Floer complex $\big(\CF_i(L,\phi_H(L)),\partial_F\big)$ is generated over $\Z/2$ by the set
of Hamiltonian chords
\beq\label{eqn:Ham_chord}
\P_L(H):=\Big\{x\in C^{\infty}([0,1],M)\;\big|\;\dot{x}(t)=X_H\big(t,x(t)\big),\;x(0),x(1)\in L,\;[x]=0\in\pi_1(M,L)\Big\}
\eeq
i.e.~$\CF(L,\phi_H(L))=\P_L(H)\otimes\Z/2$. Let us explain some notions. First, $X_H$ is the (time dependent) Hamiltonian vector field generated
by the Hamiltonian function $H:S^1\times M\pf\R$ and is defined by $\om(X_H(t,\cdot),\,\cdot\,)=-dH(t,\cdot)$. The time-1-map
of the flow $\phi_H^t$ of $X_H$ is denoted by $\phi_H\equiv\phi_H^1$.

A certain subset of the intersection points $L\cap\phi_H(L)$ is often taken to generate a chain complex.
There is an one-to-one correspondence between $\P_L(H)$ and this subset by applying the flow $\phi_H^t$
to an intersection point. Furthermore, in either approach the Hamiltonian
function $H$ is required to be non-degenerate meaning that $L\pitchfork\phi_H(L)$.

The Maslov index defines a grading on $\P_L(H)$, which is only defined modulo the minimal Maslov number $N_L$ and
up to an overall shift. Let us briefly recall the construction of the grading.
Given two elements $x,y\in\P_L(H)$ we choose a map $u:[0,1]^2\rightarrow M$ s.t.~$u(0,t)=x(t)$, $u(1,t)=y(t)$ and $u(\tau,0),u(\tau,1)\in L$. According
to \cite{Viterbo_Maslov_index_for_Lagrangian_Floer_homology,Floer_A_relative_Morse_index_for_the_symplectic_action} a Maslov index is assigned to
the map $u$ as follows. Since the symplectic vector bundle $u^*TM$ is trivial, a loop of Lagrangian subspaces in $\R^{2n}$ is obtained
by following the Lagrangian subspaces of $TL$ along the two $u(\tau,0/1)$-sides of the strip and transporting them by the Hamiltonian flow along
the $u(0/1,t)$-sides (and flipping them by 90 degrees at the corners). To this loop of Lagrangian subspaces in $\R^{2n}$ the classical
Maslov index is assigned.

This gives rise to a \emph{relative} Maslov index for $x$ and $y$ which certainly depends on the choice of $u$.
Indeed, let $v:[0,1]^2\rightarrow M$ be another choice connecting $x$ and $y$ and let
$h:\D^2_+\rightarrow M$ be a half-disk realizing a homotopy of the chord $x$ to a constant path. We can form the disk $d:=h\#u\#(-v)\#(-h)$
with boundary on $L$, where $\#$ denotes concatenation and $-v$ is the map $(\tau,t)\mapsto v(-\tau,t)$. If the relative Maslov index of $x$ and $y$
is computed with help of either $u$ or $v$, the difference is given by the Maslov index $\Mas([d])$. We note that the Maslov index of $[d]$ does not
depend on the choice of the half-disk $h$.
Thus, we can assign a number $\mu(x,y)\in\Z/N_L$ to each pair $x,y\in\P_L(H)$. By construction, this number satisfies $\mu(x,z)=\mu(x,y)+\mu(y,z)$
for all $x,y,z\in\P_L(H)$. Therefore, we artificially set $\mu(x_0):=0$ for a fixed $x_0\in\P(H)$ and define the degree $\mu(y):=\mu(y,x_0)\in\Z/N_L$
for all other $y\in\P_L(H)$. Assigning index zero to another element in $\P_L(H)$ leads to a shift of the degree. Therefore, by this procedure we define
a mod $N_L$ grading on $\P_L(H)$ up to an overall shift.

In this article we choose to fix the shifting ambiguity by requiring that the dimension of the moduli space $\M(H;x)$ (see definition
\ref{def:moduli_space_of_cut_off_holo_strips} below), is given by $\mu(x)$ mod $N_L$. Equivalently, we could demand that the space $\M(x;H)$ (see
definition \ref{def:moduli_space_of_cut_off_holo_strips_II}) has dimension $n-\mu(x)$ mod $N_L$. More details can be found
in section \ref{section:the_Lagrangian_PSS_morphisms}. This convention is consistent by a gluing argument and additivity of the Fredholm index.

\begin{figure}[ht]
\begin{center}
\input{./perturbed_strip.pstex_t}
\caption{A perturbed holomorphic strip $u\in\M_L(x,y;J,H)$.}
\label{fig:perturbed strip}
\end{center}
\end{figure}
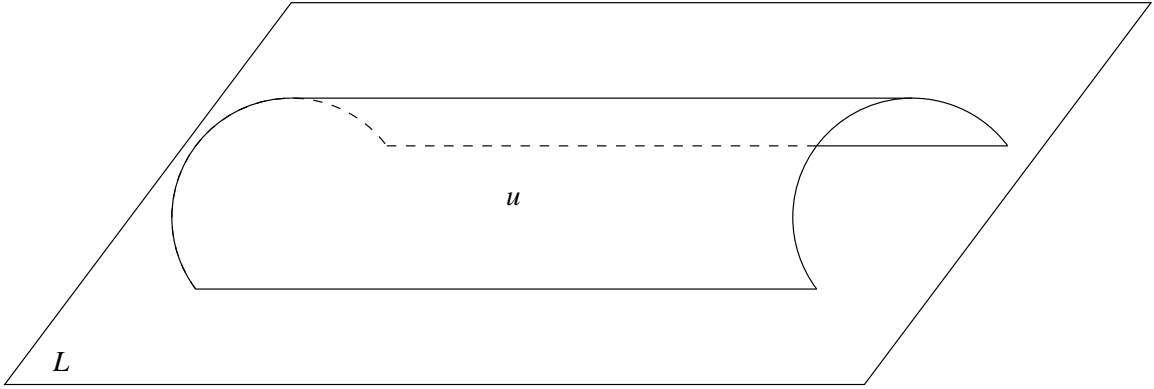

The Floer differential $\partial_F$ is defined by counting perturbed holomorphic strips
(a.k.a.~semi-tubes or Floer strips), see figure \ref{fig:perturbed strip}. For $x,y\in\P_L(H)$ we define the moduli spaces
\beq\label{eqn:moduli_space_Floer_strips}
\M_L(x,y;J,H):=\left\{\;u:\R\times[0,1]\pf M\left|\;\;
    \begin{aligned}
    &\partial_su+J(t,u)\big(\partial_tu-X_H(t,u)\big)=0\\
    &u(s,0),u(s,1)\in L\;\;\forall s\in\R\\
    &u(-\infty)=x,\;\;u(+\infty)=y
    \end{aligned}\;
\right.
\right\}
\eeq
\noindent
where $J(t,\cdot)$, $t\in[0,1]$, is a family of compatible almost complex structures on $(M,\om)$. If we would use
the intersection point $L\pitchfork\phi_H(L)$
to generate the Floer complex then the differential would be defined by counting unperturbed holomorphic strips having one boundary
component on $L$ and the other on $\phi_H(L)$. Again the flow of the Hamiltonian vector field provides a one-to-one correspondence between perturbed
and unperturbed strips.

\begin{Thm}[Floer]\label{thm:Floer_thm}
For a generic family $J$, the moduli spaces $\M_L(x,y;J,H)$ are smooth manifolds of dimension $\dim \M_L(x,y;J,H)\equiv\mu(y)-\mu(x)$ mod $N_L$, carrying
a free $\R$-action if $x\not=y$.
\end{Thm}

We note that the dimension of the moduli spaces is given by the Maslov index modulo the minimal Maslov number $N_L$. In other words, if we fix the
asymptotic data to be $x,y\in\P_L(H)$, the moduli space $\M_L(x,y;J,H)$ consists (in general) out of several connected components each of which has
dimension $\equiv\mu(y)-\mu(x)$ mod $N_L$.

\begin{Conv}
We set $\M_L(x,y;J,H)_{[d]}$ to be the union of the $d$-dimensional components.
\end{Conv}

\begin{Thm}[Floer, Oh]\label{thm:compactification_for_Floer_strips}
If the minimal Maslov number satisfies $N_L\geq2$ then for all $x,z\in\P_L(H)$ the moduli space
\beq
\Mh_L(x,z;J,H)_{[d-1]}:=\M_L(x,z;J,H)_{[d]}/\R
\eeq
is compact if $d=1$ and compact up to simple breaking if $d=2$, i.e.~it admits a compactification
(denoted by the same symbol) such that the boundary decomposes as follows
\beq
\partial\Mh_L(x,z;J,H)_{[1]}=\bigcup_{\substack{y\in\P_L(H)}}\Mh_L(x,y;J,H)_{[0]}\times\Mh_L(y,z;J,H)_{[0]}\,.
\eeq
\end{Thm}
\noindent
The boundary operator $\partial_F$ in the Floer complex is defined on generators $y\in\P_L(H)$ by
\beq
\partial_F(y):=\sum_{x\in\P_L(H)}\#_2\Mh_L(x,y;J,H)_{[0]}\cdot x
\eeq
and is extended linearly to $\CF_*(L,\phi_H(L))$.
Here, $\#_2\Mh_L(x,y;J,H)_{[0]}$ denotes the (mod 2) number of elements in $\Mh_L(x,y;J,H)_{[0]}$.

The two theorems above justify this definition of $\partial_F$, namely the sum is finite and $\partial_F\circ\partial_F=0$. The
\textit{Lagrangian Floer homology} groups are $\HF_*(L,\phi_H(L)):=\H_*(\CF(L,\phi_H(L),\partial_F)$.
It is an important feature of Floer homology that it is independent of the chosen family of
almost complex structures and invariant under Hamiltonian perturbations.
In particular, there exists an canonical isomorphism $\HF_*(L,\phi_H(L))\cong\HF_*(L,\phi_K(L))$ for any two Hamiltonian functions $H,K$.\\[0.5ex]
Floer theory is a (relative) Morse theory for the \emph{action functional} $\A_H$ defined on the space of paths in $M$ which start and end
on $L$ and are homotopic (relative $L$) to a constant path in $L$. By definition the action functional is
\beq\label{eqn:action_functional}
\A_H(x,d_x):=\int_{\D^2_+}d_x^*\om-\int_0^1H(t,x(t))dt
\eeq
where $d_x:\D^2_+\rightarrow M$ realizes a homotopy from a constant path to the path $x$. The value of the action functional depends only on the
relative homotopy class of $d_x$. Its critical points are exactly the elements of $\P_L(H)$.

Given a Hamiltonian function $H:S^1\times M\pf\R$ there exists \textit{Hamiltonian Floer homology} $\HF_*(H)$. It is generated
over $\Z/2$ by the set
\beq\label{eqn:Ham_orbits}
\P(H):=\Big\{x\in C^{\infty}(S^1,M)\;\big|\;\dot{x}(t)=X_H\big(t,x(t)\big),\;[x]=0\in\pi_1(M)\Big\}
\eeq
i.e.~$\CF_*(H)=\P(H)\otimes\Z/2$, and graded by the \textit{Conley-Zehnder index} $\CZ$. We normalize the Conley-Zehnder index by
requiring that $\CZ\big(t\mapsto a(-t)\big)=-\CZ(a)$ for all $a\in\P(H)$ and
that for $C^2$-small Morse functions $f$ we have $\CZ(a)=\frac12\dim M-\Morse(a)$ for all $a\in\P(f)=\Crit(f)$.

The set $\P(H)$ correspond to a subset of the fixed points of the time-1-map $\phi_H$. The Floer differential $\partial_F$
(we use the same notation as for Lagrangian Floer homology) is defined
by counting perturbed holomorphic \textit{cylinders} (instead of strips). For $a,b\in\P(H)$ we define the space
\beq\label{eqn:moduli_space_Floer_cylinders}
\M(a,b;J,H):=\left\{\;u:\R\times S^1\pf M\left|\;\;
    \begin{aligned}
    &\partial_su+J(t,u)\big(\partial_tu-X_H(t,u)\big)=0\\
    &u(-\infty)=a,\;\;u(+\infty)=b
    \end{aligned}\;
\right.
\right\}\;.
\eeq
Its elements are called Floer cylinders. For generic choices, $\M(a,b;J,H)$ is a smooth manifold of dimension $\CZ(b)-\CZ(a)$ modulo twice
the \textit{minimal Chern number} $N_M$ (see definition \ref{def:minimal_Maslov_and_Chern_nb}).
If $a\not=b$ then $\M(a,b;J,H)$ carries a free $\R$-action and the quotient is denoted by $\Mh(a,b;J,H)$. The zero-dimensional components of
the quotient are compact and the one-dimensional are compact up to simple breaking. In particular, the above construction can be carried out verbatim
and we obtain Hamiltonian Floer homology $\HF_*(H):=\H_*(\CF(H),\partial_F)$.

In contrast to Lagrangian Floer homology the construction of Hamiltonian Floer homology can be
accomplished in more general situations, for instance for semi-positive symplectic manifold. A very nice account with many
details of Hamiltonian Floer homology for monotone symplectic manifolds can be found in \cite{Salamon_lectures_on_floer_homology}. Extensions
to more general settings, e.g.~using Novikov rings, can be found in \cite{Hofer_Salamon_Floer_homology_and_Novikov_rings,
McDuff_Salamon_J_holomorphic_curves_and_symplectic_topology}.
%
%
%
\section{The Lagrangian Piunikhin-Salamon-Schwarz morphisms}\label{section:the_Lagrangian_PSS_morphisms}
%
%
\subsection{The construction}\label{subsection:construction_of_the_Lagrangian_PSS_morphisms}
$ $\\[1ex]
In this section we will describe the construction of the maps
\begin{gather}
\varphi_k:\CF_k(L,\phi_H(L))\pf\CM^{n-k}(L;\Z/2)\\[1ex]
\rho_k:\CM^{n-k}(L;\Z/2)\pf\CF_k(L,\phi_H(L))
\end{gather}
on the Floer resp.~Morse chain complex and show that they are well-defined and chain maps whenever $k\leq N_L-2$ respectively
$k\geq n-N_L+2$.

These maps are analogs of the Piunikhin-Salamon-Schwarz isomorphism, which is an
isomorphism $\PSS:\H^*(M)\stackrel{\cong}{\pf}\HF_*(H)$ between Hamiltonian Floer homology and singular homology. It is defined for semi-positive
symplectic manifolds $(M,\om)$. Due to transversality problems for moduli spaces of holomorphic disks, Lagrangian Floer homology
is only defined for monotone Lagrangian submanifolds $L$ with minimal Maslov number $N_L\geq 2$,
see \cite{Oh_Floer_cohomology_of_Lagrangian_intersections_and_pseudo-holomorphic_disks_I,
Oh_Addendum_Floer_cohomology_of_Lagrangian_intersections_and_pseudo-holomorphic_disks_I}.
Furthermore, a complete analogue of the PSS isomorphism cannot exist since Lagrangian Floer homology is \emph{not} isomorphic to the singular
homology in general.
\begin{Conv}\label{conv:Floer_generic_J}
Throughout the following sections
we fix a family of compatible almost complex $J_0(t,\cdot)$ for $t\in[0,1]$ for which the Fredholm operator given by the
shift-invariant Floer equation, see \eqref{eqn:moduli_space_Floer_strips} and \eqref{eqn:moduli_space_Floer_cylinders},
is surjective. That is, $J_0(t,\cdot)$ is generic in the sense of theorem \ref{thm:Floer_thm}.
\end{Conv}
The following constructions will be made with reference to the family $J_0$. Since all constructions are natural with respect to
the continuation homomorphisms in Floer homology they do not depend on the choice of $J_0$ in homology, cp.~also remark
\ref{rmk:compatibility_with_continuation}. We will suppress $J_0$ in the notation.
%
\subsubsection{Construction of $\varphi_*$}   \label{section:construction_of_varphi}
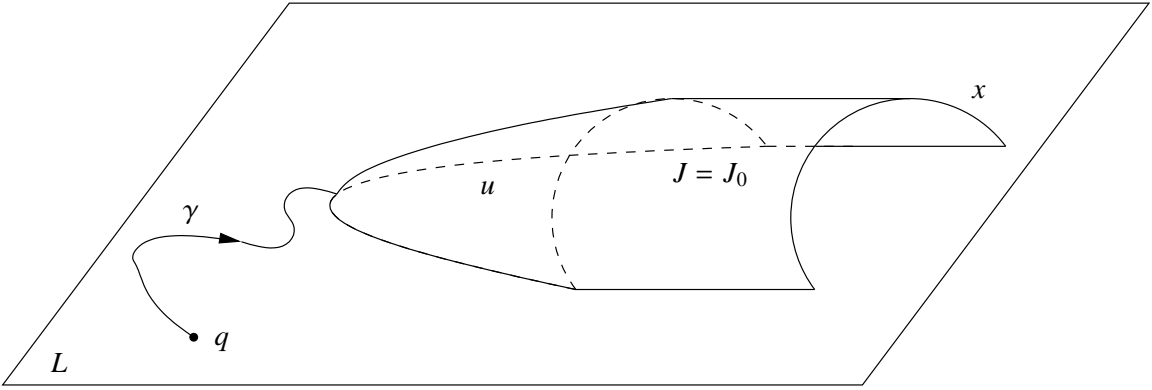
\begin{figure}[ht]
\begin{center}
\input{./PSS_strip_1.pstex_t}
\caption{An element $(\gamma,u)\in\M^{\varphi}(q,x)$.}
\label{fig:PSS_strip_1}
\end{center}
\end{figure}

%
$ $\\[1ex]
\noindent We define the moduli space $\M^{\varphi}(q,x;J,H,\beta,f,g)$ to consist of pairs $(\gamma,u)$ of maps
\begin{gather}
\gamma:(-\infty,0]\pf L\zw\text{and}\zw u:\R\times[0,1]\pf M\zw \text{with }E(u)=\!\!\int_{-\infty}^\infty\int_0^1|\partial_su|^2\,dt\,ds<\infty\\[1ex]
\text{solving}\qquad\dot{\gamma}+\nabla^{g} f\circ\gamma=0
        \quad\text{and}\quad\partial_su+J(s,t,u)\big(\partial_tu-\beta(s)X_{H}(t,u)\big)=0\,,
\end{gather}\\[-1ex]
where $\beta:\R\rightarrow[0,1]$ is a smooth cut-off function satisfying $\beta(s)=0$ for $s\leq0$ and $\beta(s)=1$ for $s\geq1$. The function
$f:L\rightarrow\R$ is a Morse function on $L$ and $\nabla^{g}$ is the gradient
\wrt a Riemannian metric $g$ on $L$. Moreover, $J(s,t,\cdot)$, $(s,t)\in\R\times[0,1]$, is a smooth family of compatible almost complex structures on $M$
which as a map $s\mapsto J(s,\cdot,\cdot)$ is constant outside $[0,1]$, and which agrees for $s\geq1$ with the fixed $J_0$ from
convention \ref{conv:Floer_generic_J}.
The pair $(\gamma,u)$ is required to satisfy the boundary conditions
\beq
\gamma(-\infty)=q\,,\quad\gamma(0)=u(-\infty)\,,\quad u(+\infty)=x\quad\text{and}\quad u(s,0),\,u(s,1)\in L\,,
\eeq
where $q\in \Crit(f)$ is a critical point of $f$ and $x\in\P_L(H)$. Due to the cut-off function $\beta$ the strip $u$ is holomorphic for $s\leq0$
and has, by assumption, finite energy $E(u)<\infty$. In particular, $u$ admits a continuous extension $u(-\infty)$, see
\cite[section 4.5]{McDuff_Salamon_J_holomorphic_curves_and_symplectic_topology}. In other words, topologically $u$ forms a half disk
$\{z\in\C: |z|\leq1\,,\;\mathrm{Re}(z)\geq0\}$ such that the boundary part $\{|z|=1\}$ is mapped to $L$ and  $\{\mathrm{Re}(z)=0\}$ is
mapped to the chord $x$ (see figure \ref{fig:PSS_strip_1}). For brevity we denote this moduli space by $\M^{\varphi}(q,x)$.

\begin{Def}\label{def:moduli_space_of_cut_off_holo_strips}
For $x\in\P_L(H)$ we set
\beq
\M(H;x):=\left\{u:\R\times[0,1]\pf M \left|\;\;
    \begin{aligned}
    &\partial_su+J(s,t,u)\big(\partial_tu-\beta(s)X_H(t,u)\big)=0\\
    &u(s,0),u(s,1)\in L\;\;\forall s\in\R\\
    &u(+\infty)=x,\;\;E(u)<+\infty
    \end{aligned}
    \right.
\right\}
\eeq
where, as above, $J$ is an element in the space of smooth families of compatible almost complex structures which as maps $s\mapsto J(s,\cdot,\cdot)$
are constant outside $[0,1]$, and which agree for $s\geq1$ with the fixed $J_0$ from
convention \ref{conv:Floer_generic_J}. We denote this space by $\J$.
\end{Def}

\begin{Rmk}
The maps in the moduli space $\M(x;H)$ are defined on strips with finite energy and we need to use the removal of singularity theorem
to fill in the point at $\pm\infty$, respectively. Instead we could have defined the maps on half-disk with a strip-like end. This
is, in fact, easily achieved by a conformal reparametrization. We will expand on this in section \ref{section:varphi_o_rho=id_transversality}.
\end{Rmk}

\begin{Def}\label{def:admissible_Hamiltonian}
A non-degenerate Hamiltonian function $H:S^1\times M\pf\R$ is called admissible if none of the (finitely many) element in $\P_L(H)\cup\P(H)$
(see equations \eqref{eqn:Ham_chord}, \eqref{eqn:Ham_orbits}) are constant.
In other words, for each $x\in M$ there exists $t\in S^1$ such that $\nabla H(t,x)\neq0$.
\end{Def}

\begin{Rmk}
Fix a Hamiltonian function $H$. If there exists a fixed $x\in M$ where $\nabla H(t,x)=0$ vanishes $\forall t\in S^1$, then an arbitrarily small
perturbation of $H$ will loose this property. In particular, being admissible is a generic property.
\end{Rmk}

\begin{Thm}\label{thm:transversality_for_M^phi}
For an admissible Hamiltonian function $H$, fixed cut-off function $\beta$ as above and a Morse-Smale pair $(f,g)$ there exists a generic subset
$\J^{\mathrm{reg}}$ of $\J$,
such that for $J\in\J^{\mathrm{reg}}$ the moduli spaces $\M^{\varphi}(q,x)=\M^{\varphi}(q,x;J,H,\beta,f,g)$ are smooth manifolds
of dimension
\beq
\dim\M^{\varphi}(q,x)=\mu(x)-n+\Morse(q)\quad\mathrm{mod}\;N_L\,,
\eeq
where $n=\dim L$. We denote by $\M^{\varphi}(q,x)_{[d]}$ the union of the $d$-dimensional components and remark that $\M^{\varphi}(q,x)$
carries no (natural) $\R$-action.
\end{Thm}

\begin{proof}
We are required to prove that the
moduli space $\M(H;x)$ given in definition \ref{def:moduli_space_of_cut_off_holo_strips} above is a smooth manifold and that the evaluation map
\beq
\big(\ev_0\times\ev_{-\infty}\big)(\gamma,u):=(\gamma(0),u(-\infty))\subset L\times L
\eeq
is transversal to the diagonal $\Delta\subset L\times L$. This is contained as special cases of the proofs
of theorems \ref{thm:tranversality_for_M^phi_o_rho} and \ref{thm:detailed_transversality} in section \ref{section:varphi_o_rho=id_transversality}.
\end{proof}
\noindent
Compactness properties are a more subtle issue. We
employ an index argument to rule out bubbling. Since the bubbling occurs for the perturbed strips we give the following definition. We recall that
$\beta:\R\rightarrow[0,1]$ is a smooth cut-off function satisfying $\beta(s)=0$ for $s\leq0$ and $\beta(s)=1$ for $s\geq1$.

\begin{Rmk}
A part of the proof of theorem \ref{thm:transversality_for_M^phi} is to show that for generic choices, $\M(H;x)$ is a smooth manifold.
We recall that by our convention, the moduli space $\M(H;x)$ has dimension $\mu(x)$ mod $N_L$.
As usual $\M(H;x)_{[d]}$ denotes the union of its $d$-dimensional components.

The moduli space $\M^{\varphi}(q,x)$ defined above is given by the intersection of the unstable manifold $W^u(q)$ of the critical point $q\in\Crit(f)$
with $\M(H;x)$. In particular, $\dim\M^{\varphi}(q,x)=\dim\M(H;x)-n+\Morse(q)$.
\end{Rmk}

\begin{Thm}\label{thm:compactness_Lagrangian_PSS}
Let $L$ be a closed, monotone Lagrangian submanifold of minimal Maslov number $N_L$. Then for all
$d<N_L$ the moduli spaces $\M(H;x)_{[d]}$ are compact up to adding broken strips.
\end{Thm}

\begin{proof}
Using Floer's equation the following inequality is easily derived for an element $u\in\M(H;x)$.
\beq
0\leq E(u)\leq\A_H(x,u)+\sup_MH\;.
\eeq
For convenience we compiled the computations in appendix \ref{appendix:energy_estimates}. We refer to the appendix for the notation and
note that for each two elements $u,v\in\M(H;x)_{[d]}$ the Maslov index of $u\#(-v)$ equals zero, $\Mas(u\#(-v))=0$. Otherwise, $u$ and $v$ could not
lie in components of the same dimension.
Since $L$ is monotone this implies $\om(u)=\om(v)$ and thus the value of the action functional $\A_H(x,u)=\A_H(x,v)$.

We obtain a uniform bound on the energy for elements in the moduli space $\M(H;x)_{[d]}$. Hence, we can apply Gromov's compactness theorem.
In particular, sequences converge up to breaking and bubbling. We claim that, under the assumption $d<N_L$, only breaking occurs.

Indeed, if a sequence $(u_n)\subset\M(H;x)_{[d]}$ converges in the Gromov-Hausdorff topology,
i.e.~$u_n\rightharpoonup(u_\infty,v_1,\dots,v_\Gamma,s_1,\dots,s_\Sigma,d_1,\dots,d_\Delta)$, where
\begin{itemize}
\item $u_\infty\in\M(H;x_0)_{[d]}$,\\[-1ex]
\item $v_\gamma\in\Mh_L(x_{\gamma-1},x_{\gamma};J,H)$, where $x_\gamma\in\P_L(H)$
        and $x_{\Gamma}=x$,\\[-1ex]
\item $\{s_\sigma\}$ are holomorphic spheres,\\[-1ex]
\item $\{d_\delta\}$ are holomorphic disks,\\[-1ex]
\end{itemize}
then (see \cite{Floer_symplectic_fixed_points_and_holomorphic_spheres},
\cite[proposition 3.7]{Oh_Floer_cohomology_of_Lagrangian_intersections_and_pseudo-holomorphic_disks_I})
the Fredholm index of the linearized Cauchy-Riemann operator $\F$ behaves as follows:
\beq\label{eqn:Fredholm_index_formula_bubbling}
\ind\F_{u_n}=\ind\F_{u_\infty}+\sum_\gamma\ind\F_{v_\gamma}+\sum_\sigma2c_1(s_\sigma)+\sum_\delta\Mas(d_\delta)\,.
\eeq
We know by transversality that $\ind\F_{v_\gamma}\geq 0$. Moreover, $2c_1(s_\sigma)\geq N_L$ and $\Mas(d_\delta)\geq N_L$ given that
the maps $s_\sigma$ and $d_\delta$ are non-constant. By assumption, we have $N_L>d=\ind\F_{u_n}$. Since $\ind\F_{u_\infty}\geq 0$ the theorem follows.
\end{proof}

Theorem \ref{thm:compactness_Lagrangian_PSS} implies that the moduli space $\M^{\varphi}(q,x)$ is compact up to breaking as long as the space
$\M(H;x)$ is compact up to breaking, i.e.~if $\dim\M(H;x)<N_L$.

Moreover, the standard gluing techniques (see for instance
\cite{Schwarz_Matthias_PhD,McDuff_Salamon_J_holomorphic_curves_and_symplectic_topology})
imply that the moduli spaces
$\M^{\varphi}(q,x)$ can be compactified by adding broken Morse trajectories or Floer strips, as long as no bubbling occurs. In particular, we can
compactify the one-dimensional moduli space $\M^{\varphi}(q,x)_{[1]}$ by adding either a Morse trajectory or a Floer strip,
see figure \ref{fig:PSS_strip_1_breaking}.

\begin{Cor}\label{cor:compactification_of_M^varphi(q,x)}
If $\dim \M(H;x)<N_L$ we conclude,
\bea\label{eqn:compactification_of_M^varphi(q,x)}
\partial\M^{\varphi}(q,x)_{[1]}=&\bigcup_{p\in\Crit(f)}\Mh(q,p;f,g)_{[0]}\times\M^{\varphi}(p,x)_{[0]}\\
                                \cup&\bigcup_{y\in\P_L(H)}\M^{\varphi}(q,y)_{[0]}\times\Mh_L(y,x;J_0,H)_{[0]}\;.
\eea
Here, $\Mh(q,p;f,g)$ is the moduli space of (unparametrized, negative) gradient-flow trajectories of the Morse function $f$ on $L$
which run from $q$ to $p$. The gradient is taken with respect to the metric $g$. Moreover, $J_0$ is as in convention \ref{conv:Floer_generic_J}.
\end{Cor}
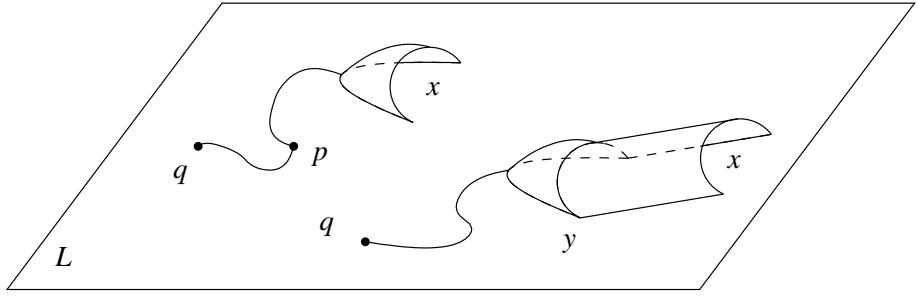
\begin{figure}[ht]
\begin{center}
\input{./PSS_strip_1_breaking.pstex_t}
\caption{Broken configurations in $\partial\M^{\varphi}(q,x)_{[1]}$.}
\label{fig:PSS_strip_1_breaking}
\end{center}
\end{figure}
\begin{Def}
For $k\leq N_L-1$ we define on generators
\bea
\varphi_k:\CF_k(L,\phi_H(L))&\pf\CM^{n-k}(L;\Z/2)\\
    x\quad&\mapsto\quad\sum_q\#_2\M^{\varphi}(q,x)_{[0]}\cdot q\;.
\eea
The map $\varphi_k$ is extended by linearity.
\end{Def}

\begin{Thm}\label{thm:varphi_well_defined_and_chain_hom}
For $k\leq N_L-1$ the homomorphism $\varphi_k$ is well-defined and for $k\leq N_L-2$ a chain morphism, that is
$\varphi_{k+1}\circ\partial_F=\delta^L\circ\varphi_k$.
\end{Thm}

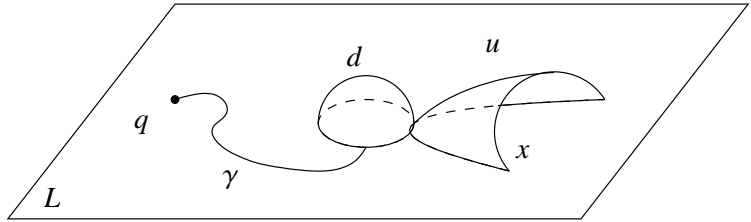
\begin{figure}[ht]
\begin{center}
\input{./PSS_strip_1_bubbling.pstex_t}
\caption{A bubbling configuration.}
\label{fig:PSS_strip_1_bubbling}
\end{center}
\end{figure}

\begin{proof}
If $k\leq N_L-1$ we know by theorem \ref{thm:compactness_Lagrangian_PSS} that all moduli spaces $\M^{\varphi}(q;x)_{[d]}$ are compact
if $x\in\CF_k(L,\phi_H(L))$ and $d=0$. In particular, the sum in the definition of $\varphi_k$ is finite and the map is therefore well-defined.

If $k\leq N_L-2$, we know by theorem \ref{thm:compactness_Lagrangian_PSS} that all moduli spaces $\M^{\varphi}(q;x)_{[d]}$ are compact
if $x\in\CF_k(L,\phi_H(L))$ and $d=0,1$.
Equation \eqref{eqn:compactification_of_M^varphi(q,x)} in corollary \ref{cor:compactification_of_M^varphi(q,x)} implies that $\varphi_k$
is a chain morphism.
\end{proof}
\noindent
In particular, in homology we obtain homomorphisms
\beq
\varphi_k:\HF_k(L,\phi_H(L))\pf\H^{n-k}(L;\Z/2)\quad\text{for }k\leq N_L-2\;.
\eeq
Before we define the map $\rho_k$ we want to give a more geometric picture of the problems leading to the restrictions $k\leq N_L-2$. The most
problematic issue arises when bubbling occurs at the continuous extension $u(-\infty)$ for $u\in\M(H;x)$. Let us assume that there exists a sequence
$(\gamma_n,u_n)\subset\M^{\varphi}(q,x)$ such that this kind of bubbling occurs. Then we obtain as limit object
(see figure \ref{fig:PSS_strip_1_bubbling}) a gradient half-trajectory
$\gamma$, an element $u\in\M(H;x)$ and a holomorphic disk $d:(\D^2,\partial\D^2)\pf(M,L)$, such that
\beq
\gamma(-\infty)=q,\quad\gamma(0)=d(z_0),\quad d(z_1)=u(-\infty),\quad u(+\infty)=x\,,
\eeq
where $z_0,z_1$ are two point in $\partial\D^2$. This also can be viewed as breaking of the strip $u$ near $-\infty$.
\begin{Rmk}\label{rmk:cornea_lalonde_use_PSS}
In general, there is no way around dealing with the bubbling-off phenomenon. The cluster complex approach by Cornea and
Lalonde \cite{Cornea_Lalonde_Cluster_Homology:_an_overview_of_the_construction_and_results} is build exactly to incorporate this phenomenon.
Namely the above configuration $(\gamma,d,u)$ is not considered a boundary point but as an interior point,
where the "other side" consists out of a configuration $(\gamma,d,\delta,u)$, where $\delta:[0,R]\pf L$ is a finite length gradient flow
line, such that $d(z_1)=\delta(0)$ and $\delta(R)=u(-\infty)$.
\end{Rmk}

%
\subsubsection{Construction of $\rho_*$}   \label{section:construction_of_rho}
%
$ $\\[1ex]
The homomorphism $\rho_k$ is the mirror image of $\varphi_k$.
The moduli space $\M^{\rho}(x,q)\equiv\M^{\rho}(x,q;J,H,\beta,f,g)$ consists of pairs $(u,\gamma)$ of maps (see figure \ref{fig:PSS_strip_2})
\begin{gather}
u:\R\times[0,1]\pf M\zw \text{with }E(u)=\!\!\int_{-\infty}^\infty\int_0^1|\partial_su|^2\,dt\,ds<\infty\zw\text{and}\zw\gamma:[0,\infty)\pf L \\[1ex]
\text{solving}\qquad\partial_su+J(-s,t,u)\big(\partial_tu-\beta(-s)X_{H}(t,u)\big)=0
        \quad\text{and}\quad\dot{\gamma}+\nabla^{g} f\circ\gamma=0\,,
\end{gather}\\[-1ex]
where $f,g,J,\beta$ are as at the beginning of section \ref{section:construction_of_varphi}. In the equation we replace
$J(s,\cdot,\cdot)$ and $\beta(s)$ by $J(-s,\cdot,\cdot)$ and $\beta(-s)$. In particular, $J(-s,\cdot,\cdot)=J_0(t,\cdot)$
for $s\leq-1$. The pair $(\gamma,u)$ is required to satisfy boundary conditions
\beq
u(s,0),\,u(s,1)\in L\,,\quad u(-\infty)=x\,,\quad u(+\infty)=\gamma(0)\quad\text{and}\quad\gamma(+\infty)=q\,,
\eeq
where $q\in \Crit(f)$ is a critical point of $f$ and $x\in\P_L(H)$. Because of the change of the sign in the cut-off function the strip $u$ is
holomorphic for $s\geq0$ and thus, $u$ admits a continuous extension $u(+\infty)$.
\begin{figure}[ht]
\begin{center}
\input{./PSS_strip_2.pstex_t}
\caption{An element $(u,\gamma)\in\M^{\rho}(x,q)$.}
\label{fig:PSS_strip_2}
\end{center}
\end{figure}
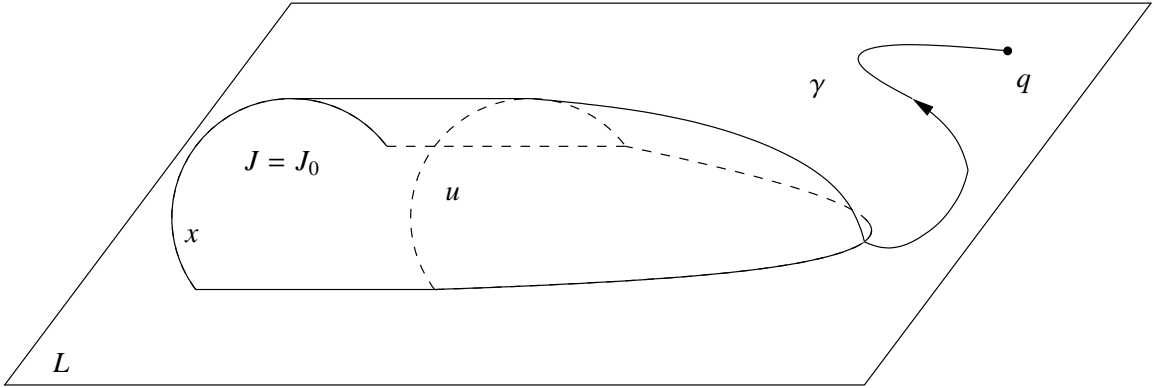
\begin{Def}\label{def:moduli_space_of_cut_off_holo_strips_II}
Analogously to $\M(H;x)$ (see definition \ref{def:moduli_space_of_cut_off_holo_strips}) we set for $x\in\P_L(H)$
\beq
\M(x;H):=\left\{u:\R\times[0,1]\pf M \left|\;\;
    \begin{aligned}
    &\partial_su+J(-s,t,u)\big(\partial_tu-\beta(-s)X_H(t,u)\big)=0\\
    &u(s,0),u(s,1)\in L\;\;\forall s\in\R\\
    &u(-\infty)=x,\;\;E(u)<+\infty
    \end{aligned}
    \right.
\right\}
\eeq
\end{Def}

As before for generic choices (see below for a precise statement) these moduli spaces are smooth manifolds. The dimension is given by
\beq
\dim\M(x;H)=n-\mu(x)\quad\mathrm{mod}\;N_L\;.
\eeq
The moduli space $\M^{\rho}(x,q)$ is the intersection of $\M(x;H)$ with the stable manifold $W^s(q)$ of the critical point $q\in\Crit(f)$.

\begin{Thm}
For an admissible Hamiltonian function $H$, fixed cut-off function $\beta$ and a Morse-Smale pair $(f,g)$ there exists a generic subset
$\J^{\mathrm{reg}}$ of $\J$,
such that for $J\in\J^{\mathrm{reg}}$ the moduli spaces$\M^{\rho}(x,q)=\M^{\rho}(x,q;J,H,\beta,f,g)$ are smooth manifolds
of dimension
\beq
\dim\M^{\rho}(x,q)=n-\mu(x)-\Morse(q)\quad\mathrm{mod}\;N_L\,,
\eeq
where $n=\dim L$. The space $\J$ and admissibility is defined in definitions
\ref{def:moduli_space_of_cut_off_holo_strips} and \ref{def:admissible_Hamiltonian}, respectively.
\end{Thm}

\begin{proof}
This is completely analogous to the proof of theorem \ref{thm:transversality_for_M^phi}.
\end{proof}

As before a subscript $[d]$ denotes the union of the $d$-dimensional components.
The compactness properties of $\M(x;H)$ are controlled in exactly the same way as for $\M(H;x)$ in theorem \ref{thm:compactness_Lagrangian_PSS}.

\begin{Thm}\label{thm:compactness_Lagrangian_PSS_partII}
Let $L$ be a closed monotone Lagrangian submanifold of minimal Maslov number $N_L$. Then for all
$d<N_L$ the moduli spaces $\M(x;H)_{[d]}$ are compact up to adding broken strips.
\end{Thm}

\begin{proof}
The proof of theorem \ref{thm:compactness_Lagrangian_PSS_partII} and of the following corollary is literally the same as for
theorem \ref{thm:compactness_Lagrangian_PSS} up to a slight change of the uniform energy estimate for a solution $u\in\M(x;H)$,
(see appendix \ref{appendix:energy_estimates}):
\beq
0\leq E(u)\leq-\A_H(x,u)-\inf_MH\;.
\eeq
\end{proof}

\begin{Cor}\label{cor:compactification_of_M^rho(x,q)}
If $\dim \M(H;x)<N_L$, the zero dimensional moduli space is compact and the
one dimensional moduli space $\M^{\rho}(x,q)_{[1]}$ can be compactified such that
\bea\label{eqn:compactification_of_M^rho(x,q)}
\partial\M^{\rho}(x,q)_{[1]}=&\bigcup_{p\in\Crit(f)}\M^{\rho}(x,p)_{[0]}\times\Mh(p,q;f,g)_{[0]}\\
                                \cup&\bigcup_{y\in\P_L(H)}\Mh_L(x,y;J_0,H)_{[0]}\times\M^{\rho}(y,q)_{[0]}\,,
\eea
where $J_0$ is as in convention \ref{conv:Floer_generic_J}.
\end{Cor}
\noindent
This leads to

\begin{Def}
For $k\geq n-N_L+1$ we define on generators
\bea
\rho_k:\CM^{n-k}(L;\Z/2)&\pf\CF_k(L,\phi_H(L))\\
    q\quad&\mapsto\quad\sum_x\#_2\M^{\varphi}(x,q)_{[0]}\cdot x
\eea
which then is extended by linearity.
\end{Def}
\noindent
The following theorem is proved in exactly the same way as theorem \ref{thm:varphi_well_defined_and_chain_hom}.

\begin{Thm}\label{thm:rho_well_defined_and_chain_hom}
For $k\geq n-N_L+1$ the homomorphism $\rho_k$ is well-defined and for $k\geq n-N_L+2$ it is a chain morphism, that is
$\rho_{k-1}\circ\delta^L=\partial_F\circ\rho_k$.
\end{Thm}
\noindent
We obtain two homomorphisms

\begin{gather}
\varphi_k:\HF_k(L,\phi_H(L))\pf\H^{n-k}(L;\Z/2)\quad\text{for }k\leq N_L-2\;,\\[1ex]
\rho_k:\H^{n-k}(L;\Z/2)\pf\HF_k(L,\phi_H(L))\quad\text{for }k\geq n-N_L+2\;.
\end{gather}

\begin{Rmk}
For trivial reasons the maps $\varphi_k$  and $\rho_k$ vanish for $k<0$ and $k>n$, namely for these degrees
the corresponding moduli spaces are empty, since there are no critical points of Morse index less than 0 or larger than $n$.
This statement has non-trivial implications when the minimal Maslov number $N_L$ exceeds $\dim L+1$, e.g.~Floer homology
vanishes in the corresponding degrees, see \cite[theorem 3.1]{Seidel_Graded_Lagrangian_submanifolds} for applications.
\end{Rmk}
We end this section with the following
\begin{Rmk}\label{rmk:compatibility_with_continuation}
The homomorphisms $\varphi_*$ and $\rho_*$ are natural with respect to the change of the Hamiltonian function.
For the moment we denote by $\varphi_*^H$ and $\varphi^{K}_*$ the homomorphisms obtained by using
two different Hamiltonian function $H,K:S^1\times M\pf\R$. Then $\varphi_*^H$ and $\varphi^{K}_*$ are intertwined by the continuation homomorphisms
$\Omega^{KH}_*:\HF_*(L,\phi_H(L))\stackrel{\cong}{\pf}\HF_*(L,\phi_{K}(L))$, i.e.~$\varphi^K_*\circ\Omega^{KH}_*=\varphi_*^H$.
The same holds true for changing the Morse function. This is proved by a suitable cobordism argument as used in the following
sections. In particular, $\varphi_*$ and $\rho_*$ do not depend on the fixed $J_0$, cp.~convention \ref{conv:Floer_generic_J}.
\end{Rmk}
%
\subsubsection{Construction of the original Piunikhin-Salamon-Schwarz isomorphism}   \label{section:construction_of_real_PSS}
%
$ $\\[1ex]
We now recall very briefly the original construction by Piunikhin-Salamon-Schwarz from
\cite{Piunikhin_Salamon_Schwarz_Symplectic_Floer_Donaldson_theory_and_quantum_cohomology} which leads to an isomorphism
$\PSS:\H^{n-k}(M)\stackrel{\cong}{\pf}\HF_k(H)$ between singular cohomology of the underlying manifold  and the Hamiltonian Floer
homology for all semi-positive symplectic manifolds. Since the maps $\varphi$ and $\rho$ are analogs of the map $\PSS$ we
will just indicate how to modify the moduli space $\M^{\rho}(x,p)$ (see definition \ref{def:moduli_space_of_cut_off_holo_strips_II} above). Let us recall
that elements $(u,\gamma)\in\M^{\rho}(x,p)$ are a Floer strip $u$ with Hamiltonian term cut off near $+\infty$ and a gradient flow half-line $\gamma$
in $L$ which satisfy $u(+\infty)=\gamma(0)$. Then the original PSS-construction is as follows. Replace Floer strips by
Floer cylinders (still with the cut-off) and gradient flow half-lines
in $L$ by such in $M$. The condition $u(+\infty)=\gamma(0)$ remains. Counting zero-dimensional components defines a chain map $\PSS$.
This gives rise to a moduli space $\M^{\mathrm{PSS}}(a,q)$.

The analogous modification of $\varphi$ gives rise to another moduli space $\M^{\mathrm{PSS,\,inv}}(q,a)$ defining a
map $\PSS^{-1}$, which actually is the inverse map to $\PSS$. The compactness properties of the moduli spaces
defining $\PSS$ are better than those of $\M^{\varphi}(q,x)$ and $\M^{\rho}(x,p)$  since the only bubbling phenomenon is bubbling-off of holomorphic
spheres (and no holomorphic disks). Holomorphic spheres are either somewhere injective or multiply covered and therefore, the theorem analogous
to theorems \ref{thm:compactness_Lagrangian_PSS} and \ref{thm:compactness_Lagrangian_PSS_partII} hold with much less restrictive assumptions and for
more general symplectic manifolds.

\subsection{The isomorphism property}\label{subsection:_the_isomorphism_property}
$ $\\[1ex]
In this section we will prove the second part of theorem \ref{thm:existence_of_Lagrangian_PSS} asserting
that, whenever both homomorphisms $\varphi_k$ and $\rho_k$ are defined simultaneously, they are inverse to each
other. This only occurs when $n-N_L+2\leq k\leq N_L-2$, that is, for $2N_L\geq n+4$.

For both compositions, $\rho_k\circ\varphi_k$ and $\varphi_k\circ\rho_k$, we will describe cobordisms which relate the counting
defining the composition to the counting defining the identity map.
%
\subsubsection{$\rho_k\circ\varphi_k=\id_{\HF_k(L,\phi_H(L))}$} \label{section:rho_o_varphi=id}
\begin{figure}[ht]
\begin{center}
\input{./composition_easy.pstex_t}
\caption{A cobordism proving $\rho_k\circ\varphi_k=\id_{\HF_k(L,\phi_H(L))}$.}
\label{fig:composition_easy}
\end{center}
\end{figure}
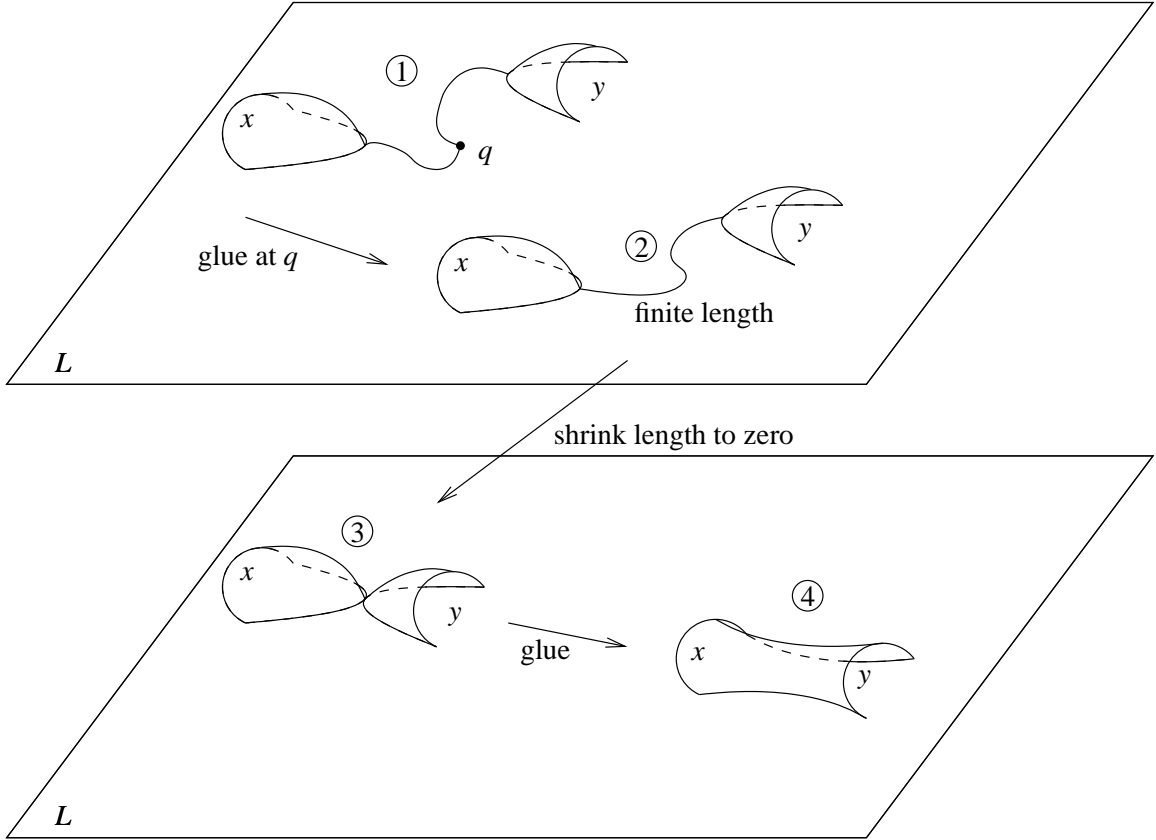

%
$ $\\[1ex]
The proof is analogous to the original approach in \cite{Piunikhin_Salamon_Schwarz_Symplectic_Floer_Donaldson_theory_and_quantum_cohomology}.
Geometrically, the cobordism relating the composition to the identity is given by the following steps, see figure \ref{fig:composition_easy}.

\begin{enumerate}
\item The composition $\rho_k\circ\varphi_k$ is a map from Floer homology to Floer homology. The coefficient
of $\rho_k\circ\varphi_k(y)$ in front of $x\in\P_L(H)$ is given by counting all zero dimensional configurations
$(u_+,\gamma_+;\gamma_-,u_-)$, where $(u_+,\gamma_+)\in\M^{\rho}(x,q)$ and $(\gamma_-,u_-)\in\M^{\varphi}(q,y)$ and $q\in\Crit(f)$ is arbitrary.

\item We glue the two gradient flow half-trajectories $\gamma_+$ and $\gamma_-$ at the critical point $q$ and obtain $(u_+,\Gamma,u_-)$, where
$\Gamma$ is a \emph{finite length} gradient flow trajectory, say parameterized by $[0,R]$, such that $u_+(+\infty)=\Gamma(0)$ and
$\Gamma(R)=u_-(-\infty)$.

\item We shrink the length $R$ of the gradient flow trajectory to zero. In the limit $R=0$ we obtain a pair $(u_+,u_-)$ of two maps
$u_+,u_-:\R\times[0,1]\pf M$ which satisfy Floer's equation on one half and are holomorphic on the $\mp$-half of the strip. Furthermore,
they satisfy $u_+(+\infty)=u_-(-\infty)$ and $u_+(-\infty)=x$ and $u_-(+\infty)=y$.

\item Since $(u_+,u_-)$ both are holomorphic around the point $u_+(+\infty)=u_-(-\infty)$ we can glue them and obtain a map
$u:\R\times[0,1]\pf M$ which satisfies Floer's equation with Hamiltonian term given by $H$ up to a compact perturbation around $s=0$.

\item We remove the compact perturbation and obtain an honest Floer strip connecting $x$ to $y$. Since we count zero dimensional configurations
(and are not dividing out the $\R$-action),
this is only non-zero if $y=x$, in which case there is exactly one such strip, namely the constant one.

In other words, in homology the coefficient $\rho_k\circ\varphi_k(y)$ in front of $x$ equals zero or one depending on whether $x=y$.
\end{enumerate}

We realize steps (1) to (3) of the above geometric picture by the following moduli space
$\M^{\rho\circ\varphi}(x,y)\cong\M^{\rho\circ\varphi}(x,y;J,H,\beta,f,g)$, which consists out of
quadruples $(R,u_+,\Gamma,u_-)$, where
\beq
R\geq0,\quad u_\pm:\R\times[0,1]\pf M,\quad\Gamma:[0,R]\pf L
\eeq
satisfying
\begin{gather}
\partial_su_\pm+J(\pm s,t,u_\pm)\big(\partial_tu_\pm-\beta(\pm s)X_H(t,u_\pm)\big)=0\,,\\[1.5ex]
u_+(-\infty)=x,\;u_-(+\infty)=y,\quad u_\pm(s,0),\,u_\pm(s,1)\in L,\quad E(u_\pm)<+\infty\,,\\[1.5ex]
\dot{\Gamma}(t)+\nabla^g f\circ\Gamma(t)=0\,,\\[1.5ex]
u_+(+\infty)=\Gamma(0),\quad\Gamma(R)=u_-(-\infty)\;.
\end{gather}
The Morse function $f$ on $L$, the metric $g$, the cut-off function $\beta$, the family of compatible almost complex structures $J$ and
the Hamiltonian function $H$ are chosen as at the beginning of section \ref{section:construction_of_varphi}.

\begin{Thm}
For an admissible Hamiltonian function $H$, fixed cut-off function $\beta$ and a Morse-Smale pair $(f,g)$ there exists a generic subset
$\J^{\mathrm{reg}}$ of $\J$ (see definition \ref{def:moduli_space_of_cut_off_holo_strips}),
such that for $J\in\J^{\mathrm{reg}}$ the moduli spaces $\M^{\rho\circ\varphi}(x,y)=\M^{\rho\circ\varphi}(x,y;J,H,\beta,f,g)$ are smooth manifolds
\beq
\dim\M^{\rho\circ\varphi}(x,y)=\mu(y)-\mu(x)+1\quad\text{mod }N_L\;,
\eeq
where the "+1" accounts for the parameter $R$. We denote by $\M^{\rho\circ\varphi}(x,y)_{[d]}$ the union of the $d$-dimensional components.
\end{Thm}

\begin{proof}
This is proved in the same way as for the other moduli spaces, cf.~e.g.~theorem \ref{thm:transversality_for_M^phi}.
\end{proof}

\begin{Thm}\label{thm:compactness_for_isomorphism_easy_case}
For $x,y\in\P_L(H)$ satisfying $\mu(y)\leq N_L-1$ and $\mu(x)\geq n-N_L+1$ the moduli space $\M^{\rho\circ\varphi}(x,y)_{[0]}$
is compact. If $\mu(y)\leq N_L-2$ and $\mu(x)\geq n-N_L+2$ the moduli space $\M^{\rho\circ\varphi}(x,y)_{[d]}$
is compact up to breaking in dimension $d=0,1$. In particular, if $d=0$ the moduli space is compact and for $d=1$ we conclude
\bea
\partial\M^{\rho\circ\varphi}(x,y)_{[1]}&=\big\{(R,u_+,\Gamma,u_-)\mid R=0\big\}\\[1.5ex]
                                      &\cup\bigcup_{x'\in\P_L(H)}\Mh_L(x,x';J_0,H)_{[0]}\times\M^{\rho\circ\varphi}(x',y)_{[0]}\\[1.5ex]
                                      &\cup\bigcup_{y'\in\P_L(H)}\M^{\rho\circ\varphi}(x,y')_{[0]}\times\Mh_L(y',y;J_0,H)_{[0]}\\[1.5ex]
                                      &\cup\bigcup_{q\in\Crit(f)}\M^{\rho}(x,q)\times\M^{\varphi}(q,y)
\eea
where $J_0$ is as in convention \ref{conv:Floer_generic_J}.
\end{Thm}

\begin{proof}
As explained in section \ref{section:construction_of_varphi} respectively \ref{section:construction_of_rho} the only source of non-compactness
apart from breaking is bubbling-off. In theorems \ref{thm:compactness_Lagrangian_PSS} respectively \ref{thm:compactness_Lagrangian_PSS_partII} we used
that the Fredholm index drops at least by the minimal Maslov number $N_L$ when a bubble appears. From this we concluded that no bubbling occurs for
sequences in the moduli spaces $\M(H;x)$ respectively $\M(x;H)$. The very same argument works for the moduli spaces
$\M^{\rho\circ\varphi}(x,y)_{[d]}$ for $d=0,1$, but now we have to exclude bubbling for both moduli spaces, $\M(H;y)$ and $\M(x;H)$.
This is guaranteed for $d=0$ by the assumptions $\mu(y)\leq N_L-1$ and $\mu(x)\geq n-N_L+1$ and for $d=1$
by $\mu(y)\leq N_L-2$ and $\mu(x)\geq n-N_L+2$, cp.~the proofs of theorems \ref{thm:varphi_well_defined_and_chain_hom} and
\ref{thm:rho_well_defined_and_chain_hom}.

Since no bubbling occurs, the only non-compactness is due to breaking. The usual gluing arguments show that we can compactify the moduli spaces by adding
broken trajectories.

This is encoded in the formula for the boundary of (the compactification) of the one dimensional moduli space as follows.
First, $R=0$ is an obvious boundary. Second, for a sequence $(R_n,u_+^{(n)},\Gamma_n,u_-^{(n)})\subset\M^{\rho\circ\varphi}(x,y)$ either $(R_n)$
converges or diverges. In the former case either $u_+^{(n)}$ or $u_-^{(n)}$ breaks, creating an perturbed holomorphic strip
in $\Mh_L(x,x';J,H)_{[0]}$ or $\Mh_L(y',y;J,H)_{[0]}$ and another solution in $\M^{\rho\circ\varphi}(x',y)$ or $\M^{\rho\circ\varphi}(x,y')$.
This comprises union two and three. In case $R_n\rightarrow\infty$ the sequence of finite length gradient trajectories $(\Gamma_n)$ breaks
into two gradient half-trajectories, creating a solution
in $\M^{\rho}(y,q)$ and $\M^{\varphi}(q,y)$ which both connect to the same critical point $q\in\Crit(f)$. This is expressed in the last union.
\end{proof}
\noindent
Theorem \ref{thm:compactness_for_isomorphism_easy_case} leads to the following definition of a chain homotopy.

\begin{Def}
For $n-N_L+1\leq k\leq N_L-2$ we define on generators the map
\bea
\Theta^{\rho\circ\varphi}_k:\CF_k(L,\phi_H(L))&\pf\CF_{k+1}(L,\phi_H(L))\\
                                y\quad&\mapsto\quad\sum_x\#_2\M^{\rho\circ\varphi}(x,y)_{[0]}\cdot x\;.
\eea
$\Theta^{\rho\circ\varphi}_k$ is extended linearly.
\end{Def}
\noindent
Before we prove that $\Theta^{\rho\circ\varphi}_k$ is a chain homotopy we set
\bea
\vartheta_k:\CF_k(L,\phi_H(L))&\pf\CF_k(L,\phi_H(L))\\
                                y\quad&\mapsto\quad\sum_x\#_2\big\{(R,u_+,\Gamma,u_-)\in\M^{\rho\circ\varphi}(x,y)_{[1]}\mid R=0\big\}\cdot x
\eea
and note that the set $\big\{(R,u_+,\Gamma,u_-)\in\M^{\rho\circ\varphi}(x,y)_{[1]}\mid R=0\big\}$ is zero dimensional and compact by the same
argument as in theorem \ref{thm:compactness_for_isomorphism_easy_case}. Since this set is a boundary, $\vartheta$ descends to homology, which also
follows immediately from the next corollary.

\begin{Cor}
For $n-N_L+1\leq k\leq N_L-2$ the map $\Theta^{\rho\circ\varphi}_k$ is well-defined and the following equality holds (where signs are arbitrary
as we use $\Z/2$-coefficients).
\beq
\partial_F\circ\Theta^{\rho\circ\varphi}_k-\Theta^{\rho\circ\varphi}_{k-1}\circ\partial_F=\rho_k\circ\varphi_k-\vartheta_k
\eeq
\end{Cor}

\begin{proof}
$\Theta^{\rho\circ\varphi}_k$ resp.~$\Theta^{\rho\circ\varphi}_{k-1}$ is well-defined since the sum is finite by theorem
\ref{thm:compactness_for_isomorphism_easy_case}
case $d=0$. The case $d=1$ and the definition of $\vartheta_k$ immediately imply the equation of the corollary.
\end{proof}

\begin{Cor}
In homology
\beq
\rho_k\circ\varphi_k=\vartheta_k:\HF_k(L,\phi_H(L))\pf\HF_k(L,\phi_H(L))
\eeq
holds for all $n-N_L+2\leq k\leq N_L-2$.
\end{Cor}
\noindent
For the rest of this section we will prove the

\begin{Claim}
In homology the identity $\vartheta_k=\id_{\HF_k(L,\phi_H(L))}$ holds.
\end{Claim}
\noindent
This is realized by the remaining steps (4) and (5) of our geometric picture from the beginning of this section.
First, we note that the set $\big\{(R,u_+,\Gamma,u_-)\in\M^{\rho\circ\varphi}(x,y)_{[1]}\mid R=0\big\}$  equals the set
\beq
\big\{(u_+,u_-)\in\M(x;H)_{[0]}\times\M(H;y)_{[0]}\mid u_+(+\infty)=u_-(-\infty)\big\}\;.
\eeq
This moduli space is cobordant to the following moduli space
\beqn
\widetilde{\M}_L(x,y;J,H):=\left\{u:\R\times[0,1]\pf M\left|\;\;
    \begin{aligned}
    &\partial_su+\widetilde{J}(s,t,u)\big(\partial_tu-\alpha(s)X_H(t,u)\big)=0\\
    &u(s,0),u(s,1)\in L\;\;\forall s\in\R\\
    &u(-\infty)=x,\;\;u(+\infty)=y
    \end{aligned}
\right.
\right\}
\eeq
where $\alpha$ is a cut-off function such that $\alpha(s)=0$ for $|s|\leq1$ and $\alpha(s)=1$ for $|s|\geq2$ and $\widetilde{J}(s,t,\cdot)$
is $s$-dependent only where $\alpha'(s)\neq0$ and equals $J_0(t,\cdot)$ for $|s|\geq2$, cp.~convention \ref{conv:Floer_generic_J}.
For the closed case a very thorough and detailed proof of the analogous statement is contained
in \cite{McDuff_Salamon_J_holomorphic_curves_and_symplectic_topology} (see theorem 10.1.2 as well as section 12.1).
For the case at hand check \cite[proposition 23.2 and lemma 23.4]{FOOO}. A detailed account is contained in
\cite[section 4]{Biran_Cornea_Quantum_Structures}.

The elements of the moduli space $\widetilde{\M}_L(x,y;J,H)$ solve almost the same equation as the elements in $\M_L(x,y;J_0,H)$, namely up to
the compact perturbation introduced by $\alpha$.
A one parameter family of Hamiltonian terms from $\alpha(s)X_H(t,\cdot)$ to $X_H(t,\cdot)$ as well as from $\widetilde{J}(s,t,\cdot)$ to the
$s$-independent family $J_0(t,\cdot)$ gives rise to another chain homotopy relating
$\vartheta_k$ to the map defined by counting elements in $\M_L(x,y;J_0,H)_{[0]}$. But the latter space carries a free $\R$-action as long
as $x\not=y$. In
particular, the zero dimensional components are empty in case $x\not=y$ and contain only the constant solution $x$, otherwise. Thus,
$\vartheta_k$ equals in homology the identity map $\id_{\HF_k(L,\phi_H(L))}$.
We will not write out explicit moduli spaces for the last chain homotopy but refer the reader to the original work of Piunikhin, Salamon and
Schwarz \cite{Piunikhin_Salamon_Schwarz_Symplectic_Floer_Donaldson_theory_and_quantum_cohomology} as well as to the book
\cite[section 12.1]{McDuff_Salamon_J_holomorphic_curves_and_symplectic_topology}.

This proves the claim and we conclude that if $n-N_L+2\leq k\leq N_L-2$ we have
\beq
\rho_k\circ\varphi_k=\id_{\HF_k(L,\phi_H(L))}\;.
\eeq
%
\subsubsection{$\varphi_k\circ\rho_k=\id_{\H^{n-k}(L;\Z/2)}$}  \label{section:varphi_o_rho=id_transversality}
%
$ $\\[1ex]
In the proof of the identity $\varphi_k\circ\rho_k=\id_{\H^{n-k}(L;\Z/2)}$  we will glue two perturbed holomorphic strips
along a common chord. Thus, the glued object differs from those appeared before and we have to argue differently to establish compactness.
Furthermore, we cannot proceed as Piunikhin, Salamon and
Schwarz did in \cite{Piunikhin_Salamon_Schwarz_Symplectic_Floer_Donaldson_theory_and_quantum_cohomology} because they use the dichotomy between somewhere
injective and multiply covered which is true for holomorphic spheres but does not hold for holomorphic disks.\\[.5ex]
Let us start again with the geometric picture.

\begin{enumerate}
\item The coefficient of $\varphi_k\circ\rho_k(p)$ in front of $q\in\Crit(f)$ is given by counting zero dimensional
configurations $(\gamma_-,u_-;u_+,\gamma_+)$ such that $(\gamma_-,u_-)\in\M^{\varphi}(q,x)$ and $(u_+,\gamma_+)\in\M^{\rho}(x,p)$
for some $x\in\P_L(H)$, see figure \ref{fig:composition_hard}.

\item We glue $u_-$ and $u_+$ at the chord $x\in\P_L(H)$ and obtain a single strip $u:\R\times[0,1]\pf M$ which is a solution of Floer's equation.
The important fact to note is, that the Hamiltonian term in the Floer equation is zero outside a compact subset of $\R\times[0,1]$
and that $u$ satisfies $\gamma_-(0)=u(-\infty)$ and $u(+\infty)=\gamma_+(0)$.
Furthermore, the Morse indices of $q$ and $p$ are equal. The set of triples $(\gamma_-,u,\gamma_+)$ as described above
is obtained by intersecting the space of maps $u$ with the unstable manifold of $q$ and the stable manifold of $p$. In particular, the space formed by
the maps $u$ has to be of dimension $n=\dim L$. This implies that the integral of the symplectic form $\om$ over $u$ vanishes: $\om(u)=0$.
At this point the monotonicity of the Lagrangian submanifold $L$ is essential.

\item The compact perturbation by the Hamiltonian term can be removed and we end up with triples $(\gamma_-,u,\gamma_+)$, where $u$ is a holomorphic
map $u:\R\times[0,1]\pf M$ (of finite energy) satisfying $\gamma_-(0)=u(-\infty)$ and $u(+\infty)=\gamma_+(0)$. Thus, $u$ is a \emph{holomorphic
disk} with boundary on the Lagrangian submanifold $L$.

\item The integral $\om(u)$ vanishes and therefore, $u$ has to be constant and $(\gamma_-,\gamma_+)$ form
an gradient flow line from $q$ to $p$. Again we are interested in zero dimensional configuration and we are not dividing by the $\R$-action.
By the same arguments as before we obtain the identity map $\id_{\H^{n-k}(L;\Z/2)}$.
\end{enumerate}
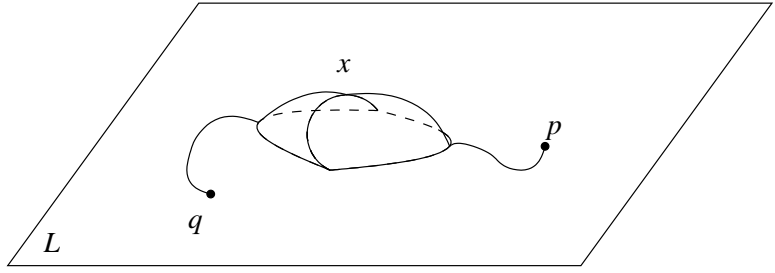
\begin{figure}[ht]
\begin{center}
\input{./composition_hard.pstex_t}
\caption{The composition $\varphi_k\circ\rho_k$.}
\label{fig:composition_hard}
\end{center}
\end{figure}
\noindent
The above geometric idea is realized as follows. We consider the following domain $\Omega$ in $\C$.
We denote by $\D_-=\{z\in\C\mid 2z-i+6\in\D^2,\;\mathrm{Re}(2z-i+6)\leq0\}$ and
$\D_+=\{z\in\C\mid 2z-i-6\in\D^2,\;\mathrm{Re}(2z-i-6)\geq0\}$. Furthermore, we denote $I:=[-3,3]\times[0,1]$ and set
\beq
\Omega:=\D_-\cup I\cup\D_+\,.
\eeq
For convenience we denote by $\pm\zeta:=\tfrac12(\pm7+i)\in\D_\pm$. We fix a sufficiently small $\varepsilon>0$ and consider
the following cut-off function $\alpha$ on $\Omega$
\beq
\alpha(z)= \begin{cases}
0&\text{for }z\in\D_-\cup\D_+\\
1 &\text{for }-2-\varepsilon\leq \mathrm{Re}(z)\leq 2+\varepsilon\\
\end{cases}
\eeq
Moreover, we require $0<\alpha(z)<1$ for all other $z\in\Omega$ and $\partial_t\alpha(s,t)=0$ for $(s,t)\in I$.
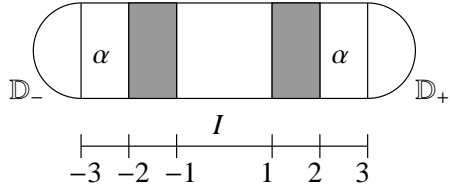
\begin{figure}[htb]
\input{omega.pstex_t}
\caption{The set $\Omega$. In the boxes containing the symbol $\alpha$ the cut-off function is non-constant.}\label{fig:omega}
\end{figure}
\begin{figure}[htb]
\input{rescale.pstex_t}
\caption{The conformal change. In the boxes containing the symbol $\alpha$ the cut-off functions are non-constant.}\label{fig:reform}
\end{figure}
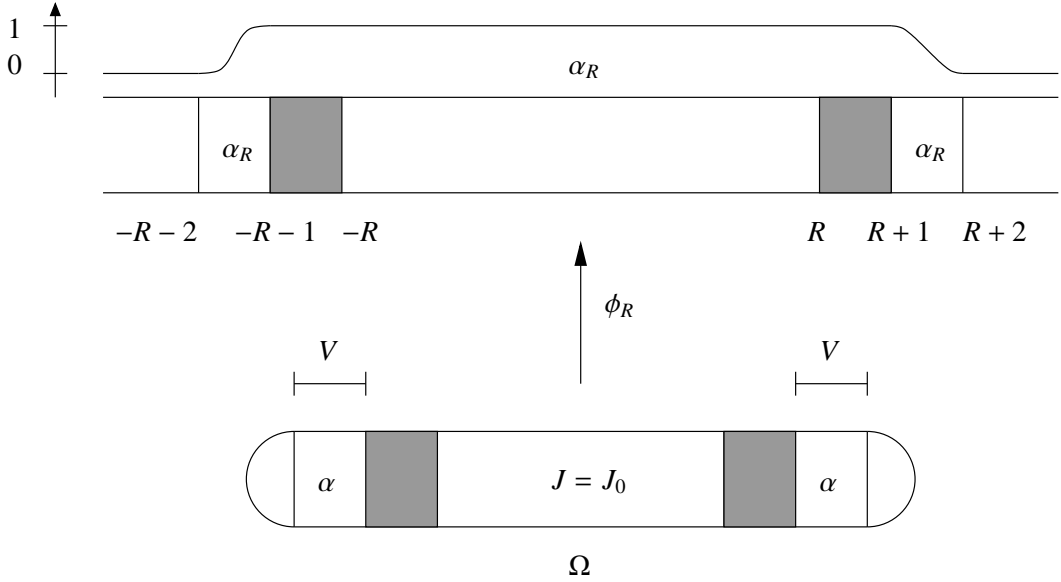
\noindent For $R\geq1$ we define a family of conformal structures on $\Omega$ by choosing diffeomorphisms
\beq
\phi_R:\Omega\setminus\{-\zeta,+\zeta\}\pf\R\times[0,1]
\eeq
satisfying:
\begin{itemize}
\item $[-1,1]\times[0,1]\subset\Omega$ is mapped onto $[-R,R]\times[0,1]$ via $(s,t)\mapsto(Rs,t)$,
\item $\D_{-}\setminus\{-\zeta\}\cup[-3,-2]\times[0,1]$ is mapped \textit{conformally} onto $(-\infty,-R-1]\times[0,1]$, moreover,
we require that $\varphi^R_-\circ\phi_R$ is independent of $R$, where $\varphi^R_-(s,t):=(s+(R-1),t)$,
\item $[2,3]\times[0,1]\cup\D_{+}\setminus\{+\zeta\}$ is mapped \textit{conformally} onto $[R+1,+\infty)\times[0,1]$, moreover,
we require that $\varphi^R_+\circ\phi_R$ is independent of $R$, where $\varphi^R_+(s,t):=(s-(R-1),t)$,
\item $\phi_1$ is a conformal throughout.
\end{itemize}
\noindent
We denote by $j_R:=\phi_R^*j$ the conformal structure obtained by pulling back the standard structure $j$ from $\R\times[0,1]$.
In other words, $j_R$ is the standard conformal structure on the set $V:=[-3,-2]\times[0,1]\cup[2,3]\times[0,1]$.
On $[-1,1]\times[0,1]$ we have $j_R\partial_s=R\partial_t$. Moreover, $j_0$ is the standard conformal structure on $\Omega$.
Finally, we can push forward the cut-off function $\alpha$ via $\phi_R$ to obtain $\alpha_R:\R\pf[0,1]$.

\begin{Rmk}\label{rmk:conformal_change}
The conformal change $(\Omega\setminus\{\pm\zeta\},j_R)\cong(\R\times[0,1],j)$ from the disk to the strip relates the
following discussion to the geometric picture we gave above. Moreover, we note that from the conformal change it is easy to see
that in the limit $R\rightarrow\infty$ the disk $\Omega$ breaks into two half disks with strip-like ends.
\end{Rmk}

We recall from convention \ref{conv:Floer_generic_J} that we fixed a smooth family $J_0(t,\cdot)$ of compatible almost complex structures
throughout.
We denote by $\Jh$ the set of all $C^k$-families $J(z,\cdot)$, $z\in\Omega$, of compatible almost complex structures which
as $C^k$-maps $z\mapsto J(z,\cdot)$ are non-constant only in the closed set $V$. Moreover, on $[-1,1]\times[0,1]\subset\Omega$
we require that $J(s,t,\cdot)=J_0(t,\cdot)$.

An element $J\in\Jh$ and $j_R$ defines a Cauchy-Riemann operator $\overline{\partial}_{j_R,J}$.

\begin{Def}\label{def:1_form_hamiltonian}
For a Hamiltonian function $H$ we denote by $\kappa(H)$ the 1-form with values
in $TM$ which is defined over $\Omega\setminus(\D_-\cup\D_+)=[-3,3]\times[0,1]$ in $(s,t)$-coordinates by
\beq
\kappa(H)(s,t,\cdot)=-dt\otimes X_H(t,\cdot)-ds\otimes J(s,t,\cdot)X_H(t,\cdot)\,.
\eeq
\end{Def}

\noindent So far we assumed $R\geq1$. We set
\beq
\tilde{\alpha}_R:=\begin{cases}\,R\alpha &\text{for }0\leq R\leq1\\\,\alpha &\text{for }R\geq1\end{cases}
\eeq
and we extend the definition of $j_R$ by defining $j_R:=j_0$ for $0\leq R\leq1$.
We recall from definition \ref{def:admissible_Hamiltonian} that a non-degenerate Hamiltonian function $H:S^1\times M\pf\R$ is called admissible
if for all $x\in L$ there exists $t\in S^1$ such that $\nabla H(t,x)\not=0$.
\begin{Def}\label{def:moduli_space_for_perturbed_hol_disks}
Let $H$ be an admissible Hamiltonian function and $J\in\Jh$.
For $a\in\R$ we define the moduli space $\M(J,H;a)$ to be the set of pairs $(R,u)$ satisfying
\beq\begin{cases}
\quad u:\Omega\pf M,\quad R\geq0\\
\quad\overline{\partial}_{j_R,J}u(z)+\tilde{\alpha}_R(z)\kappa(H)(t,u(z))=0\\
\quad u(\partial\Omega)\subset L\\
\quad \Mas(u)=a
\end{cases}\eeq
\end{Def}
\noindent
We remark that, since $\tilde{\alpha}_R(z)\not=0$ implies $z=(s,t)\in \Omega\setminus(\D_-\cup\D_+)$, the
expression $\tilde{\alpha}_R(z)\kappa(H)(t,u(z))$ is well-defined over $\Omega$.

\begin{Thm}\label{thm:detailed_transversality}
For $a\leq0$, an admissible Hamiltonian function $H$ and a generic $J\in\Jh$, the moduli space $\M(J,H;a)$ is a smooth manifold of dimension
$\dim \M(J,H;a)=\dim L+a+1$.
\end{Thm}

\begin{Rmk}
We will provide a proof of the above theorem within the category of $C^k$-maps for sufficiently large $k$. Then the usual circle of arguments
due Taubes extends this then to the $C^\infty$-setting, see \cite[Chapter 3.2]{McDuff_Salamon_J_holomorphic_curves_and_symplectic_topology}
\end{Rmk}

\begin{proof}
We set $\B:=\{u\in\H^{1,p}(\Omega,M)\mid u(\partial\Omega)\subset L,\;\Mas(u)=a\}$ for some $p>1$. The proof relies on two claims.
\begin{Claim}\label{claim:surj_lin}
For $R\geq0$ the linearization $D_{(R,J,u)}$ of the operator
\bea
\O:\R_{>0}\times\Jh\times\B&\pf L^p(\Omega,\Lambda T^*\Omega\otimes TM)\\
(R,J,u)&\mapsto\overline{\partial}_{j_R,J}u(z)+\tilde{\alpha}_R(z)\kappa(H)(t,u(z))
\eea
is surjective at solutions, i.e.~for $(R,J,u)$ s.t.~$\O(R,J,u)=0$. In particular,
the universal moduli space $\mathscr{M}:=\O^{-1}(0)$ is a Banach manifold.
\end{Claim}
\begin{proof}[Proof of the Claim \ref{claim:surj_lin}.] We will in fact prove more, namely the linearization of $\O$ at a solution $(R,J,u)$
is surjective already when restricted to $\{0\}\times T_J\Jh\times T_u\B\subset \R\times T_J\Jh\times T_u\B$.
An explicit expression for the linearization at $(R,J,u)$ is
\bea
D_{(R,J,u)}:\{0\}\times T_J\Jh\times T_u\B&\pf L^p(\Omega,\Lambda T^*\Omega\otimes u^*TM)\\
(0,Y,\xi)&\mapsto Y(z,u)\circ\Big(du-\tilde{\alpha}_R(z)dt\otimes X_H(t,u(z))\Big)\circ j_R + D_u^J\xi
\eea
which again is well-defined on all of $\Omega$ due to the cut-off function $\tilde{\alpha}_R$. The last operator equals
(when written in coordinates $(s,t)\in\Omega\setminus(\D_-\cup\D_+)$)
\beq
D_u^J\xi=\nabla_s\xi+J(s,t,u)\nabla_t\xi+\nabla_\xi J(s,t,u)\partial_tu+\nabla_\xi[\tilde{\alpha}_R(s)\nabla H(t,u)]\;.
\eeq

\noindent On $\D_-\cup\D_+$ the last term vanish by definition.
Since $D_u^J$ is a Fredholm operator, the range of $D_{(R,J,u)}$ is also closed. Thus, arguing by contradiction, we assume that there exists an element
$\eta\in L^q(\Omega,\Lambda T^*\Omega\otimes u^*TM)$ s.t.
\begin{gather}
\int_\Omega<\eta,D_u^J\xi>_{J}dz=0\qquad\text{for all }\xi\\
\int_\Omega <\eta,Y(z,u)\circ\Big(du-\tilde{\alpha}_R(z)dt\otimes X_H(t,u(z))\Big)\circ j_R>_{J}dz=0\qquad\text{for all }Y
\end{gather}
$<,>_{J}$ is the $z$-dependent Riemannian metric induced by $J(z,\cdot)$.
Choosing $Y=0$ the first equation tells us that $\eta$ is a solution of the formally adjoint of the $\bar{\partial}$-operator $D_u^J$.
Thus, by elliptic regularity, $\eta$ is of class $C^k$. To prove that $\eta$ in fact is identically zero, we need the second claim.

\begin{Claim}\label{claim:non-vanishing}
For $R>0$ there exists $z_0\in V$ s.t.
\beq
du(z_0)-\tilde{\alpha}_R(z_0)dt\otimes X_H(t_0,u(z_0))\not=0\;.
\eeq
\end{Claim}

\begin{proof}[Proof of Claim \ref{claim:non-vanishing}.]
Let us assume that $du(z)\tilde{\alpha}_R(z)dt\otimes X_H(t,u(z))=0$ for all $z\in V$. We recall that $u$ solves
\bea
0&=\overline{\partial}_{j_R,J}u+\tilde{\alpha}_R(z)\kappa(H)(t,u(z))\\
&=J\circ du(z)\circ j_R+du(z)+\tilde{\alpha}_R(z)(-dt\otimes X_H(t,u(z))-ds\otimes JX_H(t,u(z)))
\eea
We remark that, by construction, $j_R=j_0$ on the set $V$. We conclude
\beq
J\circ du(z)\circ j_0-\tilde{\alpha}_R(z)ds\otimes JX_H(t,u(z))=0
\eeq
for all $z=(s,t)\in V$. Thus, we can evaluate the last equation on $\partial_t$, entailing
\beq
[J\circ du(z)\circ j_0-\tilde{\alpha}_R(z)ds\otimes JX_H(t,u(z))]\cdot\partial_t=J\circ du(\partial_s)=0
\eeq
that is $\partial_su(s,t)=0\;\forall(s,t)\in V$, i.e.
\beq
\partial_tu(t)=\tilde{\alpha}_R(s)X_H(t,u(t)) \quad\forall (s,t)\in V\,.
\eeq
Since by assumption
\beq
\tilde{\alpha}_R(s)\;\begin{cases}=\min\{R,1\}\quad\text{for }|s|\leq 2+\varepsilon\\<\min\{R,1\}\quad\text{else} \end{cases}
\eeq
this implies
\beq
\partial_tu(t)=0\quad\text{and}\quad X_H(t,u(t))=0\quad\forall t\in[0,1]
\eeq
and since $u(0),u(1)\in L$ we conclude that there exists $x\in L$ where $\nabla H(t,x)=0$ for all $t\in S^1$. This contradicts the admissibility of $H$
(see definition \ref{def:admissible_Hamiltonian}) and finishes the proof of Claim \ref{claim:non-vanishing}.
\end{proof}

\noindent
We assume for now that $R>0$. Let us come back to the equation
\beq
\int_\Omega <\eta,Y(z,u)\circ\Big(du-\tilde{\alpha}_R(z)dt\otimes X_H(t,u(z))\Big)\circ j_R>_{J}dz=0\qquad\text{for all }Y
\eeq
We choose $z_0\in V$ provided by Claim 2 and $Y_0$
(see e.g.~\cite{Salamon_Zehnder_Morse_theory_for_periodic_solutions_of_Hamiltonian_systems_and_the_Maslov_index}) at $z_0$ such that
\beq
<\eta(z_0),Y(z_0,u(z_0))\cdot\Big(\partial_tu(z_0)-\tilde{\alpha}_R(z_0)X_H(t_0,u(z_0))\Big)>_J\not=0\;.
\eeq
We multiply $Y_0$ with a cut-off function supported in the vicinity of $z_0$ and extend it to an element $Y\in T_J\Jh$. Since we allow $J$ to be
$z$-dependent this is easy. In particular, there is no problem with multiple covers of the maps $u$.

This proves that in the neighborhood of the point $z_0$ the map $\eta$ has to vanish. As remarked above $\eta$ solves also an equation of the
type $\partial_s\eta+\tilde{J}(z,\eta)\partial_t\eta+Y(z)\cdot\eta=0$. This equation is linear in $\eta$, thus unique continuation
(see \cite[Proposition 3.1]{Floer_Hofer_Salamon_Transversality_in_elliptic_Morse_theory_for_the_symplectic_action}) implies
that $\eta=0$ throughout.

In the case $R=0$ the operator $\O(0,\cdot,\cdot)$ reduces to $\overline{\partial}_{j_0,J}$. It remains to prove that the linearization of
$\overline{\partial}_{j_0,J}$ at the solution set is surjective. We assume that the Maslov index $a$
of the maps $u\in\B$ is negative. The monotonicity of $L$ then implies the solution set of $\overline{\partial}_{j_0,J}$ either is empty ($a<0$)
or equals $L$ ($a=0$). I.e.~it remains to prove that the linearization of $\overline{\partial}_{j_0,J}$ at constant holomorphic disks
is surjective. This follows easily with help of the Schwarz reflection principle from the case of constant holomorphic spheres, wich
in turn can be found in the book \cite{Hofer_Abbas_book_preprint}.

This finishes the proof of Claim \ref{claim:surj_lin}.
\end{proof}
\noindent
Claim \ref{claim:surj_lin} implies the theorem by standard reasoning as follows, see e.g.~\cite[Proposition 4.2.5]{Schwarz_Matthias_PhD} as follows.
We denote by $\pi:\R_{>0}\times\Jh\times\B\pf\Jh$ the projection on the second factor.
From the fact that the restriction of the linearization $D_{(R,J,u)}\big|_{\{0\}\times T_J\Jh\times T_uB}$
of the operator $\O$ is surjective it follows that the operator $\sigma_J:(R,u)\mapsto\overline{\partial}_{j_R,J}u+\tilde{\alpha}_R(z)\kappa(H)(t,u)$
is surjective if and only if $J$ is a regular value the projection $\pi|_{\mathscr{M}}:\mathscr{M}\pf\Jh$ restricted
to the universal moduli space $\mathscr{M}=\O^{-1}(0)$.
Thus, by the Sard-Smale theorem, there exists a set $\Jh_{\mathrm{reg}}\subset\Jh$ of second
category with the property that for all $J\in\Jh_{\mathrm{reg}}$ the operator
$\sigma_J:(R,u)\mapsto\overline{\partial}_{j_R,J}u+\tilde{\alpha}_R(z)\kappa(H)(t,u)$ has a surjective linearization at solutions of $\sigma_J(R,u)=0$.
Thus, for $J\in\Jh_{\mathrm{reg}}$ the moduli space $\M(J,H;a)$.

The dimension formula is derived from the fact that for fixed $R$ the space of maps $u$ has dimension $\dim L+\Mas(u)$. The parameter $R$ adds $+1$ to
the dimension.
\end{proof}

\begin{Rmk}
We proved (and will need) theorem \ref{thm:detailed_transversality} only in the cases $a\leq0$. For $a>0$ we would have to include
a proof of transversality for non-constant holomorphic disks on $L$.
\end{Rmk}
\noindent
For $q,p\in\Crit(f)$ we define the moduli space $\M^{\varphi\circ\rho}(q,p;a)$ to be the set of quadruples $(R,\gamma_-,u,\gamma_+)$, where
\beq\label{eqn:def_M(R,H;a)_start}
R\geq0,\quad(R,u)\in\M(J,H;a),\quad \gamma_-:(-\infty,0]\pf L,\quad \gamma_+:[0,+\infty)\pf L
\eeq
satisfying
\begin{gather}
\dot{\gamma}_\pm(t)+\nabla^g f\circ\gamma_\pm(t)=0\,,\\[1.5ex]
\label{eqn:def_M(R,H;a)_end}
\gamma_-(-\infty)=q,\quad\gamma_-(0)=u(-\zeta),\quad u(+\zeta)=\gamma_+(0),\quad\gamma_+(+\infty)=p\;.
\end{gather}
The function $f:L\rightarrow\R$ is a Morse function on $L$ and $\nabla^{g}$ is the gradient
\wrt a Riemannian metric $g$ on $L$. The moduli space $\M(J,H;a)$ is defined in definition \ref{def:moduli_space_for_perturbed_hol_disks}.

\begin{Thm}\label{thm:tranversality_for_M^phi_o_rho}
If $a\leq0$ then for a generic choice of data the moduli spaces $\M^{\varphi\circ\rho}(q,p;a)$ are smooth manifolds and
\beq
\dim\M^{\varphi\circ\rho}(q,p;a)=\Morse(q)-\Morse(p)+a+1\;.
\eeq
We denote by $\M^{\varphi\circ\rho}(q,p;a)_{[d]}$ the union of the $d$-dimensional components.
\end{Thm}

\begin{Rmk}
We note that the dimension formula for the moduli spaces $\M^{\varphi\circ\rho}(q,p;a)$ is absolute, i.e.~not modulo the minimal Maslov
number. Equivalently, we could have defined $\M^{\varphi\circ\rho}(q,p;a)=\M(J,H;a)\pitchfork W^u(q)\pitchfork W^s(p)$.
\end{Rmk}

\begin{proof}
Transversality for the moduli space $\M(J,H;a)$, $a\leq0$, has been proved in theorem \ref{thm:detailed_transversality}. These spaces carry natural
evaluation maps
\bea
\ev_\pm:\M(J,H;a)&\pf L\\
(R,u)&\mapsto u(\pm\zeta)\;.
\eea
We need to prove that for generic data the evaluation map $\mathrm{Ev}:=\ev_-\times\ev_+$ is transverse
to $W^u(q)\times W^s(p)\subset L\times L$. For this it suffices to prove that the maps
$\ev_\pm:\R_{\geq0}\times\Jh\times\B\pf L$ are submersions when restricted to solutions of $\O(R,J,u)=0$, see for instance
\cite[Theorem 3.4.1]{McDuff_Salamon_J_holomorphic_curves_and_symplectic_topology}. For the relevant notation we refer to the proof of theorem
\ref{thm:detailed_transversality} on page \pageref{thm:detailed_transversality}.

The maps $\ev_\pm$ being submersions at a solution $(R,J,u)$  amounts to showing that for every pair of tangent vectors
$(v_-,v_+)\in T_{u(-\zeta)}L\times T_{u(+\zeta)}L$ we can find $(Y,\xi)\in T_J\Jh\times T_u\B$ such that
\beq
\xi(\pm\zeta)=v_\pm,\quad D_u^J\xi+Y(z,u)\circ\Big(du-\tilde{\alpha}_R(z)dt\otimes X_H(t,u(z))\Big)\circ j_R=0\;.
\eeq
This can be achieved in the same way as in the detailed proof of Lemma 3.4.4 in \cite{McDuff_Salamon_J_holomorphic_curves_and_symplectic_topology}.
\end{proof}

\noindent
Now we come to the core part of this section, namely the suitable compactness statement for the moduli space $\M^{\varphi\circ\rho}(q,p;a=0)_{[d]}$. The
proof has a completely different flavor than the previous compactness proofs since we cannot argue by solely considering the Fredholm index.

\begin{Thm}\label{thm:compactness_for_isomorphism_hard_case}
If the Morse indices $k=\Morse(q)$ and $l=\Morse(p)$ satisfy $n-N_L+2\leq l$ and $k\leq N_L-2$,
then the moduli space $\M^{\varphi\circ\rho}(q,p;a=0)_{[d]}$ is compact up to breaking.
In particular, if $d=0$ the moduli space is compact and if $d=1$ we conclude
\bea
\partial\M^{\varphi\circ\rho}(q,p;0)_{[1]}&=\big\{(R,\gamma_-,u,\gamma_+)\mid R=0\big\}\\[1.5ex]
                                      &\cup\bigcup_{q'\in\Crit(f)}\Mh(q,q';f,g)_{[0]}\times\M^{\varphi\circ\rho}(q',p)_{[0]}\\[1.5ex]
                                      &\cup\bigcup_{p'\in\Crit(f)}\M^{\varphi\circ\rho}(q,p')_{[0]}\times\Mh(p',p;f,g)_{[0]}\\[1.5ex]
                                      &\cup\bigcup_{x\in\P_L(H)}\M^{\varphi}(q,x)\times\M^{\rho}(x,q)
\eea
\end{Thm}

\begin{Rmk}
Since $\M^{\varphi\circ\rho}(q,p;a=0)_{[d]}\not=\emptyset$ implies $\Morse(q)-\Morse(p)+1=k-l+1\geq0$ by theorem \ref{thm:tranversality_for_M^phi_o_rho}
the above theorem is only non-empty for $2N_L\geq n+3$.
\end{Rmk}

\begin{proof}
The computation of the following uniform energy bounds on the energy of solutions $u$, where
$(R,\gamma_-,u,\gamma_+)\in\M^{\varphi\circ\rho}(q,p;0)_{[d]}$, is included in appendix \ref{appendix:energy_estimates}. The following holds,
\beq
E(u)\leq ||H||\,,
\eeq
where $||H||$ is the Hofer norm of $H$.
In particular, we know that sequences $(R_n,\gamma_-^{(n)},u_n,\gamma_+^{(n)})$ converge
modulo breaking and bubbling-off. We claim, that no bubbling occurs given $n-N_L+2\leq l$ and $k\leq N_L-2$,
where $k$ and $l$ are the Morse indices of $q$ and $p$.
We will prove this in two steps. In the first step, we will describe what type of bubbling can occur. In the second step, we then conclude that
this will (generically) not affect the moduli spaces $\M^{\varphi\circ\rho}(q,p;0)_{[d]}$.
If bubbling-off of holomorphic disks and spheres is excluded the statement of the theorem follows from the standard gluing constructions
in Morse and Floer theory.\\[1ex]
%
\textbf{Step 1}:\quad
%
Let us give an auxiliary definition which we will use throughout the proof. For $R\geq0$, $a\in\R$ we set
\beq
\M(R,H;a):=\{u\mid (R,u)\in\M(J,H;a)\}\,,
\eeq
where $\M(J,H;a)$ is defined in definition \ref{def:moduli_space_for_perturbed_hol_disks}.
If bubbling occurs for sequences $(R_n,\gamma_-^{(n)},u_n,\gamma_+^{(n)})$ in $\M^{\varphi\circ\rho}(q,p;0)$ then we can distinguish two cases.\\[0.5ex]
%
%
\textbf{(A)}\quad $\displaystyle R_n\pf R_\infty\in[0,\infty)$\;:\quad We start with the case that exactly one holomorphic disk bubbles off:
$u_n\rightharpoonup (u,d)$, where $u$ is the remaining solution and $d$ is the holomorphic disk. Since the homotopy class is preserved in the limit we
conclude $\om(u)+\om(d)=0$. Since $d$ is holomorphic, $\om(d)>0$ and $a:=\Mas([d])>0$. By definition, $u\in\M(R_\infty,H;-a)$, that is
$(R_\infty,u)\in\M(J,H;-a)$.

If $u=\mathrm{const}$, we conclude $\om(u)=0$ and thus $\om(d)=0$ and thus $d=\mathrm{const}$.
Therefore, no bubbling occurs in this case. This corresponds to the case $R_\infty=0$, where the set $\M(0,H;0)\cong L$
contains exactly the constant maps.

If we assume that the map $u\in\M(R_\infty,H;-a)$ is non-constant (in particular, $R_\infty\not=0$)
we can distinguish three different possibilities how $u$ and $d$ connect,

\begin{enumerate}
\item at $-\infty$: there exists $z_1\in \partial\D^2$ such that $u(-\infty)=d(z_1)$,
\item at $+\infty$: there exists $z_2\in \partial\D^2$ such that $u(+\infty)=d(z_2)$,
\item in between: there exists $z_3\in \partial\D^2$ and $(s,t)\in\R\times\{0,1\}$ such that $u(s,t)=d(z_3)$.
\end{enumerate}
In \textit{step 2} we will argue why this does not harm the moduli spaces $\M^{\varphi\circ\rho}(q,p;0)$. Furthermore, the case of
multiple bubbling is mentioned below.\\[1ex]
%
%
\textbf{(B)}\quad $\displaystyle R_n\pf\infty:$\quad In this case $u_n$ breaks. Again we assume that the
strip breaks once and one holomorphic disk bubbles off: $u_n\rightharpoonup (V_1,V_2,d)$. That is, $V_1\in\M(H;x)$ and $V_2\in\M(x;H)$
for some $x\in\P_L(H)$, (the moduli spaces $\M(H;x)$ and $\M(x;H)$ appear in definition \ref{def:moduli_space_of_cut_off_holo_strips} and
\ref{def:moduli_space_of_cut_off_holo_strips_II}).

As before $\om(V_1\#V_2)+\om(d)=0$ holds. Again, if both maps $V_1$ and $V_2$ are constant then so is the holomorphic
disk $d$ and we conclude that no bubbling occurred. If either $V_1$ or $V_2$ is non-constant so is the other because they converge to the same element
$x\in\P_L(H)$ and are holomorphic near $\pm\infty$.

We are left with the case that all maps $V_1$, $V_2$ and $d$ are non-constant and again we distinguish the places where the bubbling occurred:

\begin{enumerate}\setcounter{enumi}{3}
\item at $-\infty$: there exists $z_4\in \partial\D^2$ such that $V_1(-\infty)=d(z_4)$.
\item at $+\infty$: there exists $z_5\in \partial\D^2$ such that $V_2(+\infty)=d(z_5)$.
\item in between: there exists $z_6\in \partial\D^2$ and $(s,t)\in\R\times\{0,1\}$ such that $V_1(s,t)=d(z_6)$ or $V_2(s,t)=d(z_6)$.
\item at $x\in\P_L(H)$: there exists $z_7\in \partial\D^2$ and $t\in\{0,1\}$ such that $x(t)=d(z_7)$.
\end{enumerate}
What we did not mention in the cases (4) -- (7) is that multiple breaking might occur, i.e.~$u_n\rightharpoonup (V_1,V_2,\ldots,V_r,d)$. The maps
$V_2,\ldots,V_{r-1}$ are perturbed holomorphic strips, i.e. elements in moduli spaces $\M_L(y,z;J,H)$, which define the boundary operator in Floer
homology. For this case we refer the reader to the end of \textit{step 2}.

Before we proceed to \textit{step 2} we note that cases (1) -- (7) are in fact the only cases which we have to consider, that is,
\textit{no multiple bubbling can occur}. Indeed, in all above cases the bubble $d$ has to satisfy $\Mas([d])=N_L$ for the following reason.
Recall that the remaining solutions $u$ respectively $V_1,V_2$ are non-constant elements of a smooth moduli
space of dimension $n+\Mas(-[d])+1$, for instance $(R_\infty,u)\in\M(J,H;-a)$, where $a=\Mas([d])$.
Since $N_L\geq\frac{n}{2}+3$ (otherwise the statement of the theorem is empty) the dimension of these moduli
space satisfies $\dim=n-k\cdot N_L+1\leq n-k(\frac{n+3}{2})+1$. In particular, for $k\geq2$ these spaces are empty due to transversality.
We are left with either $k=0$, in which case no bubbling occurs, or $k=1$ which corresponds to $a=\Mas([d])=N_L$.

This shows that multiple bubbling-off is not possible. We conclude that either a disk or a sphere bubbles off and has Maslov index $N_L$.\\[1ex]
We conclude \textit{step 1} by the following two remarks.\\
$\bullet$\quad If $N_L\geq n+2$ not even $k=1$ is admissible and we are immediately done. In particular, Oh's result that $\HF_*(L)\cong\H_*(L)$ for
$N_L\geq n+2$ follows immediately.\\
$\bullet$\quad The case $R_\infty=0$ is completely settled, that is, for $R_\infty=0$ no bubbling occurs (see \textbf{(A)}).\\[1ex]
%
\textbf{Step 2}:\quad
%
We discuss cases (1) -- (7).

\noindent
We start with case (1): $u_n\rightharpoonup (u,d)$ and
$\exists z_1\in \partial\D^2$ such that $u(-\infty)=d(z_1)$. The map $u$ is a non-constant element in the space $\M(R_\infty,H;-N_L)$
with $R_\infty\in(0,\infty)$. In particular, the moduli space $\M(J,H;-N_L)$ is non-empty, furthermore
according, to theorem \ref{thm:detailed_transversality}, $\dim\M(J,H;-N_L)=n-N_L+1$.

We recall $\M^{\varphi\circ\rho}(q,p;0)=\M(J,H;0)\pitchfork W^u(q)\pitchfork W^s(p)$. In particular, if case (1) occurs then
$\M(J,H;-N_L)\pitchfork W^s(p)\not=\emptyset$, i.e.~the following type of configuration exists (see figure \ref{fig:bad_bubbling}):
$(\gamma_-,d,u,\gamma_+)$, where $\gamma_\pm$ are gradient flow half-lines as before and $(u,d)$ are the limit objects.
All these objects have to be connect to each other as follows:
$\gamma_-(-\infty)=q$, $\gamma_-(0)=d(z_1')$, $d(z_1)=u(-\infty)$, $u(+\infty)=\gamma_+(0)$, $\gamma_+(+\infty)=p$.
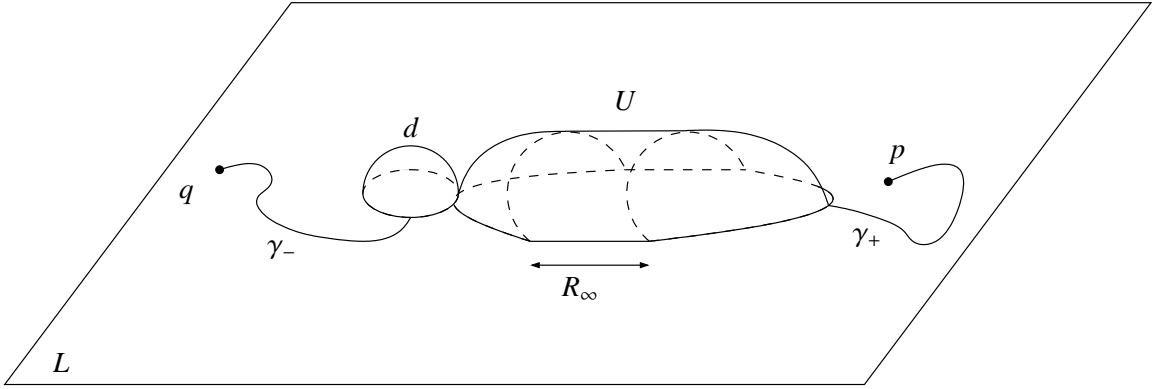
\begin{figure}[ht]
\begin{center}
\input{./bad_bubbling.pstex_t}
\caption{An bubbling configuration possibly spoiling the compactness of $\M^{\varphi\circ\rho}(q,p;0)$.}
\label{fig:bad_bubbling}
\end{center}
\end{figure}
We recall that the Morse indices $k=\Morse(q)$ and $l=\Morse(p)$ satisfy by assumption $n-N_L+2\leq l$ and $k\leq N_L-2$. Thus, the dimension of
$\M(J,H;-N_L)\pitchfork W^s(p)$ is
\bea
\dim\M(J,H;-N_L)\pitchfork W^s(p)&=n-N_L+1-\Morse(p)\\
                                &\leq n-N_L+1-(n-N_L+2)\\
                                &=-1
\eea
This implies that the intersection $\M(J,H;-N_L)\pitchfork W^s(p)=\emptyset$ by transversality. In other words generically
case (1) does not occur.
This shows that case (1) does not affect the moduli spaces $\M^{\varphi\circ\rho}(q,p;0)$ as long as $n-N_L+2\leq l$ and $k\leq N_L-2$ holds,
where $k=\Morse(q)$ and $l=\Morse(p)$.

Before we handle the other cases we want to make the following remark.
Although figure \ref{fig:bad_bubbling} looks similar to figure \ref{fig:PSS_strip_1_bubbling} on page \pageref{fig:PSS_strip_1_bubbling} there is a major
difference. Namely, in the present configuration no chord $x\in\P_L(H)$ is involved. Moreover, in the limit $R_n\pf\infty$ the sequence of
strips $(u_n)$ breaks and there is no control about which chord $x\in\P_L(H)$ appears in the limit. In particular, the arguments in sections
\ref{section:construction_of_varphi}, \ref{section:construction_of_rho} or \ref{section:rho_o_varphi=id}, where we imposed bounds on the
Maslov indices of $x\in\P_L(H)$, are not applicable.\\[1ex]
The same type of argument as for bubbling case (1) settles the other cases.
In case (2) we need to consider the intersection $\M(J,H;-N_L)\pitchfork W^u(q)$ which has
dimension
\bea
\dim\M(J,H;-N_L)\pitchfork W^u(q)&=n-N_L+1-(n-\Morse(q))\\
                                &\leq -N_L+1+(N_L-2)=-1\;.
\eea
For case (3) we can apply the argument for case (1) or (2).\\[1ex]
It remains to deal with cases (4) -- (7). We apply the same type of argument to the space
\beq
\big\{(V_1,V_2)\in\M(H;x)\times\M(x;H)\mid \Mas(V_1\#V_2)=-N_L\big\}\,,
\eeq
which replaces $\M(R_\infty,H;-N_L)$. It is smooth since by assumption the space $\M(H;x)$ and
$\M(x;H)$ are smooth. Furthermore, it has dimension $n-N_L$ and has to has a non-empty intersection with
$W^u(q)$ in case (4) and with $W^s(p)$ in case (5). The same dimension considerations apply. Cases (6) and (7) are then treated as case (3) is.

At the end of \textit{step 1} we mentioned the case of multiple breaking $u_n\rightharpoonup (V_1,V_2,\ldots,V_r,d)$. \textit{Step 2} is easily adapted,
namely the space $\big\{(V_1,V_2)\in\M(H;x)\times\M(x;H)\mid \Mas(V_1\#V_2)=-N_L\big\}$,
is replaced by the set containing $r$-tuples $(V_1,V_2,\ldots,V_{r-1},V_r)$ in the set
\beq
\M(H;x_0)\times\M_L(x_0,x_1;H,J)\times\ldots\times\M_L(x_{r-1},x_r;H,J)\times\M(x_r;H)
\eeq
with the property $\Mas(V_1\#V_2\#\ldots\#V_r)=-N_L$. Now we can argue exactly as above since the set of these $r$-tuples again has dimension
$n-N_L$.\\[1ex]
This concludes the proof of theorem \ref{thm:compactness_for_isomorphism_hard_case}.
\end{proof}
\noindent
Now we are in the same position as in the preceding section \ref{section:rho_o_varphi=id} where we proved
$\rho_k\circ\varphi_k=\id_{\H^{n-k}(L;\Z/2)}$ with help of a
chain homotopy $\Theta^{\rho\circ\varphi}$. We follow the same scheme and define another chain homotopy $\Theta^{\varphi\circ\rho}$.

\begin{Def}
For $n-N_L+1\leq k\leq N_L-2$ we define on generators the map
\bea
\Theta^{\varphi\circ\rho}_k:\CM^{n-k}(L;\Z/2)&\pf\CM^{n-k-1}(L;\Z/2)\\
                                p\quad&\mapsto\quad\sum_{
                                \substack{q\in\Crit(f)\\\Morse(q)=\Morse(p)-1}}
                                \#_2\M^{\varphi\circ\rho}(q,p;0)_{[0]}\cdot q
\eea
and extend it linearly.
\end{Def}

\begin{Rmk}
Since $n-N_L+2\leq k\leq N_L-2$ we conclude $\Morse(p)=n-k\geq n-N_L+2$ and $\Morse(q)=n-k-1\leq n-(n-N_L+1)-1=N_L-2$. Thus the assumptions
of theorem \ref{thm:compactness_for_isomorphism_hard_case} are satisfied and $\Theta^{\varphi\circ\rho}_k$ is well-defined.
We note the asymmetry in the inequalities for the Morse index $k$.
\end{Rmk}
\noindent
Before, we prove that $\Theta^{\varphi\circ\rho}_k$ is a chain homotopy we set
\bea
\vartheta_k:\CM^{n-k}(L;\Z/2)&\pf\CM^{n-k}(L;\Z/2)\\
                                q\quad&\mapsto\quad\sum_p\#_2\big\{(R,\gamma_-,u,\gamma_+)\in\M^{\varphi\circ\rho}(q,p;0)_{[1]}\mid R=0\big\}\cdot p
\eea
and note that the set $\big\{(R,\gamma_-,u,\gamma_+)\in\M^{\varphi\circ\rho}(q,p;0)_{[1]}\mid R=0\big\}$ is zero dimensional and compact by the same
argument as in theorem \ref{thm:compactness_for_isomorphism_hard_case}.
The compactness statement in theorem \ref{thm:compactness_for_isomorphism_hard_case} for one dimensional components immediately implies

\begin{Cor}
For $n-N_L+2\leq k\leq N_L-2$ the following equality holds
\beq
\Theta^{\varphi\circ\rho}_{k-1}\circ\delta^L-\delta^L\circ\Theta^{\varphi\circ\rho}_k=\varphi_k\circ\rho_k-\vartheta_k
\eeq
where $\delta^L$ denotes the Morse-differential, and thus, in homology
\beq
\varphi_k\circ\rho_k=\vartheta_k:\H^{n-k}(L;\Z/2)\pf\H^{n-k}(L;\Z/2)\;.
\eeq
\end{Cor}
\noindent
We note that $\Theta^{\varphi\circ\rho}_{k-1}$ is well-defined due to the assumption $n-N_L+2\leq k$.
The map $\vartheta_k:\CM^{n-k}(L;\Z/2)\pf\CM^{n-k}(L;\Z/2)$ actually equals the identity, even on
chain level. Indeed, we count the elements in the set
\beq
\big\{(R,\gamma_-,u,\gamma_+)\in\M^{\varphi\circ\rho}(q,p;0)_{[1]}\mid R=0\big\}\;.
\eeq
Since $R=0$, the map $u$ is a holomorphic disk satisfying $\om(u)=0$, i.e.~$u$ is constant. In particular, the two gradient half trajectories
match up: $\gamma_-(0)=\gamma_+(0)$. In other words, we count single gradient trajectories $\gamma:\R\pf M$ from $q$ to $p$. Moreover, we
are interested in zero-dimensional families \emph{without} dividing out the $\R$-action, that is
\beq
\#_2\big\{(R,\gamma_-,u,\gamma_+)\in\M^{\varphi\circ\rho}(q,p;0)_{[1]}\mid R=0\big\}=\left\{
\begin{aligned}
1\quad&\text{if } q=p\\
0\quad&\text{else}
\end{aligned}\right.\;.
\eeq
This shows that $\vartheta_k=\id_{\CM^{n-k}(L;\Z/2)}$ and concludes the section. In particular, we finally proved the full statement of theorem
\ref{thm:existence_of_Lagrangian_PSS}.
\section{The comparison homomorphisms}\label{section:construction_of_the_comparison_homs}
\noindent
In the last two sections we prove theorem \ref{thm:diagram_commutes}. We start the construction of the homomorphism
\beq
\chi_k:\HF_k(L,\phi_H(L))\pf\HF_{k-n}(H)\;.
\eeq
\noindent
As already mentioned in the introduction this map was independently considered by A.~Abbondandolo and M.~Schwarz in the context of
Floer homology of cotangent bundles and the ring-isomorphism to the homology of the loop space.

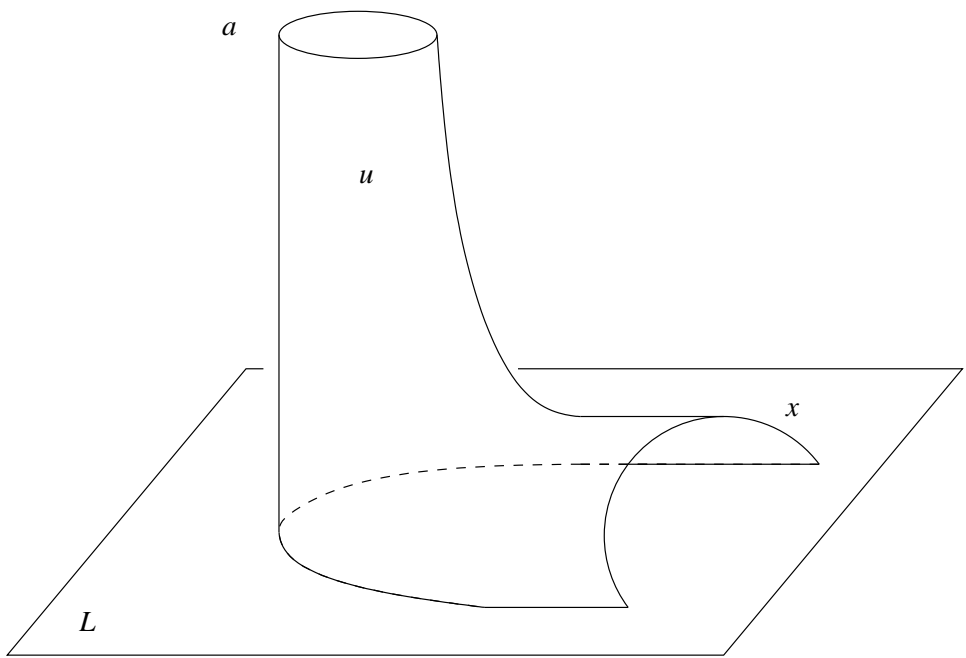
\begin{figure}[ht]
\begin{center}
\input{./chi_pic.pstex_t}
\caption{An element $u\in\M_L^\chi(a,x)$}
\label{fig:chi_pic}
\end{center}
\end{figure}

Let us recall the current setting (cf.~theorem \ref{thm:existence_of_Lagrangian_PSS}). $(M,\om)$ is a closed symplectic manifold and
$L\subset M$ a closed, monotone Lagrangian submanifold with minimal
Maslov number $N_L\geq2$. Furthermore, $H:S^1\times M\pf\R$ is a admissible Hamiltonian function, see definition \ref{def:admissible_Hamiltonian}.
Schematically, $\chi$ is defined by counting maps of the form depicted in figure \ref{fig:chi_pic}.
\noindent
A very nice description of the chimney shaped domain in figure \ref{fig:chi_pic} is due to Abbondandolo-Schwarz in
\cite{Abbo_Schwarz_Notes_on_Floer_homology_and_loop_space_homology} and goes as follows. We set
\beq
\Sigma^\chi:=\R\times[0,1]/\sim\qquad\text{where}\quad(s,0)\sim(s,1)\text{ for }\;s\leq0
\eeq
with the induced conformal structure. In other words, the interior of $\Sigma^\chi$ is a cylinder with a half-line removed.
Furthermore, the conformal structure on the interior is standard. $\Sigma^\chi$ is a Riemann surface with boundary where the conformal structure
at the point $(0,0)$ is induced by the map $z\mapsto\sqrt{z}$.
In \cite{Abbo_Schwarz_Notes_on_Floer_homology_and_loop_space_homology} further pictures and details can be found. A nice feature of this
particular description of the Riemann surface $\Sigma^\chi$ are the global conformal coordinates.

The set $\Sigma^\chi$ is the domain for the maps $u$ defining the homomorphism $\chi$. For $x\in\P_L(H)$ and $a\in\P(H)$ we define
\beq\label{def:moduli_space_for_chi}
\M_L^\chi(a,x):=\left\{u:\Sigma^\chi\pf M\left|\;\;
    \begin{aligned}
    &\partial_su+J(t,u)\big(\partial_tu-X_H(t,u)\big)=0\\
    &u(s,0),u(s,1)\in L\;\;\;\forall s\geq0\\
    &u(-\infty)=a,\;\;u(+\infty)=x
    \end{aligned}\;
\right.
\right\}\;.
\eeq
\noindent
We recall that as before $J(t,\cdot)$ is a smooth family of compatible almost complex structures and $H:S^1\times M\pf\R$ is a Hamiltonian
function on $M$.
The next proposition is taken from \cite[proposition 3.3]{Abbo_Schwarz_Notes_on_Floer_homology_and_loop_space_homology}.
\begin{Prop}\label{prop:transversality_comparison_chi}
Let $H:S^1\times M\pf\R$ be admissible (in particular non-degenerate for both, Hamiltonian and Lagrangian Floer homology) Hamiltonian function.
Then for a generic choice of a family $J(t,\cdot)$ of almost complex structures
the moduli space $\M_L^\chi(a,x)$ is a smooth manifold of dimension
\beq
\dim\M_L^\chi(a,x)=\mu(x)-\CZ(a)-n\quad\mathrm{mod}\;N_L\;.
\eeq
As before we denote by $\M_L^\chi(a,x)_{[d]}$ the union of the $d$-dimensional components.
\end{Prop}
\noindent
In \cite{Abbo_Schwarz_Notes_on_Floer_homology_and_loop_space_homology} the authors include a compactness statement in the same
proposition since they work in cotangent bundles which are exact, moreover, the involved Lagrangian submanifolds are exact, thus bubbling is
not present. In the current situation bubbling is possible but can be ruled out due to monotonicity.
\begin{Prop}
Under the hypothesis of theorem \ref{thm:existence_of_Lagrangian_PSS} the following holds.
\begin{enumerate}
\item The moduli space $\M_L^\chi(a,x)_{[0]}$ is compact and\\[-2ex]
\item the moduli space $\M_L^\chi(a,x)_{[1]}$ is compact up to breaking. Moreover, in this case it can be compactified such that
\bea
\partial\M_L^\chi(a,x)_{[1]}=&\bigcup_{a'\in\P(H)}\Mh(a,a';J,H)_{[0]}\times\M_L^\chi(a',x)_{[0]}\\[1ex]
                         \cup&\bigcup_{x'\in\P_L(H)}\M_L^\chi(a,x')_{[0]}\times\Mh_L(x',x;J,H)_{[0]}
\eea
\end{enumerate}
\end{Prop}

\begin{proof}
Bubbling-off is handled in the same manner as for previous moduli spaces.
Elements $u\in\M_L^\chi(a,x)$ satisfy the same uniform energy estimate as do connecting trajectories (see appendix \ref{appendix:energy_estimates}).
In particular, sequences converge to broken solutions and holomorphic disks resp.~holomorphic spheres. The very same argument as in
the proof of theorem \ref{thm:compactness_Lagrangian_PSS} rules out bubbling for the zero and one-dimensional moduli spaces.
\end{proof}

\noindent We obtain a map
\bea
\chi:\CF_k(L,\phi_H(L))&\pf\CF_{k-n}(H)\\
x\quad&\mapsto\!\!\!\!\!\sum_{\CZ(a)=\mu(x)-n}\!\!\!\!\!\#_2\M_L^\chi(a,x)\cdot a
\eea
which descends to homology
\beq
\chi:\HF_k(L,\phi_H(L))\pf\HF_{k-n}(H)\;.
\eeq

\begin{Rmk}
We recall that we need to reduce the grading of the Hamiltonian Floer homology $\HF_*(H)$ from modulo the minimal Chern number $N_M$ to modulo
the minimal Maslov number $N_L$ due to the dimension formula in proposition \ref{prop:transversality_comparison_chi}.
To avoid this we could work with appropriate Novikov rings.
\end{Rmk}
\noindent The map $\tau:\HF_{k}(H)\pf\HF_{k}(L,\phi_H(L))$ is defined by flipping picture \ref{fig:chi_pic}. More
precisely, we define the surface
\beq
\Sigma^\tau:=\R\times[0,1]/\sim\qquad\text{where}\quad(s,0)\sim(s,1)\text{ for }\;s\geq0
\eeq
which is obtained from $\Sigma^\chi$ by replacing $s$ by $-s$. In particular, all above result remain unchanged but the dimension formula which is
modified accordingly as follows.
\begin{Prop}
For admissible Hamiltonian function $H:S^1\times M\pf\R$ and a generic choice of a family $J(t,\cdot)$ of almost complex structures
the moduli space
\beq\label{def:moduli_space_for_commuting_diagram}
\M_L^\tau(x,a):=\left\{u:\Sigma^\tau\pf M\left|\;\;
    \begin{aligned}
    &\partial_su+J(t,u)\big(\partial_tu-X_H(t,u)\big)=0\\
    &u(s,0),u(s,1)\in L\;\;\forall s\leq0\\
    &u(-\infty)=x,\;\;u(+\infty)=a
    \end{aligned}\;
\right.
\right\}
\eeq
is a smooth manifold of dimension
\beq
\dim\M_L^\tau(a,x)=\CZ(a)-\mu(x)\quad\mathrm{mod}\;N_L\;.
\eeq
Moreover, the zero-dimensional components are compact and the one-dimensional components are compact up to simple breaking.
\end{Prop}

\begin{Rmk}
The dimension formula looks asymmetric in the sense that one would expect a summand $\pm n$ in the formula. That this is not the case is due
to our index convention as described in section \ref{section:recollection_of_Floer_homology}. Indeed, under the change $s\mapsto -s$ the indices
behave as follows: $\mu(x)\mapsto n-\mu(x)$ and $\CZ(a)\mapsto-\CZ(a)$.
\end{Rmk}

\noindent We obtain a map
\bea
\tau:\CF_{k}(H)&\pf\CF_{k}(L,\phi_H(L))\\
a\quad&\mapsto\!\!\!\!\!\sum_{\CZ(a)=\mu(x)}\!\!\!\!\!\#_2\M_L^\tau(x,a)_{[0]}\cdot x
\eea
which descends to homology
\beq
\tau:\HF_{k}(H)\pf\HF_{k}(L,\phi_H(L))\;.
\eeq
\section{The diagram commutes}\label{section:the_diagram_commutes}

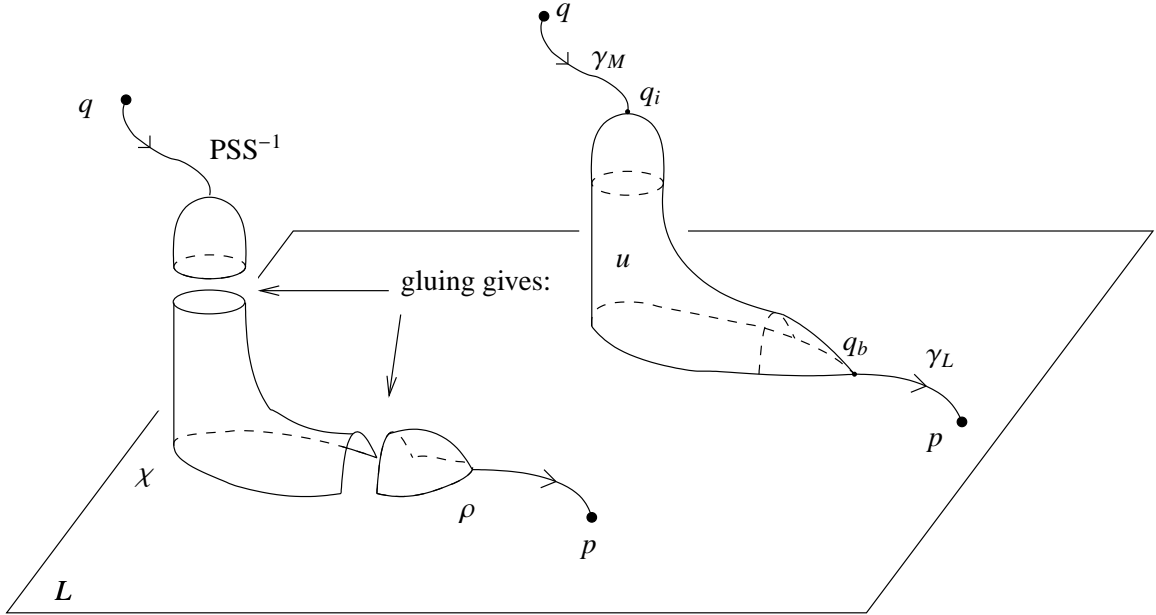
\begin{figure}[ht]
\begin{center}
\input{./diagram_commutes_1.pstex_t}
\caption{$(\gamma_M,u,\gamma_L)\in\M^!(q,p)$}
\label{fig:diagram_commutes_1}
\end{center}
\end{figure}

\noindent
In this section we prove theorem \ref{thm:diagram_commutes}, more precisely we show that the diagrams
\beqn
\xymatrix{
{\HF_k(L,\phi_H(L))}\ar[d]^{\varphi_k}\ar[r]^-{\chi_k}&{\HF_{k-n}(H)}&
             &{\HF_{k}(H)}\ar[r]^-{\tau_k}&{\HF_{k}(L,\phi_H(L))}\ar[d]^{\varphi_k}\\
{\H^{n-k}(L;\Z/2)}\ar@<1ex>[u]^{\rho_k}\ar[r]^{\iota^!}&{\H^{2n-k}(M;\Z/2)}\ar[u]^{\PSS}_\cong&
            &{\H^{n-k}(M;\Z/2)}\ar[u]^{\PSS}_\cong\ar[r]^{\iota^*}&{\H^{n-k}(L;\Z/2)}\ar@<1ex>[u]^{\rho_k} \\
}\eeq
commute. This is again achieved by inspecting suitable cobordisms. Of the four identities, which we have to prove, we will give details
for the following
\beq
\iota^!=\PSS^{-1}\circ\chi_k\circ\rho_k:\H^{n-k}(L;\Z/2)\pf\H^{2n-k}(M;\Z/2)\quad\text{for }\;k\geq n-N_L+2\;.
\eeq
The suitable moduli space $\M^!(q,p)$ will defined below.
%
%
The geometric idea behind this cobordism is depicted in figure \ref{fig:diagram_commutes_1}.

\begin{figure}[ht]
\begin{center}
\input{./comparison.pstex_t}
\caption{The domain $\Xi$}
\label{fig:comparision}
\end{center}
\end{figure}
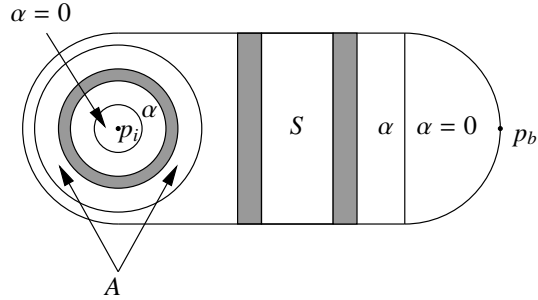

In order to carry out the transversality theory we need to conformally reparameterize the domain
$\Sigma^\chi$ as we did in section \ref{section:varphi_o_rho=id_transversality}.
We again consider a disk-like domain $\Xi$ in the plain, and fix one point $p_i$ in the interior and
another point $p_b$ on the boundary, see figure \ref{fig:comparision}.

Next we fix a rectangular region $S$ of the form $[-1,1]\times[0,1]$
and an annular region $A$ centered at $p_i$ as depicted in figure \ref{fig:comparision}. Then we use a conformal rescaling as in section
\ref{section:varphi_o_rho=id_transversality}, that is, $[-1,1]\times[0,1]$ is mapped onto $[-R,R]\times[0,1]$ via $(s,t)\mapsto(Rs,t)$.
If we use polar coordinates to write for the annulus $A=\{(r,\theta)\mid1\leq r\leq2,\,\theta\in S^1\}$ the conformal rescaling is
given by $(r,\theta)\mapsto (Rr,\theta)$. In both cases $R\geq1$.
In the grey regions in figure $\ref{fig:comparision}$ we interpolate between the conformal rescaling
and the identity map on the other parts of $\Xi$. Pulling back the standard complex structure gives rise to a family of complex structures $j_R$
on the standard disk, which in the limit $R\pf\infty$ degenerates as depicted in figure \ref{fig:diagram_commutes_1}. We note that
$\Xi\setminus\{p_i,p_b\}$ is conformally equivalent to $\Sigma^\chi$, see figure \ref{fig:chi_pic}.

We define a cut-off function $\alpha:\Xi\pf[0,1]$ by setting it to zero as depicted in figure \ref{fig:comparision}, that is,
near the two points $p_i$ and $p_b$. The region in the figure \ref{fig:comparision} containing the symbol $\alpha$
are where the cut-off function $\alpha$ is non-constant. Finally,
on the remaining regions it is equal to $1$. We choose $\alpha$ such that on the annular region $A$ it depends only on the radial coordinate $r$
and on the rectangular region $S=[-1,1]\times[0,1]$ such that it depends only on the $s\in[-1,1]$.
We define
\beq
\tilde{\alpha}_R:=\begin{cases}\,R\alpha &\text{for }0\leq R\leq1\\\,\alpha &\text{for }R\geq1\end{cases}
\eeq

We denote by $\Jh$ the set of families of compatible almost complex structures $J(z,\cdot)$, $z\in\Xi$, which depend on $z$ only inside the closed
rectangular and annular set where the cut-off function is non-constant, see figure \ref{fig:comparision}. Moreover,
on the annular region $A$ and rectangular region $S$ we demand that $J=J_0$
where $J_0$ is as in convention \ref{conv:Floer_generic_J}.
An element $J\in\Jh$ and $j_R$ determine a Cauchy-Riemann operator $\overline{\partial}_{j_R,J}$.

Finally, we choose an admissible (see definition \ref{def:admissible_Hamiltonian}) Hamiltonian function $H$ and define the corresponding
1-form $\kappa(H)$ with values in $TM$ as in definition \ref{def:1_form_hamiltonian} on $\Xi\setminus\{p_i,p_b\}\cong\Sigma^\chi$. We recall
that this set in fact admits global conformal coordinates.

With these definitions at hand
we define the moduli space $\widetilde{\M}^!$ to be the set of pairs $(R,u)$ satisfying
\beq\begin{cases}
\quad u:\Xi\pf M,\quad R\geq0\\
\quad\overline{\partial}_{j_R,J}u(z)+\tilde{\alpha}_R(z)\kappa(H)(z,u(z))=0\\
\quad u(\partial\Xi)\subset L\\
\quad \om(u)=0
\end{cases}\eeq
\noindent
Since $\tilde{\alpha}_R(z)$ vanishes near $\{p_i,p_b\}$ the form $\tilde{\alpha}_R(z)\kappa(H)(z,u(z))$ is well-defined over all of $\Xi$.

\begin{Thm}\label{thm:detailed_transversality_in_comparison}
For an admissible Hamiltonian function $H$ and a generic $J\in\Jh$, the moduli space $\widetilde{\M}^!$ is a smooth manifold of dimension
$\dim \widetilde{\M}^!=\dim L+1$.
\end{Thm}

\begin{proof}
As in the proof of theorem \ref{thm:detailed_transversality} it can be proved that the universal moduli space is a
Banach manifold, cp.~claim \ref{claim:surj_lin}. This uses the fact that the Hamiltonian function is admissible. This relies
on claim \ref{claim:non-vanishing}. The proofs of both claims are unchanged in the present situation.
Again at $R=0$ the equation reduces to constant holomorphic disks.

In the present situation we don't need the full strength of the definition \ref{def:admissible_Hamiltonian} of
admissibility. Indeed in the proof of claim \ref{claim:non-vanishing} we only use that there are no non-constant elements
in $\P_L(H)$.That elements in $\P(H)$ are non-constant is necessary for the transversality theory of the moduli spaces below.
\end{proof}

The moduli space $\M^!(q,p)$ proving the equality $\iota^!=\PSS^{-1}\circ\chi_k\circ\rho_k$, consists out of tuples $(R,\gamma_M,u,\gamma_L)$ where
\beq
R\geq0,\quad\gamma_M:(-\infty,0]\pf M,\quad u\in\widetilde{\M}^!,\quad\gamma_L:[0,+\infty)\pf L
\eeq
satisfying
\begin{gather}
\dot{\gamma}_M(t)+\nabla^{g_M} f_M\circ\gamma_M(t)=0\,,\qquad
\dot{\gamma}_L(t)+\nabla^{g_L} f_L\circ\gamma_L(t)=0\,,\\[1.5ex]
\gamma_M(-\infty)=q,\quad\gamma_M(0)=u(p_i),\quad u(p_b)=\gamma_L(0),\quad\gamma_L(+\infty)=p\,,
\end{gather}
see figure \ref{fig:diagram_commutes_1}.
\begin{Prop}
For a generic choice of data the moduli space $\M^!(q,p)$ is a smooth manifold and
\beq
\dim\M^!(q,p)=\Morse(q;f_M)-\Morse(p;f_L)-n+1\;.
\eeq
\end{Prop}

\begin{proof}
This follows from theorem \ref{thm:detailed_transversality_in_comparison} together with the fact that the evaluation maps
at $p_i$ and $p_b$ are submersions on the universal moduli spaces. In particular, for generically chosen Morse-Smale pairs
$(f_M,g_M)$ and $(f_L,g_L)$ the evaluation maps will be transverse to the unstable and stable manifolds, see the proof
of theorem \ref{thm:tranversality_for_M^phi_o_rho} for some more details.
\end{proof}
\noindent
The compactness issues are exactly the same as in section \ref{section:varphi_o_rho=id_transversality}, cf.~theorem \ref{thm:compactness_for_isomorphism_hard_case}.

\begin{Thm}\label{thm:compactification_diagram_commutes}
If the Morse index $k=\Morse(p)$ satisfies $k\leq N_L-1$, then the moduli space $\M^!(q,p)_{[d]}$
is compact up to breaking. In particular, if $d=0$ the moduli space is compact and if $d=1$ we conclude
\bean
\partial\M^!(q,p)_{[1]}&\stackrel{(1)}{=}\big\{(R,\gamma_M,u,\gamma_L)\mid R=0\big\}\\[1ex]
                        &\stackrel{(2)}{\cup}\bigcup_{q'\in\Crit(f_M)}\Mh(q,q';f_M,g_M)_{[0]}\times\M^!(q',p)_{[0]}\\[1ex]
                        &\stackrel{(3)}{\cup}\bigcup_{p'\in\Crit(f_L)}\M^!(q,p')_{[0]}\times\Mh(p',p;f_L,g_L)_{[0]}\\[1ex]
                        &\stackrel{(4)}{\cup}\bigcup_{a\in\P(H)}\M^{\mathrm{PSS,\,inv}}(q,a)_{[0]}\times\M^{\mathrm{I}}(a,p;0)_{[0]}\\[1ex]
                        &\stackrel{(5)}{\cup}\bigcup_{x\in\P_L(H)}\M^{\mathrm{I\!\:\!I}}(q,x;0)_{[0]}\times\M^\rho(x,p)_{[0]}
\eea
\begin{enumerate}
\item This is the obvious boundary of $\M^!(q,p)_{[1]}$ given by $R=0$.
\item[(2)\&(3)] Either $\gamma_M$ or $\gamma_L$ breaks.
\item[(4)\&(5)] $u:\Sigma^\chi\pf M$ breaks either at $-\infty$ or at $+\infty$.
\end{enumerate}
The moduli spaces $\M^{\mathrm{I}}(a,p;0)$ and $\M^{\mathrm{I\!\:\!I}}(q,x;0)$ are defined below. $\M^{\mathrm{PSS,\,inv}}(q,a)$ is defined
in section \ref{section:construction_of_real_PSS} and $\M^\rho(x,p)$ in section \ref{section:construction_of_rho}.
\end{Thm}

\begin{proof}
This is proved by the same methods as before, see the proof of theorem \ref{thm:compactness_for_isomorphism_hard_case}.
\end{proof}

In case (1) the map $u$ is pseudo-holomorphic since due to $R=0$ the Hamiltonian perturbation equals 0. According to the definition of
$\M^!(q,p)$ we have $\om([u])=0$ and thus $u$ is constant. We conclude
\beqn
\big\{(R,\gamma_M,u,\gamma_L)\mid R=0\big\}=\left\{(\gamma_M,\gamma_L)\left|\;\;
    \begin{aligned}
    &\dot{\gamma}_\blacklozenge(t)+\nabla^{g_\blacklozenge} f_\blacklozenge\circ\gamma_\blacklozenge(t)=0,\quad\blacklozenge\in\{L,M\}\\
    &\gamma_M(-\infty)=q,\;\gamma_L(+\infty)=p\\
    &\gamma_M(0)=\gamma_L(0)
    \end{aligned}\;
\right.
\right\}\;.
\eeq
\noindent
It is a straight-forward exercise in Morse homology to prove that counting zero dimensional components of the right hand side
defines the map $i^!:\H^{n-k}(L)\pf\H^{2n-k}(M)$.\\[1ex]
For any real number $T\in\R$ we define the moduli space $\M^{\mathrm{I}}(a,p;T)$ as follows
\beqn
\left\{(u,\gamma_L)\left|\;\;
    \begin{aligned}
    &u:\Sigma^\chi\pf M\,,\quad\partial_su+J(s,t,u)\big(\partial_tu-\beta(-s-T)\cdot X_H(t,u)\big)=0\,,\quad u(-\infty)=a\\
    &\gamma_L:[0,\infty)\pf L\,,\quad\dot{\gamma}_L(t)+\nabla^{g_L} f_L\circ\gamma_L(t)=0\,,\quad\gamma_L(-\infty)=p\\
    &E(u)<+\infty\,,\quad u(+\infty)=\gamma_L(0)\\
    \end{aligned}\;
\right.
\right\}\;.
\eeq
The cut-off function $\beta$ satisfies $\beta(s)=0$ for $s\leq0$ and $\beta(s)=1$ for $s\geq1$. As before the finite energy condition
$E(u)<+\infty$ together with the cut-off of the Hamiltonian term guarantees the existence of an extension of $u$ to $u(+\infty)$. For $T=0$
the moduli spaces in theorem \ref{thm:compactification_diagram_commutes} is obtained.\\[1ex]
The moduli space $\M^{\mathrm{I\!\:\!I}}(q,x;T)$ is defined as follows
\beqn
\left\{(\gamma_M,u)\left|\;\;
    \begin{aligned}
    &\gamma_M:(-\infty,0]\pf M\,,\quad\dot{\gamma}_M(t)+\nabla^{g_M} f_M\circ\gamma_M(t)=0\,,\quad\gamma_M(-\infty)=q\\
    &u:\Sigma^\chi\pf M\,,\quad\partial_su+J(s,t,u)\big(\partial_tu-\beta(s+T)\cdot X_H(t,u)\big)=0\,,\quad u(+\infty)=x\\
    &E(u)<+\infty\,,\quad\gamma_M(0)=u(-\infty)
    \end{aligned}\;
\right.
\right\}\;.
\eeq
We note that in the latter definition the cut-off function $\beta(s+T)$ is replaced by $\beta(-s-T)$ and $\gamma_M$ is a gradient flow line in $M$.
\noindent
With the following notations
\bea
\Theta:\CM^{n-k}(f_L,g_L)&\pf\CM^{2n-k-1}(f_M,g_M)\\
p\quad&\mapsto \sum_{q} \#_2\M_L^!(q,p)_{[0]}\cdot q
\eea
\bea
\mathrm{I}:\CM^{n-k}(f_L,g_L)&\pf\CF_{k-n}(H)&\quad\mathrm{I\!\:\!I}:\CF_k(L,\phi_H(L))&\pf\CM^{2n-k}(f_M,g_M)\\
p\quad&\mapsto \sum_{a} \#_2\M^{\mathrm{I}}(a,p)_{[0]}\cdot a&x\quad&\mapsto \sum_{q} \#_2\M^{\mathrm{I\!\:\!I}}(q,x)_{[0]}\cdot q
\eea
theorem \ref{thm:compactification_diagram_commutes} implies
\beq
i^!-(\PSS^{-1}\circ\mathrm{I}+\mathrm{I\!\:\!I}\circ\rho_k)=\delta^M\circ\Theta+\Theta\circ\delta^L\;:\,\CM^{n-k}(f_L,g_L)\pf\CM^{2n-k}(f_M,g_M)\;.
\eeq
\noindent
Thus, in homology the identity
\beq
i^!=\PSS^{-1}\circ\mathrm{I}+\mathrm{I\!\:\!I}\circ\rho_k\;:\,\H^{n-k}(f_L,g_L)\pf\H^{2n-k}(f_M,g_M)
\eeq
holds. We need to prove that in homology the identity $\PSS^{-1}\circ\mathrm{I}+\mathrm{I\!\:\!I}\circ\rho_k=\PSS^{-1}\circ\chi_k\circ\rho_k$ holds.
This is again achieved by two cobordisms.
Indeed, if we define the maps $\mathrm{I^T}$ and $\mathrm{I\!\:\!I}^T$ with help of the moduli spaces $\M^{\mathrm{I}}(a,p;T)$ and
$\M^{\mathrm{I\!\:\!I}}(q,x;T)$ then for $T\rightarrow\infty$ the maps $\mathrm{I^T}$ and $\mathrm{I\!\:\!I}^T$
degenerate into $\chi_k\circ\rho_k$ and $\PSS^{-1}\circ\chi_k$, respectively.
The cobordisms are the appropriate compactification of $\cup_{T\geq0}\{T\}\times\M^{\mathrm{I}}(a,p;T)$ and
$\cup_{T\geq0}\{T\}\times\M^{\mathrm{I\!\:\!I}}(a,p;T)$, respectively.
That the map $\PSS^{-1}\circ\chi_k\circ\rho_k$ is recovered as a sum of two cobordism arguments is due to the fact that the corresponding
gluing procedure in the compactifications reflect either $\PSS^{-1}\circ\big(\chi_k\circ\rho_k\big)$ or $\big(\PSS^{-1}\circ\chi_k\big)\circ\rho_k$.
This proves $i^!=\PSS^{-1}\circ\chi_k\circ\rho_k$. The other three equalities from theorem \ref{thm:diagram_commutes} are proved analogously.
We leave the details to the reader.
\begin{Rmk}
To achieve transversality we again need to conformally reparametrize $\Sigma^\chi\sim\Xi\setminus\{p_i,p_b\}$. For $\M^{\mathrm{I}}$ the suitable
model is a disk with one interior puncture $\Xi\setminus\{p_i\}$ and for $\M^{\mathrm{I\!\:\!I}}$ the suitable
model is a disk with one boundary puncture $\Xi\setminus\{p_b\}$. To prove transversality one again needs to check claims
\ref{claim:surj_lin} and \ref{claim:non-vanishing}. In the case $\M^{\mathrm{I}}$ this can be copied verbatim. For $\M^{\mathrm{I\!\:\!I}}$
the argument in claim \ref{claim:non-vanishing} is the same only that now it is used that there is no constant element in $\P(H)$
as opposed to $\P_L(H)$.
\end{Rmk}
\begin{figure}[ht]
\begin{center}
\input{./diagram_commutes_2.pstex_t}
\caption{}
\label{fig:diagram_commutes_2}
\end{center}
\end{figure}
\begin{Rmk}
Another way of proving the identity $i^!=\PSS^{-1}\circ\chi_k\circ\rho_k$ is by noting that $\chi_k\circ\rho_k$ is cobordant to
counting perturbed half-cylinders with boundary on $L$ plus a gradient half-trajectory, see figure \ref{fig:diagram_commutes_2}. That
counting these half-cylinders gives rise to $\PSS\circ i^!$ is proved in \cite[theorem 3.1]{Albers_On_the_extrinsic_topology_of_Lagrangian_submanifolds}.
\end{Rmk}

\appendix
\section{Energy estimates}\label{appendix:energy_estimates}
%
\noindent
For the convenience of the reader we compile in this appendix some standard computations in Floer homology which
lead to uniform energy bounds for solutions in various moduli spaces.
These moduli spaces depend on some boundary data such as $x\in\P_L(H)$ or a homotopy class $\alpha\in\pi_2(M,L)$ etc.. The term \textit{uniform
energy bounds} indicates that the energy for all solutions in the moduli space can be bounded in quantities only involving the boundary data,
e.g.~the action value $\A_H(x,d_x)$ or the $\om$-integral $\om(\alpha)$.
We recall that the \textit{energy} $E(u)$ of a map $u:\R\times[0,1]\pf M$ is
\beq
E(u)=\!\int_{-\infty}^\infty\int_0^1|\partial_su|^2\,dt\,ds\;.
\eeq

\noindent
For convenience we set
\beq
\sup_MH:=\int_0^1\sup_MH(t,\cdot)\,dt\qquad\text{and}\qquad\inf_MH:=\int_0^1\inf_MH(t,\cdot)\,dt\;.
\eeq
In particular, in this notation the Hofer norm $||H||$ of a Hamiltonian function $H:S^1\times M\pf\R$ reads $\displaystyle||H||=\sup_MH-\inf_MH$.
We begin with energy estimates for solutions in the moduli spaces $\M(H;x)$ and $\M(x;H)$ (cf.~definitions
\ref{def:moduli_space_of_cut_off_holo_strips} and \ref{def:moduli_space_of_cut_off_holo_strips_II}).
\begin{Lemma}\label{lemma:energy_estimate_M(x;H)_and_M(H;x)}
For a solution $u\in\M(H;x)$ the following inequality holds:
\beq
0\leq E(u)\leq\A_H(x,u)+\sup_MH\,,
\eeq
cf.~equation \eqref{eqn:action_functional} for the definition of the action functional. Furthermore, for an element $u\in\M(x;H)$
\beq
0\leq E(u)\leq-\A_H(x,-u)-\inf_MH\,,
\eeq
where $-u$ denotes the map $(s,t)\mapsto u(-s,t)$.
\end{Lemma}

\begin{proof}
We start with an element $u\in\M(H;x)$ and denote by $\left<\cdot,\cdot\right>:=\om(\cdot,J(s,t,\cdot)\cdot)$ the $(s,t)$-dependent
inner product induced by family of compatible almost complex structures $J$.

Before we perform the computation we recall that $u$ solves Floer equation: $\partial_su+J(s,t,u)\big(\partial_tu-\beta(s)X_H(t,u)\big)=0$,
where $\beta:\R\rightarrow[0,1]$ is a smooth cut-off function satisfying $\beta(s)=0$ for $s\leq0$ and $\beta(s)=1$ for $s\geq1$.
We recall the sign convention $\om(X_H,\cdot)=-dH$ for the Hamiltonian vector field $X_H$. Moreover, we point out
that in the computations below, the $(s,t)$-dependence of the almost complex structure and of the metric drops out, see for instance line 3 in
the following computation.

\bean
E(u)=&\;\int_{-\infty}^{+\infty}\int_0^1|\partial_s u|^2\,dt\,ds\\[2ex]
=&\;\int_{-\infty}^{+\infty}\int_0^1
            \left<\,\partial_s u,\Big(\underbrace{-J(s,t,u)\,\big(\partial_t u-\beta(s)X_H(t,u)\big)}_{=\partial_su}\Big)\,\right>\,dt\,ds\\[2ex]
=&\;\int_{-\infty}^{+\infty}\int_0^1\Big[\om(\partial_s u,\partial_t u)
    -\beta(s)\cdot\underbrace{\om\big(\partial_su,X_H(t,u)\big)}_{=+dH(t,\partial_su)}\Big]\,dt\,ds\\[2ex]
=&\;\int_{-\infty}^{+\infty}\int_0^1\Big[\om(\partial_s u,\partial_t u)-\beta(s)\cdot dH(t,\partial_su)\Big]\,dt\,ds\\[2ex]
=&\;\om(u)-\int_{-\infty}^{+\infty}\int_0^1\beta(s)\cdot dH(t,\partial_su)\,dt\,ds\\[2ex]
=&\;\om(u)+\int_{-\infty}^{+\infty}\int_0^1\beta'(s)\cdot H(t,u)\,dt\,ds
        -\int_{-\infty}^{+\infty}\int_0^1\frac{d}{ds}\Big[\beta(s)\cdot H(t,u)\Big]\,dt\,ds\\[2ex]
=&\;\om(u)+\int_{-\infty}^{+\infty}\beta'(s)\int_0^1H(t,u)\,dt\,ds\\
        &\qquad\qquad-\int_0^1\Big[\underbrace{\beta(+\infty)}_{=1}\cdot H(t,\underbrace{u(+\infty)}_{=x(t)})-
                            \underbrace{\beta(-\infty)}_{=0}\cdot H(t,u(-\infty))\Big]\,dt\\[2ex]
=&\;\underbrace{\om(u)-\int_0^1H(t,x(t))\,dt}_{=\A_H(x,u)}+\int_{-\infty}^{+\infty}\beta'(s)\int_0^1H(t,u)\,dt\,ds\\[2ex]
\leq&\;\A_H(x,u)+\underbrace{\int_{-\infty}^{+\infty}\beta'(s)\,ds}_{=1}\;\cdot\int_0^1\sup_MH(t,\cdot)\,dt\\[2ex]
=&\;\A_H(x,u)+\sup_MH
\eea
The computation for an element $u\in\M(x;H)$ is the same up to the fact, that $\beta(s)$ becomes replaced
by $\beta(-s)$. Furthermore, the map $u$ is not a homotopy from a constant path to $x$ but the other way round.
This leads to the following changes in the above computation.
\bean
E(u)=&\;\int_{-\infty}^{+\infty}\int_0^1|\partial_s u|^2\,dt\,ds\\[2ex]
=&\;\om(u)-\int_{-\infty}^{+\infty}\int_0^1\beta(-s)\cdot dH(t,\partial_su)\,dt\,ds\\[2ex]
=&\;-\om(-u)-\int_{-\infty}^{+\infty}\int_0^1\beta'(-s)\cdot H(t,u)\,dt\,ds
        -\int_{-\infty}^{+\infty}\int_0^1\frac{d}{ds}\Big[\beta(-s)\cdot H(t,u)\Big]\,dt\,ds\\[2ex]
=&\;-\om(-u)-\int_{-\infty}^{+\infty}\beta'(-s)\int_0^1H(t,u)\,dt\,ds\\
        &\qquad\qquad-\int_0^1\Big[\underbrace{\beta(-\infty)}_{=0}\cdot H(t,u(+\infty))-
                            \underbrace{\beta(+\infty)}_{=1}\cdot H(t,\underbrace{u(-\infty)}_{=x(t)})\Big]\,dt\\[2ex]
=&\;\underbrace{-\om(-u)+\int_0^1H(t,x(t))\,dt}_{=-\A_H(x,-u)}-\int_{-\infty}^{+\infty}\beta'(-s)\int_0^1H(t,u)\,dt\,ds\\[2ex]
\leq&\;-\A_H(x,-u)-\underbrace{\int_{-\infty}^{+\infty}\beta'(-s)\,ds}_{=1}\;\cdot\int_0^1\inf_MH(t,\cdot)\,dt\\[2ex]
=&\;-\A_H(x,-u)+\inf_MH
\eea
\end{proof}

\begin{Rmk}\label{rmk:energy_equality_for_Floer_strips}
If the cut-off function $\beta$ is removed in the above computation we obtain for $u\in\M_L(x,y;J,H)$ the standard equality
\beq
E(u)=\A_H(y,d_x\#u)-\A_H(x,d_x)\,,
\eeq
where $d_x\#u$ denotes the concatenation of the half-disk $d_x$ with the Floer strip $u$.
Furthermore, the same energy equality holds for maps $u\in\M^\chi(a,x)$ (see equation \eqref{def:moduli_space_for_chi})
which are defined on the Riemann surface $\Sigma^\chi$, since $\Sigma^\chi$ is expressed in global conformal coordinates.
Again the same holds for elements in $\M^\tau(x,a)$.
\end{Rmk}
\noindent
Finally, we give the uniform energy estimate for elements in the moduli space $\M^{\varphi\circ\rho}(q,p;a)$, which are defined in equations
\eqref{eqn:def_M(R,H;a)_start} -- \eqref{eqn:def_M(R,H;a)_end}.

\begin{Lemma}\label{lemma:energy_estimate_M(R,H;a)}
For $(R,\gamma_-,u,\gamma_+)\in\M^{\varphi\circ\rho}(q,p;a)$ holds
\beq
0\leq E(u)\leq \lambda a+||H||\,,
\eeq
where $\lambda$ is the constant from the monotonicity assumption (see definition \ref{def:minimal_Maslov_and_Chern_nb}).
\end{Lemma}

\begin{proof}
Since the energy is invariant under conformal changes, we can assume according to remark \ref{rmk:conformal_change} that
an element $u$ is map $u:\R\times[0,1]\pf M$ solving Floer's equation (with $R\geq0$):
$\partial_su+J(s,t,u)\big(\partial_tu-\tilde{\alpha}_R(s)\cdot X_H(t,u)\big)=0$, where $\tilde{\alpha}_R$ is a cut-off function such that
$\tilde{\alpha}_R(s)=1$ for $|s|\leq R$ and $\tilde{\alpha}_R(s)=0$ for $|s|\geq R+1$. Furthermore, we require
$\tilde{\alpha}_R'(s)\leq0$ for $s\geq0$ and $\tilde{\alpha}_R'(s)\geq0$ for $s\leq0$ and $\Mas(u)=a$.

With this we can redo the first lines in the above computations with $\tilde{\alpha}_R(s)$ replacing $\beta(s)$.
\bean
E(u)=&\;\int_{-\infty}^{+\infty}\int_0^1|\partial_s u|^2\,dt\,ds\\[2ex]
=&\;\om(u)+\int_{-\infty}^{+\infty}\tilde{\alpha}_R'(s)\int_0^1H(t,u)\,dt\,ds\\
        &\qquad\qquad-\int_0^1\Big[\underbrace{\tilde{\alpha}_R(+\infty)}_{=0}\cdot H(t,u(+\infty))-
                            \underbrace{\tilde{\alpha}_R(-\infty)}_{=0}\cdot H(t,u(-\infty))\Big]\,dt\\[2ex]
=&\;\underbrace{\;\om(u)\;}_{=\lambda a}+\int_{-\infty}^{+\infty}\int_0^1\tilde{\alpha}_R'(s)\cdot H(t,u)\,dt\,ds\\[2ex]
=&\;\lambda a+\int_{-\infty}^{0}\int_0^1\underbrace{\tilde{\alpha}_R'(s)}_{0\leq\;\bullet\;\leq1}\cdot H(t,u)\,dt\,ds
+\int_{0}^{+\infty}\int_0^1\underbrace{\tilde{\alpha}_R'(s)}_{-1\leq\;\bullet\;\leq0\;\;\;}\cdot H(t,u)\,dt\,ds\\[2ex]
\leq&\;\lambda a+\sup_MH-\inf_MH\;=\;\lambda a+||H||
\eea
\end{proof}
\noindent\hrulefill

%
\bibliographystyle{amsalpha}
\bibliography{../../../Bibtex/bibtex_paper_list}
\end{document}

%% file: perturbed_strip.pstex_t
\begin{picture}(0,0)%
\includegraphics{perturbed_strip.pstex}%
\end{picture}%
\setlength{\unitlength}{3947sp}%
\begingroup\makeatletter\ifx\SetFigFont\undefined%
\gdef\SetFigFont#1#2#3#4#5{%
  \reset@font\fontsize{#1}{#2pt}%
  \fontfamily{#3}\fontseries{#4}\fontshape{#5}%
  \selectfont}%
\fi\endgroup%
\begin{picture}(7224,2424)(1789,-4573)
\put(4951,-3436){\makebox(0,0)[lb]{\smash{{\SetFigFont{11}{13.2}{\rmdefault}{\mddefault}{\updefault}{\color[rgb]{0,0,0}$u$}%
}}}}
\put(2100,-4480){\makebox(0,0)[lb]{\smash{{\SetFigFont{11}{13.2}{\rmdefault}{\mddefault}{\updefault}{\color[rgb]{0,0,0}$L$}%
}}}}
\end{picture}%

%% file: PSS_strip_1.pstex_t
\begin{picture}(0,0)%
\includegraphics{PSS_strip_1.pstex}%
\end{picture}%
\setlength{\unitlength}{3947sp}%
\begingroup\makeatletter\ifx\SetFigFont\undefined%
\gdef\SetFigFont#1#2#3#4#5{%
  \reset@font\fontsize{#1}{#2pt}%
  \fontfamily{#3}\fontseries{#4}\fontshape{#5}%
  \selectfont}%
\fi\endgroup%
\begin{picture}(7224,2424)(1789,-4573)
\put(6001,-3286){\makebox(0,0)[lb]{\smash{{\SetFigFont{11}{13.2}{\rmdefault}{\mddefault}{\updefault}{\color[rgb]{0,0,0}$J=J_0$}%
}}}}
\put(2100,-4480){\makebox(0,0)[lb]{\smash{{\SetFigFont{11}{13.2}{\rmdefault}{\mddefault}{\updefault}{\color[rgb]{0,0,0}$L$}%
}}}}
\put(7881,-2756){\makebox(0,0)[lb]{\smash{{\SetFigFont{11}{13.2}{\rmdefault}{\mddefault}{\updefault}{\color[rgb]{0,0,0}$x$}%
}}}}
\put(2934,-3514){\makebox(0,0)[lb]{\smash{{\SetFigFont{11}{13.2}{\rmdefault}{\mddefault}{\updefault}{\color[rgb]{0,0,0}$\gamma$}%
}}}}
\put(3131,-4311){\makebox(0,0)[lb]{\smash{{\SetFigFont{11}{13.2}{\rmdefault}{\mddefault}{\updefault}{\color[rgb]{0,0,0}$q$}%
}}}}
\put(4801,-3361){\makebox(0,0)[lb]{\smash{{\SetFigFont{11}{13.2}{\rmdefault}{\mddefault}{\updefault}{\color[rgb]{0,0,0}$u$}%
}}}}
\end{picture}%

%% file: PSS_strip_1_breaking.pstex_t
\begin{picture}(0,0)%
\includegraphics{PSS_strip_1_breaking.pstex}%
\end{picture}%
\setlength{\unitlength}{3947sp}%
\begingroup\makeatletter\ifx\SetFigFont\undefined%
\gdef\SetFigFont#1#2#3#4#5{%
  \reset@font\fontsize{#1}{#2pt}%
  \fontfamily{#3}\fontseries{#4}\fontshape{#5}%
  \selectfont}%
\fi\endgroup%
\begin{picture}(5724,1824)(2239,-3973)
\put(4161,-3146){\makebox(0,0)[lb]{\smash{{\SetFigFont{11}{13.2}{\rmdefault}{\mddefault}{\updefault}{\color[rgb]{0,0,0}$p$}%
}}}}
\put(3292,-3263){\makebox(0,0)[lb]{\smash{{\SetFigFont{11}{13.2}{\rmdefault}{\mddefault}{\updefault}{\color[rgb]{0,0,0}$q$}%
}}}}
\put(5742,-3692){\makebox(0,0)[lb]{\smash{{\SetFigFont{11}{13.2}{\rmdefault}{\mddefault}{\updefault}{\color[rgb]{0,0,0}$y$}%
}}}}
\put(6766,-3199){\makebox(0,0)[lb]{\smash{{\SetFigFont{11}{13.2}{\rmdefault}{\mddefault}{\updefault}{\color[rgb]{0,0,0}$x$}%
}}}}
\put(4212,-3589){\makebox(0,0)[lb]{\smash{{\SetFigFont{11}{13.2}{\rmdefault}{\mddefault}{\updefault}{\color[rgb]{0,0,0}$q$}%
}}}}
\put(2551,-3811){\makebox(0,0)[lb]{\smash{{\SetFigFont{11}{13.2}{\rmdefault}{\mddefault}{\updefault}{\color[rgb]{0,0,0}$L$}%
}}}}
\put(4876,-2741){\makebox(0,0)[lb]{\smash{{\SetFigFont{11}{13.2}{\rmdefault}{\mddefault}{\updefault}{\color[rgb]{0,0,0}$x$}%
}}}}
\end{picture}%

%% file: PSS_strip_1_bubbling.pstex_t
\begin{picture}(0,0)%
\includegraphics{PSS_strip_1_bubbling.pstex}%
\end{picture}%
\setlength{\unitlength}{3947sp}%
\begingroup\makeatletter\ifx\SetFigFont\undefined%
\gdef\SetFigFont#1#2#3#4#5{%
  \reset@font\fontsize{#1}{#2pt}%
  \fontfamily{#3}\fontseries{#4}\fontshape{#5}%
  \selectfont}%
\fi\endgroup%
\begin{picture}(4674,1374)(2539,-3523)
\put(2776,-3436){\makebox(0,0)[lb]{\smash{{\SetFigFont{11}{13.2}{\rmdefault}{\mddefault}{\updefault}{\color[rgb]{0,0,0}$L$}%
}}}}
\put(3346,-2936){\makebox(0,0)[lb]{\smash{{\SetFigFont{11}{13.2}{\rmdefault}{\mddefault}{\updefault}{\color[rgb]{0,0,0}$q$}%
}}}}
\put(3901,-3286){\makebox(0,0)[lb]{\smash{{\SetFigFont{11}{13.2}{\rmdefault}{\mddefault}{\updefault}{\color[rgb]{0,0,0}$\gamma$}%
}}}}
\put(5551,-2461){\makebox(0,0)[lb]{\smash{{\SetFigFont{11}{13.2}{\rmdefault}{\mddefault}{\updefault}{\color[rgb]{0,0,0}$u$}%
}}}}
\put(4679,-2552){\makebox(0,0)[lb]{\smash{{\SetFigFont{11}{13.2}{\rmdefault}{\mddefault}{\updefault}{\color[rgb]{0,0,0}$d$}%
}}}}
\put(5734,-3127){\makebox(0,0)[lb]{\smash{{\SetFigFont{11}{13.2}{\rmdefault}{\mddefault}{\updefault}{\color[rgb]{0,0,0}$x$}%
}}}}
\end{picture}%

%% file: PSS_strip_2.pstex_t
\begin{picture}(0,0)%
\includegraphics{PSS_strip_2.pstex}%
\end{picture}%
\setlength{\unitlength}{3947sp}%
\begingroup\makeatletter\ifx\SetFigFont\undefined%
\gdef\SetFigFont#1#2#3#4#5{%
  \reset@font\fontsize{#1}{#2pt}%
  \fontfamily{#3}\fontseries{#4}\fontshape{#5}%
  \selectfont}%
\fi\endgroup%
\begin{picture}(7224,2424)(1789,-4573)
\put(3301,-3211){\makebox(0,0)[lb]{\smash{{\SetFigFont{11}{13.2}{\rmdefault}{\mddefault}{\updefault}{\color[rgb]{0,0,0}$J=J_0$}%
}}}}
\put(4569,-3399){\makebox(0,0)[lb]{\smash{{\SetFigFont{11}{13.2}{\rmdefault}{\mddefault}{\updefault}{\color[rgb]{0,0,0}$u$}%
}}}}
\put(6856,-2717){\makebox(0,0)[lb]{\smash{{\SetFigFont{11}{13.2}{\rmdefault}{\mddefault}{\updefault}{\color[rgb]{0,0,0}$\gamma$}%
}}}}
\put(8154,-2687){\makebox(0,0)[lb]{\smash{{\SetFigFont{11}{13.2}{\rmdefault}{\mddefault}{\updefault}{\color[rgb]{0,0,0}$q$}%
}}}}
\put(2926,-3661){\makebox(0,0)[lb]{\smash{{\SetFigFont{11}{13.2}{\rmdefault}{\mddefault}{\updefault}{\color[rgb]{0,0,0}$x$}%
}}}}
\put(2100,-4480){\makebox(0,0)[lb]{\smash{{\SetFigFont{11}{13.2}{\rmdefault}{\mddefault}{\updefault}{\color[rgb]{0,0,0}$L$}%
}}}}
\end{picture}%

%% file: composition_easy.pstex_t
\begin{picture}(0,0)%
\includegraphics{composition_easy.pstex}%
\end{picture}%
\setlength{\unitlength}{3947sp}%
\begingroup\makeatletter\ifx\SetFigFont\undefined%
\gdef\SetFigFont#1#2#3#4#5{%
  \reset@font\fontsize{#1}{#2pt}%
  \fontfamily{#3}\fontseries{#4}\fontshape{#5}%
  \selectfont}%
\fi\endgroup%
\begin{picture}(7224,5274)(1789,-7423)
\put(2100,-7330){\makebox(0,0)[lb]{\smash{{\SetFigFont{11}{13.2}{\rmdefault}{\mddefault}{\updefault}{\color[rgb]{0,0,0}$L$}%
}}}}
\put(2100,-7330){\makebox(0,0)[lb]{\smash{{\SetFigFont{11}{13.2}{\rmdefault}{\mddefault}{\updefault}{\color[rgb]{0,0,0}$L$}%
}}}}
\put(3261,-5791){\makebox(0,0)[lb]{\smash{{\SetFigFont{11}{13.2}{\rmdefault}{\mddefault}{\updefault}{\color[rgb]{0,0,0}$x$}%
}}}}
\put(2100,-4480){\makebox(0,0)[lb]{\smash{{\SetFigFont{11}{13.2}{\rmdefault}{\mddefault}{\updefault}{\color[rgb]{0,0,0}$L$}%
}}}}
\put(2100,-4480){\makebox(0,0)[lb]{\smash{{\SetFigFont{11}{13.2}{\rmdefault}{\mddefault}{\updefault}{\color[rgb]{0,0,0}$L$}%
}}}}
\put(5476,-2741){\makebox(0,0)[lb]{\smash{{\SetFigFont{11}{13.2}{\rmdefault}{\mddefault}{\updefault}{\color[rgb]{0,0,0}$y$}%
}}}}
\put(4761,-3146){\makebox(0,0)[lb]{\smash{{\SetFigFont{11}{13.2}{\rmdefault}{\mddefault}{\updefault}{\color[rgb]{0,0,0}$q$}%
}}}}
\put(6776,-3631){\makebox(0,0)[lb]{\smash{{\SetFigFont{11}{13.2}{\rmdefault}{\mddefault}{\updefault}{\color[rgb]{0,0,0}$y$}%
}}}}
\put(3261,-2941){\makebox(0,0)[lb]{\smash{{\SetFigFont{11}{13.2}{\rmdefault}{\mddefault}{\updefault}{\color[rgb]{0,0,0}$x$}%
}}}}
\put(4606,-3851){\makebox(0,0)[lb]{\smash{{\SetFigFont{11}{13.2}{\rmdefault}{\mddefault}{\updefault}{\color[rgb]{0,0,0}$x$}%
}}}}
\put(5787,-3760){\makebox(0,0)[b]{\smash{{\SetFigFont{11}{13.2}{\familydefault}{\mddefault}{\updefault}{\color[rgb]{0,0,0}2}%
}}}}
\put(4005,-5552){\makebox(0,0)[b]{\smash{{\SetFigFont{11}{13.2}{\familydefault}{\mddefault}{\updefault}{\color[rgb]{0,0,0}3}%
}}}}
\put(4576,-6041){\makebox(0,0)[lb]{\smash{{\SetFigFont{11}{13.2}{\rmdefault}{\mddefault}{\updefault}{\color[rgb]{0,0,0}$y$}%
}}}}
\put(5236,-4936){\makebox(0,0)[lb]{\smash{{\SetFigFont{11}{13.2}{\rmdefault}{\mddefault}{\updefault}{\color[rgb]{0,0,0}shrink length to zero}%
}}}}
\put(3001,-3811){\makebox(0,0)[lb]{\smash{{\SetFigFont{11}{13.2}{\rmdefault}{\mddefault}{\updefault}{\color[rgb]{0,0,0}glue at $q$}%
}}}}
\put(5741,-4166){\makebox(0,0)[lb]{\smash{{\SetFigFont{11}{13.2}{\rmdefault}{\mddefault}{\updefault}{\color[rgb]{0,0,0}finite length}%
}}}}
\put(5026,-6286){\makebox(0,0)[lb]{\smash{{\SetFigFont{11}{13.2}{\rmdefault}{\mddefault}{\updefault}{\color[rgb]{0,0,0}glue}%
}}}}
\put(4283,-2657){\makebox(0,0)[b]{\smash{{\SetFigFont{11}{13.2}{\familydefault}{\mddefault}{\updefault}{\color[rgb]{0,0,0}1}%
}}}}
\put(6094,-6296){\makebox(0,0)[lb]{\smash{{\SetFigFont{11}{13.2}{\rmdefault}{\mddefault}{\updefault}{\color[rgb]{0,0,0}$x$}%
}}}}
\put(7148,-6439){\makebox(0,0)[lb]{\smash{{\SetFigFont{11}{13.2}{\rmdefault}{\mddefault}{\updefault}{\color[rgb]{0,0,0}$y$}%
}}}}
\put(6833,-5957){\makebox(0,0)[b]{\smash{{\SetFigFont{11}{13.2}{\familydefault}{\mddefault}{\updefault}{\color[rgb]{0,0,0}4}%
}}}}
\end{picture}%

%% file: composition_hard.pstex_t
\begin{picture}(0,0)%
\includegraphics{composition_hard.pstex}%
\end{picture}%
\setlength{\unitlength}{3947sp}%
\begingroup\makeatletter\ifx\SetFigFont\undefined%
\gdef\SetFigFont#1#2#3#4#5{%
  \reset@font\fontsize{#1}{#2pt}%
  \fontfamily{#3}\fontseries{#4}\fontshape{#5}%
  \selectfont}%
\fi\endgroup%
\begin{picture}(4824,1674)(2389,-3823)
\put(5776,-2986){\makebox(0,0)[lb]{\smash{{\SetFigFont{11}{13.2}{\rmdefault}{\mddefault}{\updefault}{\color[rgb]{0,0,0}$p$}%
}}}}
\put(2622,-3724){\makebox(0,0)[lb]{\smash{{\SetFigFont{11}{13.2}{\rmdefault}{\mddefault}{\updefault}{\color[rgb]{0,0,0}$L$}%
}}}}
\put(4461,-2595){\makebox(0,0)[lb]{\smash{{\SetFigFont{11}{13.2}{\rmdefault}{\mddefault}{\updefault}{\color[rgb]{0,0,0}$x$}%
}}}}
\put(3532,-3574){\makebox(0,0)[lb]{\smash{{\SetFigFont{11}{13.2}{\rmdefault}{\mddefault}{\updefault}{\color[rgb]{0,0,0}$q$}%
}}}}
\end{picture}%

%% file: omega.pstex_t
\begin{picture}(0,0)%
\includegraphics{omega.pstex}%
\end{picture}%
\setlength{\unitlength}{3947sp}%
\begingroup\makeatletter\ifx\SetFigFont\undefined%
\gdef\SetFigFont#1#2#3#4#5{%
  \reset@font\fontsize{#1}{#2pt}%
  \fontfamily{#3}\fontseries{#4}\fontshape{#5}%
  \selectfont}%
\fi\endgroup%
\begin{picture}(2559,1177)(451,-926)
\put(1726,-586){\makebox(0,0)[lb]{\smash{{\SetFigFont{11}{13.2}{\rmdefault}{\mddefault}{\updefault}{\color[rgb]{0,0,0}$I$}%
}}}}
\put(1426,-886){\makebox(0,0)[lb]{\smash{{\SetFigFont{11}{13.2}{\rmdefault}{\mddefault}{\updefault}{\color[rgb]{0,0,0}$-1$}%
}}}}
\put(1126,-886){\makebox(0,0)[lb]{\smash{{\SetFigFont{11}{13.2}{\rmdefault}{\mddefault}{\updefault}{\color[rgb]{0,0,0}$-2$}%
}}}}
\put(2026,-886){\makebox(0,0)[lb]{\smash{{\SetFigFont{11}{13.2}{\rmdefault}{\mddefault}{\updefault}{\color[rgb]{0,0,0}$1$}%
}}}}
\put(2326,-886){\makebox(0,0)[lb]{\smash{{\SetFigFont{11}{13.2}{\rmdefault}{\mddefault}{\updefault}{\color[rgb]{0,0,0}$2$}%
}}}}
\put(451,-361){\makebox(0,0)[lb]{\smash{{\SetFigFont{11}{13.2}{\rmdefault}{\mddefault}{\updefault}{\color[rgb]{0,0,0}$\D_-$}%
}}}}
\put(3001,-361){\makebox(0,0)[lb]{\smash{{\SetFigFont{11}{13.2}{\rmdefault}{\mddefault}{\updefault}{\color[rgb]{0,0,0}$\D_+$}%
}}}}
\put(826,-886){\makebox(0,0)[lb]{\smash{{\SetFigFont{11}{13.2}{\rmdefault}{\mddefault}{\updefault}{\color[rgb]{0,0,0}$-3$}%
}}}}
\put(2626,-886){\makebox(0,0)[lb]{\smash{{\SetFigFont{11}{13.2}{\rmdefault}{\mddefault}{\updefault}{\color[rgb]{0,0,0}$3$}%
}}}}
\put(976,-136){\makebox(0,0)[lb]{\smash{{\SetFigFont{11}{13.2}{\rmdefault}{\mddefault}{\updefault}{\color[rgb]{0,0,0}$\alpha$}%
}}}}
\put(2476,-136){\makebox(0,0)[lb]{\smash{{\SetFigFont{11}{13.2}{\rmdefault}{\mddefault}{\updefault}{\color[rgb]{0,0,0}$\alpha$}%
}}}}
\end{picture}%

%% file: rescale.pstex_t
\begin{picture}(0,0)%
\includegraphics{rescale.pstex}%
\end{picture}%
\setlength{\unitlength}{3947sp}%
\begingroup\makeatletter\ifx\SetFigFont\undefined%
\gdef\SetFigFont#1#2#3#4#5{%
  \reset@font\fontsize{#1}{#2pt}%
  \fontfamily{#3}\fontseries{#4}\fontshape{#5}%
  \selectfont}%
\fi\endgroup%
\begin{picture}(6612,3657)(601,-2806)
\put(4009,-2221){\makebox(0,0)[lb]{\smash{{\SetFigFont{11}{13.2}{\rmdefault}{\mddefault}{\updefault}{\color[rgb]{0,0,0}$J=J_0$}%
}}}}
\put(4351,-1111){\makebox(0,0)[lb]{\smash{{\SetFigFont{11}{13.2}{\rmdefault}{\mddefault}{\updefault}{\color[rgb]{0,0,0}$\phi_R$}%
}}}}
\put(6601,-661){\makebox(0,0)[lb]{\smash{{\SetFigFont{11}{13.2}{\rmdefault}{\mddefault}{\updefault}{\color[rgb]{0,0,0}$R+2$}%
}}}}
\put(601,389){\makebox(0,0)[lb]{\smash{{\SetFigFont{11}{13.2}{\rmdefault}{\mddefault}{\updefault}{\color[rgb]{0,0,0}$0$}%
}}}}
\put(601,614){\makebox(0,0)[lb]{\smash{{\SetFigFont{11}{13.2}{\rmdefault}{\mddefault}{\updefault}{\color[rgb]{0,0,0}$1$}%
}}}}
\put(6001,-661){\makebox(0,0)[lb]{\smash{{\SetFigFont{11}{13.2}{\rmdefault}{\mddefault}{\updefault}{\color[rgb]{0,0,0}$R+1$}%
}}}}
\put(2701,-661){\makebox(0,0)[lb]{\smash{{\SetFigFont{11}{13.2}{\rmdefault}{\mddefault}{\updefault}{\color[rgb]{0,0,0}$-R$}%
}}}}
\put(5626,-661){\makebox(0,0)[lb]{\smash{{\SetFigFont{11}{13.2}{\rmdefault}{\mddefault}{\updefault}{\color[rgb]{0,0,0}$R$}%
}}}}
\put(2551,-2236){\makebox(0,0)[lb]{\smash{{\SetFigFont{11}{13.2}{\rmdefault}{\mddefault}{\updefault}{\color[rgb]{0,0,0}$\alpha$}%
}}}}
\put(5701,-2236){\makebox(0,0)[lb]{\smash{{\SetFigFont{11}{13.2}{\rmdefault}{\mddefault}{\updefault}{\color[rgb]{0,0,0}$\alpha$}%
}}}}
\put(1276,-661){\makebox(0,0)[lb]{\smash{{\SetFigFont{11}{13.2}{\rmdefault}{\mddefault}{\updefault}{\color[rgb]{0,0,0}$-R-2$}%
}}}}
\put(2026,-661){\makebox(0,0)[lb]{\smash{{\SetFigFont{11}{13.2}{\rmdefault}{\mddefault}{\updefault}{\color[rgb]{0,0,0}$-R-1$}%
}}}}
\put(4126,-2761){\makebox(0,0)[lb]{\smash{{\SetFigFont{11}{13.2}{\rmdefault}{\mddefault}{\updefault}{\color[rgb]{0,0,0}$\Omega$}%
}}}}
\put(4126,389){\makebox(0,0)[lb]{\smash{{\SetFigFont{11}{13.2}{\rmdefault}{\mddefault}{\updefault}{\color[rgb]{0,0,0}$\alpha_R$}%
}}}}
\put(1951,-136){\makebox(0,0)[lb]{\smash{{\SetFigFont{11}{13.2}{\rmdefault}{\mddefault}{\updefault}{\color[rgb]{0,0,0}$\alpha_R$}%
}}}}
\put(6301,-136){\makebox(0,0)[lb]{\smash{{\SetFigFont{11}{13.2}{\rmdefault}{\mddefault}{\updefault}{\color[rgb]{0,0,0}$\alpha_R$}%
}}}}
\put(5701,-1411){\makebox(0,0)[lb]{\smash{{\SetFigFont{11}{13.2}{\rmdefault}{\mddefault}{\updefault}{\color[rgb]{0,0,0}$V$}%
}}}}
\put(2551,-1411){\makebox(0,0)[lb]{\smash{{\SetFigFont{11}{13.2}{\rmdefault}{\mddefault}{\updefault}{\color[rgb]{0,0,0}$V$}%
}}}}
\end{picture}%

%% file: bad_bubbling.pstex_t
\begin{picture}(0,0)%
\includegraphics{bad_bubbling.pstex}%
\end{picture}%
\setlength{\unitlength}{3947sp}%
\begingroup\makeatletter\ifx\SetFigFont\undefined%
\gdef\SetFigFont#1#2#3#4#5{%
  \reset@font\fontsize{#1}{#2pt}%
  \fontfamily{#3}\fontseries{#4}\fontshape{#5}%
  \selectfont}%
\fi\endgroup%
\begin{picture}(7224,2424)(1789,-4573)
\put(2100,-4480){\makebox(0,0)[lb]{\smash{{\SetFigFont{11}{13.2}{\rmdefault}{\mddefault}{\updefault}{\color[rgb]{0,0,0}$L$}%
}}}}
\put(5626,-2836){\makebox(0,0)[lb]{\smash{{\SetFigFont{11}{13.2}{\rmdefault}{\mddefault}{\updefault}{\color[rgb]{0,0,0}$U$}%
}}}}
\put(7351,-3136){\makebox(0,0)[lb]{\smash{{\SetFigFont{11}{13.2}{\rmdefault}{\mddefault}{\updefault}{\color[rgb]{0,0,0}$p$}%
}}}}
\put(2896,-3386){\makebox(0,0)[lb]{\smash{{\SetFigFont{11}{13.2}{\rmdefault}{\mddefault}{\updefault}{\color[rgb]{0,0,0}$q$}%
}}}}
\put(3451,-3736){\makebox(0,0)[lb]{\smash{{\SetFigFont{11}{13.2}{\rmdefault}{\mddefault}{\updefault}{\color[rgb]{0,0,0}$\gamma_-$}%
}}}}
\put(5303,-4006){\makebox(0,0)[lb]{\smash{{\SetFigFont{11}{13.2}{\rmdefault}{\mddefault}{\updefault}{\color[rgb]{0,0,0}$R_\infty$}%
}}}}
\put(7126,-3661){\makebox(0,0)[lb]{\smash{{\SetFigFont{11}{13.2}{\rmdefault}{\mddefault}{\updefault}{\color[rgb]{0,0,0}$\gamma_+$}%
}}}}
\put(4304,-3002){\makebox(0,0)[lb]{\smash{{\SetFigFont{11}{13.2}{\rmdefault}{\mddefault}{\updefault}{\color[rgb]{0,0,0}$d$}%
}}}}
\end{picture}%

%% file: chi_pic.pstex_t
\begin{picture}(0,0)%
\includegraphics{chi_pic.pstex}%
\end{picture}%
\setlength{\unitlength}{3947sp}%
\begingroup\makeatletter\ifx\SetFigFont\undefined%
\gdef\SetFigFont#1#2#3#4#5{%
  \reset@font\fontsize{#1}{#2pt}%
  \fontfamily{#3}\fontseries{#4}\fontshape{#5}%
  \selectfont}%
\fi\endgroup%
\begin{picture}(6024,4066)(2989,-4273)
\put(7881,-2756){\makebox(0,0)[lb]{\smash{{\SetFigFont{11}{13.2}{\rmdefault}{\mddefault}{\updefault}{\color[rgb]{0,0,0}$x$}%
}}}}
\put(5211,-1291){\makebox(0,0)[lb]{\smash{{\SetFigFont{11}{13.2}{\rmdefault}{\mddefault}{\updefault}{\color[rgb]{0,0,0}$u$}%
}}}}
\put(3451,-4111){\makebox(0,0)[lb]{\smash{{\SetFigFont{11}{13.2}{\rmdefault}{\mddefault}{\updefault}{\color[rgb]{0,0,0}$L$}%
}}}}
\put(4351,-361){\makebox(0,0)[lb]{\smash{{\SetFigFont{11}{13.2}{\rmdefault}{\mddefault}{\updefault}{\color[rgb]{0,0,0}$a$}%
}}}}
\end{picture}%

%% file: diagram_commutes_1.pstex_t
\begin{picture}(0,0)%
\includegraphics{diagram_commutes_1.pstex}%
\end{picture}%
\setlength{\unitlength}{3947sp}%
\begingroup\makeatletter\ifx\SetFigFont\undefined%
\gdef\SetFigFont#1#2#3#4#5{%
  \reset@font\fontsize{#1}{#2pt}%
  \fontfamily{#3}\fontseries{#4}\fontshape{#5}%
  \selectfont}%
\fi\endgroup%
\begin{picture}(7224,3894)(1789,-4573)
\put(7576,-2986){\makebox(0,0)[lb]{\smash{{\SetFigFont{11}{13.2}{\rmdefault}{\mddefault}{\updefault}{\color[rgb]{0,0,0}$\gamma_L$}%
}}}}
\put(3076,-1711){\makebox(0,0)[lb]{\smash{{\SetFigFont{11}{13.2}{\rmdefault}{\mddefault}{\updefault}{\color[rgb]{0,0,0}$\PSS^{-1}$}%
}}}}
\put(4651,-3961){\makebox(0,0)[lb]{\smash{{\SetFigFont{11}{13.2}{\rmdefault}{\mddefault}{\updefault}{\color[rgb]{0,0,0}$\rho$}%
}}}}
\put(5776,-1336){\makebox(0,0)[lb]{\smash{{\SetFigFont{11}{13.2}{\rmdefault}{\mddefault}{\updefault}{\color[rgb]{0,0,0}$q_i$}%
}}}}
\put(5476,-1111){\makebox(0,0)[lb]{\smash{{\SetFigFont{11}{13.2}{\rmdefault}{\mddefault}{\updefault}{\color[rgb]{0,0,0}$\gamma_M$}%
}}}}
\put(5251,-811){\makebox(0,0)[lb]{\smash{{\SetFigFont{11}{13.2}{\rmdefault}{\mddefault}{\updefault}{\color[rgb]{0,0,0}$q$}%
}}}}
\put(7051,-2911){\makebox(0,0)[lb]{\smash{{\SetFigFont{11}{13.2}{\rmdefault}{\mddefault}{\updefault}{\color[rgb]{0,0,0}$q_b$}%
}}}}
\put(2100,-4480){\makebox(0,0)[lb]{\smash{{\SetFigFont{11}{13.2}{\rmdefault}{\mddefault}{\updefault}{\color[rgb]{0,0,0}$L$}%
}}}}
\put(2100,-4480){\makebox(0,0)[lb]{\smash{{\SetFigFont{11}{13.2}{\rmdefault}{\mddefault}{\updefault}{\color[rgb]{0,0,0}$L$}%
}}}}
\put(5401,-4186){\makebox(0,0)[lb]{\smash{{\SetFigFont{11}{13.2}{\rmdefault}{\mddefault}{\updefault}{\color[rgb]{0,0,0}$p$}%
}}}}
\put(7576,-3511){\makebox(0,0)[lb]{\smash{{\SetFigFont{11}{13.2}{\rmdefault}{\mddefault}{\updefault}{\color[rgb]{0,0,0}$p$}%
}}}}
\put(2251,-1411){\makebox(0,0)[lb]{\smash{{\SetFigFont{11}{13.2}{\rmdefault}{\mddefault}{\updefault}{\color[rgb]{0,0,0}$q$}%
}}}}
\put(5626,-2386){\makebox(0,0)[lb]{\smash{{\SetFigFont{11}{13.2}{\rmdefault}{\mddefault}{\updefault}{\color[rgb]{0,0,0}$u$}%
}}}}
\put(4276,-2536){\makebox(0,0)[lb]{\smash{{\SetFigFont{11}{13.2}{\rmdefault}{\mddefault}{\updefault}{\color[rgb]{0,0,0}gluing gives:}%
}}}}
\put(2626,-3736){\makebox(0,0)[lb]{\smash{{\SetFigFont{11}{13.2}{\rmdefault}{\mddefault}{\updefault}{\color[rgb]{0,0,0}$\chi$}%
}}}}
\end{picture}%

%% file: comparison.pstex_t
\begin{picture}(0,0)%
\includegraphics{comparison.pstex}%
\end{picture}%
\setlength{\unitlength}{3947sp}%
\begingroup\makeatletter\ifx\SetFigFont\undefined%
\gdef\SetFigFont#1#2#3#4#5{%
  \reset@font\fontsize{#1}{#2pt}%
  \fontfamily{#3}\fontseries{#4}\fontshape{#5}%
  \selectfont}%
\fi\endgroup%
\begin{picture}(3534,1880)(226,-1146)
\put(1981,-97){\makebox(0,0)[lb]{\smash{{\SetFigFont{10}{12.0}{\rmdefault}{\mddefault}{\updefault}{\color[rgb]{0,0,0}$S$}%
}}}}
\put(1051, 14){\makebox(0,0)[lb]{\smash{{\SetFigFont{10}{12.0}{\rmdefault}{\mddefault}{\updefault}{\color[rgb]{0,0,0}$\alpha$}%
}}}}
\put(2776,-97){\makebox(0,0)[lb]{\smash{{\SetFigFont{10}{12.0}{\rmdefault}{\mddefault}{\updefault}{\color[rgb]{0,0,0}$\alpha=0$}%
}}}}
\put(226,614){\makebox(0,0)[lb]{\smash{{\SetFigFont{10}{12.0}{\rmdefault}{\mddefault}{\updefault}{\color[rgb]{0,0,0}$\alpha=0$}%
}}}}
\put(2533,-97){\makebox(0,0)[lb]{\smash{{\SetFigFont{10}{12.0}{\rmdefault}{\mddefault}{\updefault}{\color[rgb]{0,0,0}$\alpha$}%
}}}}
\put(816,-1106){\makebox(0,0)[lb]{\smash{{\SetFigFont{10}{12.0}{\rmdefault}{\mddefault}{\updefault}{\color[rgb]{0,0,0}$A$}%
}}}}
\put(901,-116){\makebox(0,0)[lb]{\smash{{\SetFigFont{10}{12.0}{\rmdefault}{\mddefault}{\updefault}{\color[rgb]{0,0,0}$p_i$}%
}}}}
\put(3376,-136){\makebox(0,0)[lb]{\smash{{\SetFigFont{10}{12.0}{\rmdefault}{\mddefault}{\updefault}{\color[rgb]{0,0,0}$p_b$}%
}}}}
\end{picture}%

%% file: diagram_commutes_2.pstex_t
\begin{picture}(0,0)%
\includegraphics{diagram_commutes_2.pstex}%
\end{picture}%
\setlength{\unitlength}{3947sp}%
\begingroup\makeatletter\ifx\SetFigFont\undefined%
\gdef\SetFigFont#1#2#3#4#5{%
  \reset@font\fontsize{#1}{#2pt}%
  \fontfamily{#3}\fontseries{#4}\fontshape{#5}%
  \selectfont}%
\fi\endgroup%
\begin{picture}(7224,3095)(1789,-4573)
\put(4051,-3286){\makebox(0,0)[lb]{\smash{{\SetFigFont{11}{13.2}{\rmdefault}{\mddefault}{\updefault}{\color[rgb]{0,0,0}gluing}%
}}}}
\put(2100,-4480){\makebox(0,0)[lb]{\smash{{\SetFigFont{11}{13.2}{\rmdefault}{\mddefault}{\updefault}{\color[rgb]{0,0,0}$L$}%
}}}}
\put(2100,-4480){\makebox(0,0)[lb]{\smash{{\SetFigFont{11}{13.2}{\rmdefault}{\mddefault}{\updefault}{\color[rgb]{0,0,0}$L$}%
}}}}
\put(2626,-3736){\makebox(0,0)[lb]{\smash{{\SetFigFont{11}{13.2}{\rmdefault}{\mddefault}{\updefault}{\color[rgb]{0,0,0}$\chi$}%
}}}}
\put(4651,-3961){\makebox(0,0)[lb]{\smash{{\SetFigFont{11}{13.2}{\rmdefault}{\mddefault}{\updefault}{\color[rgb]{0,0,0}$\rho$}%
}}}}
\end{picture}%

%% file: Lagrangian_PSS_and_comparison.bbl
\providecommand{\bysame}{\leavevmode\hbox to3em{\hrulefill}\thinspace}
\providecommand{\MR}{\relax\ifhmode\unskip\space\fi MR }
\providecommand{\MRhref}[2]{%
  \href{http://www.ams.org/mathscinet-getitem?mr=#1}{#2}
}
\providecommand{\href}[2]{#2}
\begin{thebibliography}{FOOO}

\bibitem[Alb05]{Albers_On_the_extrinsic_topology_of_Lagrangian_submanifolds}
P.~Albers, \emph{On the extrinsic topology of {L}agrangian submanifolds},
  Int.~Math.~Res.~Not. (2005), no.~38, 2341--2371.

\bibitem[AS06]{Abbo_Schwarz_Notes_on_Floer_homology_and_loop_space_homology}
A.~Abbondandolo and M.~Schwarz, \emph{Notes on {F}loer homology and loop space
  homology}, Morse theoretic methods in nonlinear analysis and in symplectic
  topology, NATO Sci. Ser. II Math. Phys. Chem., vol. 217, Springer, Dordrecht,
  2006, pp.~75--108.

\bibitem[BC06]{Barraud_Cornea_Homotopical_dynamics_in_symplectic_topology}
J.-F. Barraud and O.~Cornea, \emph{{Homotopical dynamics in symplectic
  topology}}, {Morse theoretic methods in nonlinear analysis and in symplectic
  topology}, {Proceedings of the NATO Advanced Study Institute, Montréal,
  Canada, July 2004}, Springer, Dordrecht, 2006.

\bibitem[BC07]{Biran_Cornea_Quantum_Structures}
P.~Biran and O.~Cornea, \emph{{Quantum structures for Lagrangian
  submanifolds}}, preprint 2007, arXiv:0708.4221.

\bibitem[CL06]{Cornea_Lalonde_Cluster_Homology:_an_overview_of_the_constructio%
n_and_results}
O.~Cornea and F.~Lalonde, \emph{Cluster homology: an overview of the
  construction and results}, Electron. Res. Announc. Amer. Math. Soc.
  \textbf{12} (2006), 1--12 (electronic).

\bibitem[FHS95]{Floer_Hofer_Salamon_Transversality_in_elliptic_Morse_theory_fo%
r_the_symplectic_action}
A.~Floer, H.~Hofer, and D.~Salamon, \emph{Transversality in elliptic {M}orse
  theory for the symplectic action}, Duke Math. J. \textbf{80} (1995), no.~1,
  251--292.

\bibitem[Flo88a]{Floer_Morse_theory_for_Lagrangian_intersections}
A.~Floer, \emph{Morse theory for {L}agrangian intersections}, J. Differential
  Geom. \textbf{28} (1988), no.~3, 513--547.

\bibitem[Flo88b]{Floer_A_relative_Morse_index_for_the_symplectic_action}
\bysame, \emph{A relative {M}orse index for the symplectic action}, Comm. Pure
  Appl. Math. \textbf{41} (1988), no.~4, 393--407.

\bibitem[Flo89]{Floer_symplectic_fixed_points_and_holomorphic_spheres}
\bysame, \emph{Symplectic fixed points and holomorphic spheres}, Comm. Math.
  Phys. \textbf{120} (1989), no.~4, 575--611.

\bibitem[FOOO]{FOOO}
K.~Fukaya, Y.-G. Oh, H.~Ohta, and K.~Ono, \emph{{Lagrangian intersection Floer
  homology - anomaly and obstruction}}, Kyoto University preprint, 2000.

\bibitem[HA]{Hofer_Abbas_book_preprint}
H.~Hofer and C.~Abbas, \emph{preprint}, see
  http://www.math.nyu.edu/$\sim$hofer/lecture.html.

\bibitem[HS95]{Hofer_Salamon_Floer_homology_and_Novikov_rings}
H.~Hofer and D.~Salamon, \emph{Floer homology and {N}ovikov rings}, The Floer
  memorial volume, Progr. Math., vol. 133, Birkh\"auser, Basel, 1995,
  pp.~483--524.

\bibitem[KM05]{Katic_Milinkovic_PSS_ismomorphism_for_Lagrangian_intersections}
J.~Kati{\'c} and D.~Milinkovi{\'c}, \emph{Piunikhin-{S}alamon-{S}chwarz
  isomorphisms for {L}agrangian intersections}, Differential Geom. Appl.
  \textbf{22} (2005), no.~2, 215--227.

\bibitem[MS04]{McDuff_Salamon_J_holomorphic_curves_and_symplectic_topology}
D.~McDuff and D.~Salamon, \emph{{$J$}-holomorphic curves and symplectic
  topology}, American Mathematical Society Colloquium Publications, vol.~52,
  American Mathematical Society, Providence, RI, 2004.

\bibitem[Oh93]{Oh_Floer_cohomology_of_Lagrangian_intersections_and_pseudo-holo%
morphic_disks_I}
Y.-G. Oh, \emph{Floer cohomology of {L}agrangian intersections and
  pseudo-holomorphic disks {I}}, Comm. Pure Appl. Math. \textbf{46} (1993),
  no.~7, 949--993.

\bibitem[Oh95]{Oh_Addendum_Floer_cohomology_of_Lagrangian_intersections_and_ps%
eudo-holomorphic_disks_I}
\bysame, \emph{Addendum to: ``{F}loer cohomology of {L}agrangian intersections
  and pseudo-holomorphic disks {I}.''}, Comm. Pure Appl. Math. \textbf{48}
  (1995), no.~11, 1299--1302.

\bibitem[Oh96]{Oh_Floer_cohomology_spectral_sequences_and_the_Maslov_class}
\bysame, \emph{Floer cohomology, spectral sequences, and the {M}aslov class of
  {L}agrangian embeddings}, Internat. Math. Res. Notices (1996), no.~7,
  305--346.

\bibitem[Pol91]{Polterovich_Monotone_Lagrange_submanifolds_of_linear_spaces}
L.~Polterovich, \emph{Monotone {L}agrange submanifolds of linear spaces and the
  {M}aslov class in cotangent bundles}, Math. Z. \textbf{207} (1991), no.~2,
  217--222.

\bibitem[PSS96]{Piunikhin_Salamon_Schwarz_Symplectic_Floer_Donaldson_theory_an%
d_quantum_cohomology}
S.~Piunikhin, D.~Salamon, and M.~Schwarz, \emph{Symplectic {F}loer-{D}onaldson
  theory and quantum cohomology}, Contact and symplectic geometry (Cambridge,
  1994), Publ. Newton Inst., vol.~8, Cambridge Univ. Press, Cambridge, 1996,
  pp.~171--200.

\bibitem[Sal99]{Salamon_lectures_on_floer_homology}
D.~Salamon, \emph{Lectures on {F}loer homology}, Symplectic geometry and
  topology (Park City, UT, 1997), IAS/Park City Math. Ser., vol.~7, Amer. Math.
  Soc., Providence, RI, 1999, pp.~143--229.

\bibitem[Sch93]{Schwarz_Matthias_Morse_homology}
M.~Schwarz, \emph{Morse homology}, Progress in Mathematics, vol. 111,
  Birkh\"auser Verlag, Basel, 1993.

\bibitem[Sch95]{Schwarz_Matthias_PhD}
\bysame, \emph{Cohomology operations from ${S}^1$-cobordisms in {F}loer
  homology}, Ph.{D}.-thesis, Swiss Federal Inst. of Techn. Zurich, Diss. {ETH}
  No. 11182, 1995.

\bibitem[Sei00]{Seidel_Graded_Lagrangian_submanifolds}
P.~Seidel, \emph{Graded {L}agrangian submanifolds}, Bull. Soc. Math. France
  \textbf{128} (2000), no.~1, 103--149.

\bibitem[SZ92]{Salamon_Zehnder_Morse_theory_for_periodic_solutions_of_Hamilton%
ian_systems_and_the_Maslov_index}
D.~Salamon and E.~Zehnder, \emph{Morse theory for periodic solutions of
  {H}amiltonian systems and the {M}aslov index}, Comm. Pure Appl. Math.
  \textbf{45} (1992), no.~10, 1303--1360.

\bibitem[Vit87]{Viterbo_Maslov_index_for_Lagrangian_Floer_homology}
C.~Viterbo, \emph{Intersection de sous-vari\'et\'es lagrangiennes,
  fonctionnelles d'action et indice des syst\`emes hamiltoniens}, Bull. Soc.
  Math. France \textbf{115} (1987), no.~3, 361--390.

\end{thebibliography}
